\documentclass[reqno,twoside,11pt,english]{amsart}
\usepackage{amsmath,amsfonts,amssymb,amsthm,epsfig}
\newcommand{\RNum}[1]{\uppercase\expandafter{\romannumeral #1\relax}}
\usepackage{comment}

\usepackage{color}
\usepackage[final,colorlinks,linkcolor=blue,anchorcolor=red,citecolor=blue]{hyperref}
\usepackage{graphicx}
\usepackage{titletoc,epsf}
\usepackage{amsmath,amsfonts,latexsym,amsthm,amsxtra,amssymb,bbm}
\allowdisplaybreaks
\usepackage{graphicx,float}

\allowdisplaybreaks
\let\f=\frac

\let\Om=\Omega

\let\th=T
\let\pa=\partial

\def\R{\mathbf R}

\DeclareMathOperator{\CGL}{CGL}
\DeclareMathOperator{\NLS}{NLS}
\DeclareMathOperator{\diam}{diam}
\DeclareMathOperator{\dist}{dist}
\DeclareMathOperator{\supp}{supp}

\newtheorem{theorem}{Theorem}[section]
\newtheorem{lemma}[theorem]{Lemma}
\newtheorem{proposition}[theorem]{Proposition}
\newtheorem{corollary}[theorem]{Corollary}
\newtheorem{definition}[theorem]{Definition}

\newtheorem{thmA}{Theorem}[section]

\numberwithin{equation}{section}

\numberwithin{equation}{section}

\topmargin - .23 in

\begin{document}

\title[Global well-posedness of CGL in exterior domains]{Global well-posedness of the energy-critical complex Ginzburg-Landau equation in exterior domains}

\author{Xing Cheng, Chang-Yu Guo, Jiqiang Zheng and Yunrui Zheng}

\address[X. Cheng]{School of mathematics, Hohai University, Nanjing 210098, Jiangsu, P. R. China}
\email{{\tt chengx@hhu.edu.cn}}

\address[C.-Y. Guo]{Research Center for Mathematics and Interdisciplinary Sciences, Shandong University, 266237, Qingdao, P. R. China, and Department of Physics and Mathematics, University of Eastern Finland, 80101, Joensuu, Finland}
\email{{\tt changyu.guo@sdu.edu.cn}}

\address[J. Zheng]{Institute of Applied Physics and Computational Mathematics, Beijing, 100088 and National Key Laboratory of Computational Physics, Beijing 100088, China}
\email{{\tt  zheng\_jiqiang@iapcm.ac.cn}}

\address[Y. Zheng]{ School of Mathematics, Shandong University, Shandong 250100, Jinan, P. R. China.}
\email{{\tt yunrui\_zheng@sdu.edu.cn}}

\date{\today}
\thanks{X. Cheng has been partially supported by the NSF of Jiangsu Province (Grant No.~BK20221497). C.-Y. Guo is supported by the Young Scientist Program of the Ministry of Science and Technology of China (No.~2021YFA1002200), the NSF of China (No.~12101362), the Taishan Scholar Project and the NSF of Shandong Province (No.~ZR2022YQ01). J. Zheng was supported by National key R\&D program of China: 2021YFA1002500 and NSF grant of China (No. 12271051).  Y. Zheng was supported by NSF of China under Grant No.~11901350 and No.~12371172.}

\begin{abstract}
In this article, we consider an energy-critical complex Ginzburg-Landau equation in the exterior of a smooth compact strictly convex obstacle.
We prove the global well-posedness of energy-critical complex Ginzburg-Landau equation in an  exterior domain by the concentration-compactness/rigidity theorem method. As corollaries of our main result, we establish both the existence of global weak solutions to energy-critical defocusing nonlinear Schr\"odinger equations and the global well-posedness theory for energy-critical semi-linear heat equation in exterior domains.

\bigskip

\bigskip

\end{abstract}

\maketitle
\tableofcontents

\section{Introduction}
The class of complex Ginzburg-Landau equations has drawn the attentions of many mathematicians and physicists, because of its close relation to many physical phenomenons in fluids and superconductivity \cite{BBH-book,GL,Guo96,WGZ,PR-book}.
Especially, the complex Ginzburg-Landau equations could be derived from the B\'enard convection and turbulence; see for instance \cite{GJL}. The physical meaning of equations in fluids in an exterior domain (see for instance \cite{G2011}) and the close relation between complex Ginzburg-Landau equation and equations in fluids such as Navier-Stokes equations motivate us to consider the following energy-critical complex Ginzburg-Landau (CGL) equation in an exterior domain $\Omega$ of a smooth compact strictly convex obstacle in $\mathbb{R}^N$, $N \ge 3$, with Dirichlet boundary condition:
\begin{align} \label{eq:gl}
\left\{
\begin{aligned}
&  u_t - z \Delta u=-z|u|^{\f {4}{N -2}} u\quad \text{in } \ \mathbb{R}_+\times\Omega,\\
& u(0, x) = u_0 (x)  \text{ in } \Omega,\\
& u(t,x)|_{x\in \partial \Omega} = 0,
\end{aligned}
\right.
\end{align}
where $z \in \mathbb{C}$, $\Re z > 0$ and $u$ is a complex-valued function of $(t,x) \in \mathbb{R}_+ \times \Omega$.
This equation is the $L^2$-gradient flow of the energy functional
\begin{equation*}
E(u) = \int_{\Omega}  \frac12 |\nabla u|^2 + \frac{N -2}{2 N }|u|^\frac{2N }{N -2} \,\mathrm{d}x.
\end{equation*}

The complex Ginzburg-Landau equation is closely related to the nonlinear Schr\"odinger equation (NLS). As one easily observes, taking $z=i$ in \eqref{eq:gl}, we obtain the energ-critical defocusing NLS with Dirichlet boundary condition:
\begin{align} \label{eq:nls}
\left\{
\begin{aligned}
&  iu_t+\Delta u = |u|^{\f {4}{N -2}} u\quad \text{in } \ \mathbb{R}_+\times\Omega,\\
& u(0, x) = u_0 (x)  \text{ in } \Omega,\\
& u(t,x)|_{x\in \partial \Omega} = 0,
\end{aligned}
\right.
\end{align}
Note that the Schr\"odinger equation with the above Dirichlet boundary condition is properly intepretated in mathematical quantum mechanics; see \cite{RS75}.
It is well-known that the inviscid limit of the complex Ginzburg-Landau equation is the nonlinear Schr\"odinger equation, and a rigorous limit theory has been established in \cite{CGZ2,HW,W0}.
%One reason to consider such a  boundary condition as in \eqref{eq:gl} is that in the case of Schr\"odinger equation, this condition can be properly intepretated in mathematical quantum mechanics; see \cite{RS75}.
Furthermore, in physics, concerning a superfluid passing through an obstacle leads to the study of defocusing (cubic) nonlinear Schr\"odinger equation in an exterior domain. Indeed, the model mathematical problem for a superfluid passing an obstacle is of the form $$i\partial_tu+\Delta u + |u|^2u=0\quad \text{in}\quad \mathbb{R}_+\times \Omega,$$
together with some Neumann boundary condition; see for instance \cite{FPR1992,JP2000}.
In the fundamental work \cite{Lin-Zhang-2005}, F.-H. Lin and P. Zhang consider the semiclassical limit of the above model in an exterior domain.
In an important recent work \cite{KVZ1}, R. Killip, M. Visan and X. Zhang establish the global well-posedness and scattering of finite energy solutions to the defocusing nonlinear Schr\"odinger equation in an exterior domain $\Omega\subset \mathbb{R}^3$; see also the work of D. Li, G. Xu, and X. Zhang \cite{LXZ} on the study of defocusing energy-critical nonlinear Schr\"odigner equation in the exterior of the unit ball in five and seven dimensions and the work \cite{LSZ12} on the global well-posedness and scattering of defocusing energy-critical NLS in the exterior of unit ball in three dimension. The result \cite{KVZ1} was soon extended to the focusing case in dimension three in \cite{KVZ2}.

Consider the $\CGL$ in $\R^N $ in dimension $N \ge 3$
 \begin{align}\label{eq1.3v27}
\begin{cases}
u_t - z \Delta u = -z|u|^\frac4{N -2} u  \quad \text{in } \ \mathbb{R}_+\times\mathbb{R}^N , \\
u(0,x) = u_0(x).
\end{cases}
\end{align}
C. Huang and B. Wang \cite{HW} establish the following global well-posedness theory of \eqref{eq1.3v27}. We also refer to \cite{CGZ} for an alternative proof when $N\in \{3,4\}$ based on the concentration-compactness/rigidity theorem argument developed in \cite{KK11,KM}.

\begin{thmA}[\cite{HW}] \label{th1.1v37}
Let $u_0 \in \dot{H}^1( \mathbb{R}^N )$ with $N \ge 3 $. Then the solution $u$ of \eqref{eq1.3v27} is global and satisfies $\lim\limits_{t \to \infty } \|u(t) \|_{\dot{H}^1} = 0$.
\end{thmA}

Motivated by its physical background and by the corresponding work on NLS, in this article, we study the global well-posedness of $\CGL$ in an exterior domain in $ \dot{H}_D^1$; see Definition \ref{de1.1}. Our main result reads as follows.
\begin{theorem}\label{th1.3}
Fix $N  \in \left\{ 3, 4, 5 \right\}$ and $u_0 \in  \dot{H}_D^1( \Omega )$. Then there exists a unique solution $u$ to \eqref{eq:gl} which is global in time and satisfies $\lim\limits_{t \to \infty } \|u(t) \|_{\dot{H}^1_D(\Omega)} = 0$.
\end{theorem}

Note that if  $\Omega=\mathbb{R}^N $, then the problem is scaling invariant. More precisely, for each $\lambda>0$, the scaled mapping
\begin{equation*}
	u(t,x) \mapsto u^\lambda (t,x) : = \frac1{ \lambda^\frac{N -2}2} u \left( \frac{t}{ \lambda^2}, \frac{x} \lambda \right)
\end{equation*}
leaves the class of solutions to $\CGL_{\mathbb{R}^N }$ invariant. This scaling invariance also makes $\dot{H}_x^1$ as the critical space. Since the presence of obstacle does not change the intrinsic dimensionality of the problem, we may regard equation \eqref{eq:gl} as being critical for data in $\dot{H}_D^1( \Omega)$. In the Euclidean setting, the problem is invariant under spatial translations; this means that one may employ the full power of harmonic and analytic tools. Indeed, much of the recent surge of progress in the analysis of dispersive equations is based on the incorporation of this powerful technology.
When we consider \eqref{eq:gl} in exterior domains, the loss of spatial translation invariance  becomes a serious issue and many of the tools developed for the case $\Omega=\mathbb{R}^N $ fails.
Working in exterior domains also destroys the scaling symmetry. Due to the presence of a boundary, suitable scaling and spatial translations lead to the study of $\CGL$ in different geometries.

As an application of our main result, we obtain the following global existence of  weak solutions to the energy critical defocusing NLS via an limiting argument from \cite{CGZ2}. Recall that a complex-valued function $u$ on a time interval $I \subseteq \mathbb{R}$ is a weak-${H}^1_D$ solution to \eqref{eq:nls} if $u\in L^\infty_t {H}_D^{1}(I \times \Om)$, $\partial_tu \in L^\infty_t H_x^{-1}(I \times \Om)$ such that \eqref{eq:nls} holds for almost every $t\in I$ in $H^{-1}_x$.

\begin{corollary}[Global weak solutions and weak-strong uniqueness for NLS on exterior domain]\label{cor:nls}
Fix $N  \in \{ 3, 4, 5 \}$ and $u_0 \in  {H}_D^1( \Omega )$. Then
\begin{itemize}
			\item[1)] Global existence: There exists a global weak solution $u\in L_t^\infty {H}_D^1(\mathbb{R} \times \Om ) \cap C_t^0 L^2_x(\mathbb{R} \times \Om )$ of \eqref{eq:nls} satisfying the energy inequality
			\begin{align}\label{eq1.4v26}
				E(u(t)) \le E(u_0)
			\end{align}
			and mass conservation
			\begin{align*}
				M(u(t)) = M(u_0).
			\end{align*}
			
			\item[2)] Weak-strong uniqueness: There exists a universal constant $C= C \left(T\right)>0$ such that if $\tilde{u}  $ is the global strong solution to  \eqref{eq:nls} with the initial data $u_0 \in C_0^\infty(\Om)$ and if ${u}$ is the global weak solution to \eqref{eq:nls} with initial datum $u_0$, then $u\equiv \tilde{u} $.
\end{itemize}
		
\end{corollary}

As was pointed out earlier, the global well-posedness for NLS in exterior domains has been proved in dimension three by  R. Killip, M. Visan and X. Zhang \cite{KVZ1}. However, due to the loss of Strichartz estimates in exterior domains in higher dimensions $N\ge 4$, the existence of global strong solutions remains open in higher dimensions $N\geq 4$. On the other hand, the above result implies that the global weak solution obtained in Corollary \ref{cor:nls} is indeed the only possible strong solutions to \eqref{eq:nls}.

Observe that when $z=1$ (or more generally, $z \in \mathbb{R}_+$), \eqref{eq:gl}  reduces to the energy-critical semi-linear heat equation
\begin{align} \label{eq:heat}
  \left \{
  \begin{aligned}
  & u_t - \Delta u =  -|u|^{\f {4}{N -2}} u \quad \text{in}\quad  \mathbb{R}_+\times \Omega,\\
  & u(0, x ) = u_0 (x)   , \\
  & u(t,x)|_{x\in \partial \Omega} = 0.
  \end{aligned}
  \right.
\end{align}
Semi-linear heat equations with Dirichlet boundary values have been extensively studied in the literature; see for instance \cite{W80,NS85,BC96,R98,miao-zhang-2004,Jiang-Lin-2024} and the references therein for more results. As a special case of Theorem \ref{th1.3}, we obtain the following GWP result for \eqref{eq:heat}.

\begin{corollary}[Global well-posedness for energy-critical semilinear heat equation]\label{coro:heat exterior}
For any $N \in \{3,4, 5 \}$, let $u_0 \in  \dot{H}_D^1( \Omega )$. Then there exists a unique solution $u$ to \eqref{eq:heat} which is global in time and satisfies $\lim\limits_{t \to \infty } \|u(t) \|_{\dot{H}^1_D(\Omega)} = 0$.
\end{corollary}

Somewhat surprisingly, this result seems to be overlooked in the literature. Our argument does not rely on the maximal principle and it is not clear for us whether one can give a proof of Corollary \ref{coro:heat exterior} with the maximal princple; see \cite{BC96,R98} for more discussions.

It is worth pointing out that our argument can be used to study the following slightly more general complex Ginzburg-Landau equation considered in \cite{HW}
	\begin{align}\label{eq1.4v40}
		\begin{cases}
			u_t -  ( b+ i ) \Delta u= - ( a + i) |u|^\frac4{N -2} u \text{ in } \mathbb{R}_+ \times \Omega, \\
			u(0,x) = u_0(x) \text{ in } \Omega, \\
			u(t,x)|_{x \in \Omega} = 0,
		\end{cases}
	\end{align}
where $a, b > 0$. More precisely, fix $N \in \{3, 4, 5 \}$ and $u_0 \in  \dot{H}_D^1( \Omega )$. An easy adoption of our argument shows that there exists a unique solution $u$ to \eqref{eq1.4v40} which is global in time and satisfies $\lim\limits_{t \to \infty } \|u(t) \|_{\dot{H}^1_D(\Omega)} = 0$. We leave the details to the interested readers.
	
Before introducing the methodology of our approach, we would like to point out that studying the GWP/scattering of CGL/NLS on exterior domains can be viewed as an intermediate step towards a similar theory on general Riemannian manifolds, where the geometry of manifolds plays a crucial role; see for instance the nice survey in \cite{BGT07}. Another natural question is to study the CGL with other type of boundary conditions, such as period/almost period initial--boundary data. This leads to some related works in $\mathbb{R}^2\times \mathbb{T}$ or the cylinder; see for instance \cite{CGYZ20,CGZ20}. In our later works, we shall try to extend the theory to more (non-Euclidean) geometric settings or with different boundary conditions.

\subsection{Methodology}
To the best of the authors' knowledge, our paper is the first result to study global well-posedness of the energy-critical nonlinear dissipative equations, especially the energy-critical semi-linear heat equations, in exterior domains in the energy space. In our proof, we shall use the concentration-compactness/rigidity theorem argument in \cite{KM,KK11}, together with some further development in \cite{CGZ} on the CGL in $\mathbb{R}^N $. In particular, our argument for the global well-posedness of the semi-linear heat equation does not rely on the maximum principle.

The general strategy to prove Theorem \ref{th1.3} consists of two main steps. In the first step, we use the concentration-compactness argument to construct the critical element. First of all, adopting the concentration-compactness argument to the problem in exterior domains will cause us a great deal of trouble, as $\CGL$ in the exterior domain $\Omega$ does not enjoy the scaling invariance or the translation invariance. It is easy to see that the scaling and translation of $\Omega$ can be $\mathbb{R}^N$, $\mathbb{R}_+^N$ or just $\Omega$.
 To applying the concentration-compactness argument, we first need to establish the linear profile decomposition of CGL in exterior domains, which describes the lack of compactness of the embedding $e^{tz\Delta_\Omega}: \dot{H}_D^1( \Omega) \to L^{\frac{2(N+2)}{N-2}}_{t,x}( \mathbb{R}_+ \times \Omega) $. This technique has been applied to NLS and nonlinear wave equation (NLW) on $\mathbb{R}^N $ (\cite{BG,Ke}), and to NS equation on $\mathbb{R}^N $ (\cite{G}). The proof for this is inspired by the linear profile decomposition of NLS in exterior domains \cite[Theorem 5.6]{KVZ1} and NLW  in exterior domains \cite{GG}. This is provided in Subsection \ref{subse3.1}. 
 %{\color{blue}briefly introduce the main technical difficulty in this step}
%The linear profile decomposition of $\CGL$ in exterior domains follows from the same argument as the linear profile decomposition of NLS in exterior domains \cite[Theorem 5.6]{KVZ1}; see also the case of NLW in exterior domains \cite{GG}.
%The proof follows the argument of the linear profile decomposition of NLS %the nonlinear Schr\"odinger equation %(NLS)
% and nonlinear wave equation (NLW) on $\mathbb{R}^N $ (\cite{BG,Ke}), and to Navier-Stokes equation (NS) on $\mathbb{R}^N $ (\cite{G}), together with the linear profile decomposition of NLS in exterior domains \cite[Theorem 5.6]{KVZ1} and NLW in exterior domains \cite{GG}.
 Once we have obtained the linear profile decomposition, we need to establish the nonlinear profile decomposition. For notational simplicity, we write below $\CGL_\Omega$ to indicate CGL in the domain $\Omega$. The nonlinear profiles are defined to be solutions of the $\CGL_\Omega$, with initial datum being each profile in the linear profile decomposition. As we do not have scaling or translation invariance for the exterior domain $\Omega$, to find a good approximation of the nonlinear profiles, we have to consider the following three different scenarios arising from scaling/translation operations on $\Omega$:
\begin{itemize}
	\item The rescaled domains $\Omega_n$ (where $\Omega_n = \lambda_n^{- 1}  \left( \Omega - \left\{ x_n \right\} \right)$) expand to fill $\mathbb{R}^N$ as $\lambda_n\to \infty$;
	\item The rescaled domains $\Omega_n^c$
	(where $\Omega_n = \lambda_n^{- 1}  \left( \Omega - \left\{ x_n \right\} \right)$) are retreating to infinity as $\frac{\dist(x_n,\Omega^c)}{\lambda_n}\to \infty$;
	\item The rescaled domains $\Omega_n$ (where $\Omega_n : = \lambda_n^{- 1} R_n^{- 1}  \left( \Omega -  \left\{  x_n^* \right\}  \right)$)
	expand to fill a half space as $\lambda_n\to 0$ and $\frac{\dist(x_n,\Omega^c)}{\lambda_n}\to \beta>0$. Here, $x_n^\ast \in \partial \Omega$ be such that $ \left|x_n - x_n^* \right| = \dist(x_n,\Omega^c)$ and $R_n \in SO( N )$ satisfies
	$R_n e_N  = \frac{x_n - x_n^*}{ \left|x_n - x_n^* \right|}$.
\end{itemize}
We will treat each of the above three cases in $\CGL_\Omega$ setting separately in Subsection \ref{subse3.2}. In particular, we construct the limiting approximated equations to $\CGL_\Omega$ in each scenarios, which is different from the case of NLS and important to show the pre-compactness of the sequence of critical elements. Then a variant of the argument in \cite{CGZ,KK11,KM} allows us to establish the existence of the critical element. This completes the proof of first step.

In the second step, we use the Morawetz estimate established on the effect of dissipation to exclude the possibility of blowup in finite time. There are a couple of technical issues one has to overcome. In particular, we construct the approximate solutions to $\CGL_\Omega$ for different scenarios, which is different from the dispersive equations.
Secondly, we develop a Morawetz type estimate for complex Ginzburg-Landau equations to exclude the existence of the critical element.
 In this step, we shall use several techniques, such as the Littlewood-Paley theory and Bernstein estimates in exterior domains.

The structure of the paper is as follows. In Section \ref{se2}, we will present the local well-posedness theory. In Section \ref{se3}, we show the existence of critical element. In Section \ref{se6v23}, we prove our main Theorem \ref{th1.3}.

\textbf{Notations.}
Throughout this paper, unless specified, $\Omega$ shall be the exterior domain of a smooth compact strictly convex obstacle in $\mathbb{R}^N $  and $-\Delta_\Omega$ is the Dirichlet Laplacian (associated to $\Omega$), which is a self-adjoint operator on $L^2(\Omega)$ with domain $H_D^1( \Omega)$.

Let $\mathbb{N}$ be the set of natural numbers. We will introduce the halfspace
\begin{align*}
\mathbb{R}^N_+ =  \left\{ x =  \left(x_1, \cdots, x_N \right) \in \mathbb{R}^N : x_N  > 0  \right\}.
\end{align*}
For any time interval $I$, we define
\begin{align*}
& S^0(I \times \Omega ): = L_t^\infty L_x^2 \cap L_t^2 L_x^\frac{2N }{N -2}   (I \times \Omega)
 , \\
& \dot{S}^1(I \times \Omega ) : =  \left\{u : I \times \Omega \to \mathbb{C},  \left( - \Delta_\Omega \right)^\frac12 u \in S^0(I) \right\}.
\end{align*}
We define $N^0(I \times \Omega )$ to be the dual Strichartz space and
\begin{align*}
\dot{N}^1(I \times \Omega ) : =  \left\{ F: I \times \Omega \to \mathbb{C},
\left( - \Delta_\Omega \right)^\frac12 F \in N^0(I) \right\}.
\end{align*}
Similarly, we can define these spaces for functions defined on $I \times \mathbb{R}^N $ and $I \times \mathbb{R}_+^N $.

If $u$ is a function defined on $\mathbb{R}^N $, then for  each $\lambda_0 >0$ and $x_0 \in \mathbb{R}^N $, we will write
\begin{align*}
	u_{[ \lambda_0 , x_0 ]}(x) = \frac1{\lambda_0^{\frac{N -2}{2}}} u \left( \frac{x- x_0 }{\lambda_0  } \right).
\end{align*}
We introduce the Littlewood-Paley projectors adapted to the Dirichlet Laplacian on $\Omega$.
Let $\phi:[0, \infty ) \to [0,1]$ be a smooth function with
\begin{align*}
\phi(\xi) : =
\begin{cases}
1,\ 0 \le \xi \le 1, \\
0,\ \xi \ge 2.
\end{cases}
\end{align*}
For each dyadic number $ L   \in 2^{\mathbb{Z}}$, we set
\begin{align*}
\phi_L (\xi) : = \phi \left( \frac{\xi}L \right)
\quad\text{and}\quad \psi_L (\xi) : = \phi_L (\xi) - \phi_{\frac{L }2} ( \xi).
\end{align*}
Then the Littlewood-Paley projectors adapted to the Dirichlet Laplacian on $\Omega$ are defined to be
\begin{align*}
P_{\le L }^\Omega : = \phi_L  \left( \sqrt{- \Delta_\Omega}  \right),
\
P_L^\Omega : = \psi_L   \left( \sqrt{ - \Delta_\Omega}  \right), \text{ and } P_{> L }^\Omega : = I - P_{\le L }^\Omega.
\end{align*}
We will also use $P_L^{\mathbb{R}^N }$ to represent the Littlewood-Paley projectors adapted to the Laplacian on $\mathbb{R}^N $, and use $P_L^{\mathbb{R}_+^N }$ to represent the Littlewood-Paley projectors adapted to the Laplacian on $\mathbb{R}_+^N $.

The following two basic estimates are well-known; {see for instance \cite{KVZ3}}.
\begin{lemma}[Bernstein estimates]\label{lea.6v65}
{For any $f\in C_c^\infty(\Omega)$, we have
\begin{align*}
\left\| |\nabla |^s P^\Omega_L  f  \right\|_{L_x^p} \sim&  L^s  \left\|P_L^\Omega  f  \right\|_{L_x^p}, \\
\left\| P_{\le  L }^\Omega  f  \right\|_{L_x^q } \lesssim&  L^{ \frac{N}p - \frac{N}q }  \left\|P_{\le  L }^\Omega f  \right\|_{L_x^p },
\end{align*}
for any $1<p<\infty$ and $s\in\R$. Moreover,
\begin{align*}
\left\|P_L^\Omega   f  \right\|_{L_x^q } \lesssim L^{\frac{N}p - \frac{N}q }  \left\|P_L^\Omega   f  \right\|_{L_x^p},
\end{align*}
for all $1\leq p<q\leq\infty$. The implicit constants depend only on $N,p,q$ and $s$.}

\end{lemma}

\begin{lemma}[Square function estimates]\label{lea.7v65}
Given a Schwartz function $f$ on $\Omega$, we have for $1 < p < \infty$,
\begin{align*}
\left\|  \left( \sum\limits_L   \left|P_L^\Omega  f(x) \right|^2  \right)^\frac12  \right\|_{L_x^p} \sim \|f \|_{L_x^p}.
\end{align*}

\end{lemma}

\section{The local well-posedness theory}\label{se2}

In this section, we present the local well-posedness theory of \eqref{eq:gl}.

\begin{definition}\label{de1.1}
For $s \ge 0$ and $1 < p < \infty$, the homogeneous fractional Sobolev space $\dot{H}_D^{s, p}(\Omega)$ is the completion of $C_c^\infty( \Omega)$ with respect to the norm
\begin{align*}
\|f\|_{\dot{H}_D^{s, p} ( \Omega)} : =  \left\| \left( - \Delta_\Omega \right)^\frac{s}2 f  \right\|_{L^p(\Omega)}.
\end{align*}
For notational simplicity, when $p=2$, we shall write $\dot{H}_D^{s}$ instead of  $\dot{H}_D^{s,p}$.
\end{definition}
We need the following equivalence of fractional Sobolev spaces.
\begin{theorem}[\cite{KVZ3}] \label{th2.3}
 Let $1 < p < \infty$. If $0 \le s < \min \left\{ 1 + \frac1p, \frac{N }p \right\}$, then
\begin{align*}
 \left\| \left(- \Delta_{\mathbb{R}^N }  \right)^\frac{s}2 f \right\|_{L^p} \sim_{N , p, s}  \left\| \left( - \Delta_\Omega \right)^\frac{s}2 f  \right\|_{L^p}
\end{align*}
for all $f \in C_c^\infty( \Omega)$.
\end{theorem}

The following Leibniz estimate for fractional derivatives is well-known.
\begin{lemma}\label{le2.4}
For all $f, g \in C_c^\infty( \Omega)$, we have
\begin{align*}
\left\|  \left( - \Delta_\Omega \right)^\frac{s}2 ( f g)  \right\|_{L^p }
\lesssim  \left\| \left( - \Delta_\Omega \right)^\frac{s}2 f \right\|_{L^{p_1}} \|g \|_{L^{p_2}} + \|f \|_{L^{q_1}}  \left\| \left( - \Delta_\Omega \right)^\frac{s}2 g \right\|_{L^{q_2}},
\end{align*}
with the exponents satisfying $1 < p, p_1, q_2 < \infty$, {$1 < p_2 , q_1 \le \infty$}, $\frac1p = \frac1{p_1} + \frac1{p_2} = \frac1{q_1} + \frac1{q_2}$, and
$0 < s < \min  \left\{ 1 + \frac1{p_1}, 1 + \frac1{q_2}, \frac{N}{p_1}, \frac{N}{q_2} \right\}. $
\end{lemma}
From now on, without loss of generality, we always assume that $0 \in \Omega^c$. We shall use $d_{\Omega}(x,y)$ to denote the geodesic distance between $x$ and $y$ in $\Omega$, $\diam := \diam  \left(\Omega^c \right) $ to denote the diameter of the obstacle and $d(x) : =\dist(x, \Omega^c)$ to denote the distance of a point $x \in \mathbb{R}^N $ to the obstacle.

For the linear part of CGL $\partial_t u =  z \Delta_\Omega u$, it is well-known that the operator $z\Delta_\Omega$ with domain $H_0^2( \Omega)$ generates a semi-group $ \left(e^{z \Delta_\Omega} \right)_{z\in \mathbb{C}_+}$ on $L^2( \Omega)$, where $\mathbb{C}_+$ is the complex half-plane $\Re z\geq 0$. Moreover, the semi-group $ \left(e^{z \Delta_\Omega} \right)_{z \in \mathbb{C}_+}$ is analytic with heat kernel $p_z(x,y)= p(z,x,y)$. Furthermore, when $z=t>0$, it holds
\begin{align*}
	e^{t \Delta_\Omega} \phi(x) =  p_t(x,y) \ast \phi = \int p_t(x,y) \phi(y) \,\mathrm{d}y.
\end{align*}

\begin{theorem}[Upper bound on the kernel of Dirichlet Laplacian in exterior domains, \cite{Coulhon-Sikora-2008,Zhang2003}]\label{th2.2}
For all $x,y \in \Omega$ and $z\in \mathbb{C}$ with $\Re z \ge 0$, there exists a constant $c>0$ such that the following upper bound for the heat kernal holds:
\begin{align*}
	\left| p_z(x,y) \right|
	\lesssim  \left(1+\frac{d_{\Omega}^2(x,y)}{4\Re z}\right)^{(N -1)/2} e^{- \frac{d_{\Omega}(x,y)^2}{4\Re z} } (\Re z)^{- \frac{N}2}\lesssim \left(1+\frac{d(x)^2+d(y)^2}{c\Re z}\right)^{(N -1)/2} e^{-\frac{|x-y|^2}{c\Re z} } (\Re z)^{- \frac{N}2}.
\end{align*}
In the particular case when $z = t>0$, we have
\begin{align*}
 \left|p_t(x, y) \right| \lesssim  \min  \left\{\frac{d(x)}{ \min  \left\{  \sqrt{t} , \diam  \right\}}, 1 \right\} \cdot \min  \left\{\frac{d(y)}{ \min  \left\{ \sqrt{t}, \diam  \right\}}, 1  \right\} e^{- \frac{ |x- y|^2}{ct}} t^{- \frac{N }2}.
\end{align*}

\end{theorem}
By the arguments in \cite{KT,miao2018}, we have the following Strichartz estimate.
\begin{theorem}[Strichartz estimates]\label{th2.8}
Let $I$ be a time interval. Then the solution $u$ to the forced $\CGL$
equation $ \partial_t u - z  \Delta_\Omega u =  F $ satisfies the estimate
\begin{align*}
\| u \|_{S^0(I \times \Omega )} \lesssim \|u(t_0) \|_{L^2( \Omega)} + \|F \|_{N^0(I \times \Omega )}
\end{align*}
for any $t_0 \in I$. In addition, we have
\begin{align*}
\|u \|_{\dot{S}^1(I \times \Omega )} \lesssim \|u(t_0) \|_{\dot{H}_D^1 ( \Omega)} + \|F\|_{\dot{N }^1(I \times \Omega )}
\end{align*}
for any $t_0 \in I$.

\end{theorem}

Embedding functions on the half-space $\mathbb{R}_+^N $ as functions on $\mathbb{R}^N $ that are odd under reflection in $\partial \mathbb{R}_+^N $, we immediately see that the whole range of Strichartz estimates, also holds for the free propagator $e^{tz \Delta_{\mathbb{R}_+^N  }}$.

\begin{theorem}[Local well-posedness, \cite{miao2018}]\label{th2.1v3}
Assume $u_0 \in \dot{H}_D^1( \Omega)$. Then the following conclusions hold:
\begin{itemize}
\item Local existence: There exists a unique, maximal-lifespan solution to \eqref{eq:gl} on $[0, T_{\max} (u_0) )$.
\item Blow-up criterion: If $T_{\max} < \infty$, then $ \|u \|_{L_{t,x}^\frac{2(N +2)}{N -2} ([0, T_{\max} (u_0))  \times \Omega)}  = \infty$.
\item Small data global existence: There is $\epsilon_0 > 0$ such that if $ \left\|e^{tz \Delta_\Omega}  \right\|_{L_{t,x}^\frac{2( N +2)}{ N  - 2}(\mathbb{R}_+ \times \Omega )} \le \epsilon_0,$ then the solution $u$ is global and $\|u\|_{L_{t,x}^\frac{2(N +2)}{N - 2} (\mathbb{R}_+ \times \Omega )} \lesssim \epsilon_0$. This holds in particular when $\|u_0 \|_{\dot{H}_D^1}$ is small enough.
\end{itemize}
\end{theorem}

As a consequence of the Strichartz estimate and the local well-posedness, we have a stability result for $\CGL$ on $\Omega$. This follows the argument of the proof of the stability result of $\CGL$ on  $\mathbb{R}^N $; see \cite[Proposition 2.4]{CGZ}. Note that we also have stability result for $\CGL$ on $\mathbb{R}_+^N $, which follows directly by embedding solutions on the half space as solutions on $\mathbb{R}^N $ that are odd under reflection in $\partial \mathbb{R}_+^N $. In particular, we can regard $-\Delta_{\mathbb{R}_+^N }$ as the restriction of $- \Delta_{\mathbb{R}^N }$ to odd functions. For simplicity, we shall write $\CGL_{\Omega}$, $\CGL_{\mathbb{R}^N }$ and $\CGL_{\mathbb{R}^{N}_+}$ below to denote $\CGL$ on $\Omega$, $\mathbb{R}^N $ and $\mathbb{R}^N_+$, respectively.

\begin{theorem}[Stability for $\CGL$] \label{th2.9}
Let $I$ be a compact time interval and let $\tilde{u}$ be an approximate solution to \eqref{eq:gl} on $I \times \mathcal{O}$, with $\mathcal{O}$ being one of the three domains $\Omega$, $\mathbb{R}^N $, or $\mathbb{R}^N_+$, in the sense that
\begin{align*}
 \partial_t  \tilde{u}  - z  \Delta_{\mathcal{O}} \tilde{u} =   - z   \left|\tilde{u} \right|^\frac4{N -2}  \tilde{u } + e
\end{align*}
for some function $e$. Assume that
\begin{align*}
 \left\|\tilde{u}  \right\|_{L_t^\infty \dot{H}_D^1( I \times \mathcal{O})} \le A
 \quad \text{and}\quad \left\|\tilde{u } \right\|_{L_{t,x}^{ \frac{2(N+2)}{N -2} }( I \times \mathcal{O})} \le  B
\end{align*}
for some positive constants $A$ and $B$.

Let $t_0 \in I$ and let $u_0 \in \dot{H}_D^1(\mathcal{O})$ satisfy
\begin{align*}
\left\| u_0 - \tilde{u}(t_0)  \right\|_{\dot{H}_D^1} \lesssim  A'
\end{align*}
for some positive constant $A'$.

Assume also the smallness condition holds: for some $0 < \epsilon < \epsilon_1 = \epsilon_1 \left(A, A',  B  \right)$, it holds
\begin{align*}
 \left\| \left( - \Delta_\mathcal{O} \right)^\frac12  e^{(t - t_0) {z} \Delta_\mathcal{O}}  \left(u(t_0) - \tilde{u}_0 \right)  \right\|_{L_t^{ \frac{2( N +2)}{N - 2} } L_x^\frac{2 N ( N +2) }{N^2 + 4 }( I \times \mathcal{O})}
+  \left\| \left( - \Delta_\mathcal{O} \right)^\frac12 e \right\|_{N^0(I)} \le \epsilon.
\end{align*}
Then there exists a unique strong solution $u: I \times \mathcal{O} \to \mathbb{C}$ to \eqref{eq:gl} with initial data $u_0$ at $t = t_0$ satisfying
\begin{align*}
\left\|u - \tilde{u}  \right\|_{L_{t,x}^{ \frac{2( N +2)}{N  - 2} } (I \times \mathcal{O})} \le C \left( A , A',  B  \right) \epsilon, \\
 \left\|  \left( - \Delta_\mathcal{O} \right)^\frac12 \left(u - \tilde{u} \right) \right\|_{S^0(I\times \mathcal{O})} \le C \left(A, A',  B  \right) E', \\
\left\| \left( - \Delta_\mathcal{O}
\right)^\frac12 u  \right\|_{S^0(I \times \mathcal{O})} \le C \left(A, A', B \right).
\end{align*}
\end{theorem}
In our later proofs, we shall need the following persistence of regularity result for $\CGL_{\mathbb{R}^N }$ and $\CGL_{\mathbb{R}^N_+}$. {The case for $\CGL_{\mathbb{R}^N }$ follows from \cite{CGZ}}. The other case follows by embedding solutions on the half space as solutions on $\mathbb{R}^N  $ that are odd under reflection in $\partial \mathbb{R}^N_+ $. In particular, one may regard $- \Delta_{\mathbb{R}^N_+ }$ as the restriction of $- \Delta_{\mathbb{R}^N  }$ to odd functions.
\begin{lemma}[Persistence of regularity for $\CGL_{\mathbb{R}^N }$ and $\CGL_{\mathbb{R}^N_+}$] \label{le2.10}
{Let $s\geq0$ when $N\in\{3,4\}$, and $0\leq s<1+\frac{4}{N-2}$ when $N\geq5$}, and let $I$ be a compact time interval and $\mathcal{O}$ be either $\mathbb{R}^N $ or $\mathbb{R}_+^N $. Assume $u: I \times \mathcal{O} \to \mathbb{C}$ is a solution to $\CGL_{\mathcal{O}}$ satisfying
\begin{align*}
E(u) \le  A < \infty
\intertext{ and }
\| u \|_{L_{t,x}^{ \frac{2(N +2)}{N  -2 } } ( I \times \mathcal{O})} \le B  < \infty.
\end{align*}
If $u(t_0) \in \dot{H}^s (\mathcal{O})$ for some $t_0 \in I$, then
\begin{align*}
 \left\|  \left( - \Delta_{\mathcal{O}}  \right)^\frac{s}2 u  \right\|_{S^0(I \times \mathcal{O} )}
\le C( A, B  ) \|u(t_0) \|_{\dot{H}^s (\mathcal{O})} .
\end{align*}
\end{lemma}
We conclude this section with an important estimate to be used later in Section \ref{se3}.
\begin{proposition}\label{co2.14}
Given $w_0 \in \dot{H}_D^1(\Omega)$,
\begin{align*}
\left\| \nabla e^{   t {z} \Delta_\Omega} w_0 \right\|_{L_{t,x}^\frac{ N +2}{ N -1 }  ([\tau, \tau + T] \times \{ |x- x_0  | \le R\} )} %\\
&\lesssim T^{\left(1-\frac{ N (4- N )(2-\alpha)}{2( N -1)}\right)\frac{\beta( N -1)}{ N +2} + \frac{2}{( N +2)(1+2\theta)}\frac{\alpha( N -1)}{ N +2} } R^{\frac{ N (1+\theta)-2(1-\theta) + \theta ( N +2)}{( N +2)(1+2\theta)}\frac{\alpha(  N -1)}{ N +2}}  \\
&\quad \times \left\| e^{   t {z} \Delta_\Omega} w_0 \right\|_{L_{t,x}^{\frac{2( N +2)}{ N -2}}}^{\frac{\theta}{1+2\theta}\frac{\alpha( N -1)}{ N +2}} \|w_0\|_{\dot{H}_D^1} ^{\frac{\beta( N -1)}{ N +2}+\frac{1+\theta}{1+2\theta}\frac{\alpha( N -1)}{ N +2}},
\end{align*}
uniformly in $w_0$ and the parameters $R, \tau, T> 0$,  and $x_0 \in \mathbb{R}^N  $, where $\alpha + \beta = \frac{N +2}{N -1} $ and $\theta \in (0, 1)$.
\end{proposition}

\begin{proof}
For simplicity, we only prove the estimates on the time interval $[0,T]$.
 Applying H\"older's inequality with $\alpha + \beta = \frac{N +2}{N -1 }$, we obtain
    \[
    \left\| \nabla \left(  e^{   t {z} \Delta_\Omega} w_0 \right)  \right\|_{L_{t,x}^\frac{N +2}{N -1 }  ([0, T] \times \{ |x- x_0  | \le R\} )} \le \left\| \nabla
    \left( e^{   t {z} \Delta_\Omega} w_0 \right)  \right\|_{L_{t,x}^2  ([0, T] \times \{ |x- x_0  | \le R\} )}^{\frac{\alpha( N -1)}{ N +2}}
     \left\| \nabla \left(  e^{   t {z} \Delta_\Omega} w_0  \right) \right\|_{L_{t,x}^{\frac{2\beta}{2-\alpha}}  ([0, T] \times \Omega )}^{\frac{\beta( N -1)}{ N +2}}.
    \]
By H\"older's inequality and the Strichartz inequality, we have
    \[
    \begin{aligned}
        \left\| \nabla  \left( e^{   t {z} \Delta_\Omega} w_0 \right)  \right\|_{L_{t,x}^{\frac{2\beta}{2-\alpha}}  ([0, T] \times \Omega )} & \lesssim T^{1-\frac{ N (4- N )(2-\alpha)}{2( N -1)}} \left\| \nabla  \left( e^{   t {z} \Delta_\Omega} w_0  \right) \right\|_{ L_t^{ \frac{4\beta( N -1)}{ N (4- N )} } L_{x}^{\frac{2\beta}{2-\alpha}}  ([0, T] \times \Omega )} \\
        & \lesssim T^{1-\frac{ N (4- N )(2-\alpha)}{2( N -1)}} \|w_0\|_{\dot{H}^1_D}.
    \end{aligned}
    \]

    Given $L >0$ and $\theta \in (0, 1)$, using the H\"older, Bernstein, and Strichartz inequalities, we have the lower frequency estimate
    \[
    \begin{aligned}
        &\quad\left\| \nabla \left(  e^{   t {z} \Delta_\Omega} P^{\Omega}_{<  L } w_0 \right)  \right\|_{L_{t,x}^2  ([0, T] \times \{ |x- x_0  | \le R\} )} \\
        & \lesssim  T^{\frac{2}{ N +2}}R^{\frac{2 N (1+\theta)-4(1-\theta)}{2( N +2)}}
        \left\|\nabla  \left( e^{   t {z} \Delta_\Omega} P^{\Omega}_{<  L } w_0 \right) \right\|_{L_t^{\frac{2( N +2)}{( N -2)}}  L_x^{\frac{2 N ( N +2)}{ N^2+4 - 2 N \theta-4\theta}} }\\
        & \lesssim T^{\frac{2}{ N +2}}R^{\frac{2 N (1+\theta)-4(1-\theta)}{2( N +2)}}  L^{\theta}\left\| \left(-\Delta_\Omega \right)^{(1-\theta)/2} e^{   t {z} \Delta_\Omega} P^{\Omega}_{<  L }w_0 \right\|_{L_t^{\frac{2( N +2)}{( N -2)}}  L_x^{\frac{2 N ( N +2)}{ N^2+4 - 2 N \theta-4\theta}} }\\
        & \lesssim T^{\frac{2}{ N +2}}R^{\frac{2 N (1+\theta)-4(1-\theta)}{2( N +2)}}  L^{\theta} \left\| e^{   t {z} \Delta_\Omega} w_0 \right\|_{L_{t,x}^{\frac{2( N +2)}{ N -2}}}^{\theta} \left\| (-\Delta_\Omega)^{1/2} e^{   t {z} \Delta_\Omega} w_0 \right\|_{L_t^{\frac{2( N +2)}{ N -2}}  L_x^{\frac{2 N ( N +2)}{ N^2+4}} }^{1-\theta} \\
        & \lesssim T^{\frac{2}{ N +2}}R^{\frac{2 N (1+\theta)-4(1-\theta)}{2( N +2)}}  L^{\theta} \left\| e^{   t {z} \Delta_\Omega} w_0 \right\|_{L_{t,x}^{\frac{2( N +2)}{ N -2}}}^{\theta} \|w_0\|_{\dot{H}^1_D} ^{1-\theta}.
    \end{aligned}
    \]
     We now turn to estimate the term of high frequencies. By functional calculus, $e^{t z\Delta_\Omega} = e^{   it \Im{z} \Delta_\Omega}  e^{   t \Re{z} \Delta_\Omega}$. By the local smoothing of NLS on the exterior domain in \cite{KVZ1}, we have
   \[
    \left\| \nabla e^{   it \Im{z} \Delta_\Omega}\left(  e^{   t \Re{z} \Delta_\Omega} P^{\Omega}_{\ge  L } w_0  \right) \right\|_{L_{t,x}^2  ([0, T] \times \{ |x- x_0  | \le R\} )}^2 \lesssim R \| e^{   T \Re{z} \Delta_\Omega} P^{\Omega}_{\ge  L }w_0\|_{L^2( \Omega)} \| e^{   T \Re{z} \Delta_\Omega} P^{\Omega}_{\ge  L }w_0\|_{ \dot{H}_D^1( \Omega)}.
    \]
    Due to the fact that $w_0 \in \dot{H}^1_D$, which enables us to do zero extension for $w_0$ to $\mathbb{R}^N$ so that $e^{t \Delta_\Omega} = e^{t \Delta_{ \mathbb{R}^N } }$, we use the Bernstein inequality to get that
    \[
    \left\|  e^{   T \Re{z} \Delta_\Omega} P^{\Omega}_{\ge  L } w_0 \right\|_{ \dot{H}_D^1  (  \Omega )} \lesssim R^{1/2} L^{-1/2} \|w_0\|_{\dot{H}^1_D}.
    \]
    Choosing $ L  = \left(T^{\frac{2}{ N +2}}R^{\frac{ N (1+2\theta)-2(3-2\theta)}{2( N +2)}} \left\| e^{   t {z} \Delta_\Omega} w_0 \right\|_{L_{t,x}^{\frac{2( N +2)}{ N -2}}}^{\theta} \|w_0\|_{\dot{H}^1_D} ^{-\theta}\right)^{-\frac{2}{1+2\theta}}$ leads to the optimal estimate
    \[
    \begin{aligned}
       \left\| \nabla e^{   t {z} \Delta_\Omega} w_0 \right\|_{L_{t,x}^2  ([0, T] \times \{ |x- x_0  | \le R\} )}\lesssim T^{\frac{2}{(N+2)(1+2\theta)}}R^{\frac{N(1+\theta)-2(1-\theta)+\theta(N+2)}{(N+2)(1+2\theta)}}  \left\| e^{   t {z} \Delta_\Omega} w_0 \right\|_{L_{t,x}^{\frac{2(N+2)}{N-2}}}^{\frac{\theta}{1+2\theta}} \|w_0\|_{\dot{H}^1_D} ^{\frac{1+\theta}{1+2\theta}}.
    \end{aligned}
    \]
\end{proof}

\section{Existence of critical elements}\label{se3}

In this section, we show the existence of a critical element if the global well-posedness of \eqref{eq:gl} does not hold.

Set
\begin{align*}
L(E) : = \sup  \,   \|u\|_{L_{t,x}^\frac{2( N +2)}{N -2} ( I \times \Omega)} ,
\end{align*}
where the supremum is taken over all solutions $u$ to \eqref{eq:gl} defined on the lifespan $I$ {and energy $E(u_0)\leq E$}. If Theorem \ref{th1.3} were false, then there would be a critical energy $E_c\in (0,\infty)$ so that
\begin{align}\label{eq7.2}
\text{ $L(E) < \infty$ for $E(u_0) < E_c$\ and \ $L(E)= \infty$ for $E(u_0) \ge E_c$.}
\end{align}
The positivity of $E_c$ follows from the global well-posedness with small data. Indeed, one can show the stronger statement
\begin{align}\label{eq7.3}
\|u \|_{\dot{X}^1( \mathbb{R}_+ \times \Omega)} \lesssim E(u_0)^\frac12 \quad\text{ for all data with } E(u_0)\le \eta_0,
\end{align}
where $\eta_0$ is a small constant and $\dot{X}^1 : = L_{t,x}^{\frac{2(N +2)}{N  -2 } } \cap {L_t^\frac{N +2}{N - 2}  \dot{W}_x^{1, \frac{2 N ( N +2) }{N^2 - 2 N  + 8 }}}$.

\subsection{Profile decomposition}\label{subse3.1}

Before presenting the existence of a critical element, we prove in this subsection a linear profile decomposition for the propagator $e^{tz \Delta_\Omega}$ with data in $\dot{H}_D^1( \Omega)$,  which is a powerful tool to study the defect of  embedding of  Sobolev spaces \cite{Ge}.

Let $\Theta: \mathbb{R}^N  \to [0, 1]$ be a smooth function such that
\begin{align*}
\Theta(x) =
\begin{cases}
0, |x| \le \frac14 , \\
1, |x| \ge \frac12.
\end{cases}
\end{align*}
Our aim of this subsection  is to prove the following linear profile decomposition in $\dot{H}_D^1(\Omega)$, which is a consequence of Proposition \ref{pr5.2v65} proved later on.

\begin{theorem}[Linear profile decomposition] \label{th5.6v65}
Let $\{f_n\}$ be a bounded sequence in $\dot{H}_D^1$. After passing to a subsequence, there exist $J^* \in  \left\{0, 1, 2, \cdots, \infty \right\}$, $ \left\{ \phi_n^j \right\}_{j = 1}^{J^*} \subseteq \dot{H}_D^1( \Omega)$, $ \left\{\lambda_n^j \right\}_{j = 1}^{J^*} \subseteq (0, \infty)$, and $ \left\{ x_n^j \right\}_{j = 1}^{J^*} \subseteq \Omega$ conforming to one of the following cases for each $j$:

{\upshape Case 1.}
$\lambda_n^j \equiv \lambda_\infty^j$, $x_n^j \to x_\infty^j$ as $n \to \infty$, and there is a $\phi^j \in \dot{H}_D^1(\Omega)$ so that
$\phi_n^j = \phi^j$.
We define
\begin{align*}
\left(g_n^j f \right)(x) =   \frac1{ \left( \lambda_n^j \right)^{ \frac{ N -2}2} }  f \left( \frac{x - x_n^j}{\lambda_n^j} \right) \quad \text{and}\quad
\Omega_n^j : = \left( \lambda_n^j \right)^{-  1} \left( \Omega -  \left\{x_n^j \right\} \right).
\end{align*}

{\upshape Case 2.}
$\lambda_n^j \to \infty$, $- \frac{x_n^j}{ \lambda_n^j} \to x_\infty^j \in \mathbb{R}^N  $ as $n\to \infty$, and there is a $\phi^j \in \dot{H}^1( \mathbb{R}^N )$ so that
\begin{align*}
& \phi_n^j ( x) = g_n^j  \left( \chi_n^j \phi^j \right)(x)\quad \text{ with }
\left(g_n^j f \right)(x) : =  \frac1{  \left( \lambda_n^j \right)^{ \frac{N - 2}2} }  f \left( \frac{x - x_n^j}{ \lambda_n^j} \right),\\
&\Omega_n^j : =  \left( \lambda_n^j \right)^{- 1}  \left( \Omega - \left\{x_n^j  \right\} \right),  \text{ and}\quad
\chi_n^j(x) =  \Theta \left( \frac{d\left( \lambda_n^j x + x_n^j \right) }{\diam \left( \Omega^c \right)} \right).
\end{align*}

{\upshape Case 3.}
$\beta_n^j: = \frac{ d \left(x_n^j \right)}{ \lambda_n^j} \to \infty$ as $n \to \infty$, and there is a $\phi^j \in \dot{H}^1(\mathbb{R}^N  )$, so that
\begin{align*}
& \phi_n^j(x) = g_n^j  \left( \chi_n^j \phi^j \right)(x)\quad \text{ with }
\left(g_n^j f \right)(x) : =   \frac1{  \left( \lambda_n^j \right)^{ \frac{N-2}2}  } f \left( \frac{x- x_n^j}{ \lambda_n^j} \right), \\
&\Omega_n^j : =  \left( \lambda_n^j \right)^{-1 }  \left( \Omega -  \left\{ x_n^j  \right\} \right),
\text{ and}\quad  \chi_n^j (x ) = 1 - \Theta \left( \frac{ %\lambda_n^j
|x|}{  \beta_n^j
d(x_n^j)
} \right).
\end{align*}

{\upshape Case 4.}
$\lambda_n^j \to 0$, $\beta_n^j: = \frac{d \left(x_n^j \right)}{ \lambda_n^j} \to \beta_\infty^j > 0$ as $n \to \infty$, and there is a $\phi^j \in \dot{H}_D^1( \mathbb{R}_+^N  )$ so that
\begin{align*}
& \phi_n^j(x) =  \left( g_n^j  \phi^j \right)(x)\quad \text{ with }\left(g_n^j f \right)(x) : =
\frac1{ \left( \lambda_n^j \right)^{ \frac{N -2}2 }  }  f \left( \frac{  \left(R_n^j \right)^{- 1}  \left( x -  \left(x_n^j \right)^* \right)}{ \lambda_n^j} \right), \\
&\Omega_n^j : =  \left( \lambda_n^j  \right)^{- 1}  \left(R_n^j \right)^{- 1}  \left( \Omega -  \left\{  \left( x_n^j \right)^*  \right\}  \right), \left(x_n^j \right)^* \in \partial \Omega \text{ is defined to satisfy }
d \left(x_n^j \right) = \left| x_n^j -  \left( x_n^j \right)^* \right|,
\intertext{ and }
& R_n^j \in SO( N  ) \text{ satisfies }
R_n^j e_N = \frac{ x_n^j -  \left( x_n^j  \right)^*}{  \left|x_n^j -  \left( x_n^j \right)^* \right|}, \text{ where } e_N  = (0, \cdots, 0, 1).
\end{align*}
Further, for any finite $0 \le J \le J^*$, we have
\begin{align*}
f_n = \sum\limits_{j = 1}^J \phi_n^j + w_n^J,
\end{align*}
with $w_n^J \in \dot{H}_D^1 (\Omega)$ satisfying
\begin{align}
& \lim\limits_{J \to J^*} \limsup\limits_{n\to \infty}  \left\| e^{tz \Delta_\Omega}  w_n^J  \right\|_{L_{t,x}^{ \frac{2 ( N +2)  }{ N - 2} }(  \mathbb{R}_+ \times  \Omega)} = 0, \label{eq5.25v65new} \\
& \lim\limits_{n \to \infty}  \left( \|f_n \|_{ \dot{H}_D^1(\Omega)}^2 - \sum\limits_{j = 1}^J  \left\|\phi_n^j \right\|_{\dot{H}_D^1(\Omega)}^2 -  \left\|w_n^J \right\|_{\dot{H}_D^1( \Omega)}^2 \right) = 0,\label{eq5.26v65new} \\
& \lim\limits_{n \to \infty}  \left( \|f_n \|_{L_x^\frac{2 N }{ N -2}  ( \Omega)}^\frac{2 N }{N -2}  - \sum\limits_{j = 1}^J  \left\|\phi_n^j \right\|_{L_x^\frac{2 N }{ N -2}}^\frac{2 N }{N -2}  -  \left\|w_n^J  \right\|_{L_x^\frac{2 N }{ N -2} }^\frac{2 N }{ N -2}   \right) = 0,\label{eq5.27v65new}  \\
&  \left(g_n^J \right)^{- 1} w_n^J \rightharpoonup 0 \quad\text{in } \dot{H}^1 \text{ as } n \to \infty ,\label{eq5.28v65new}
\end{align}
and for all $j \ne k$, we have the asymptotic orthogonality property
\begin{align}\label{eq5.29v65new}
\lim\limits_{n \to \infty}  \left( \frac{ \lambda_n^j }{\lambda_n^k} + \frac{\lambda_n^k}{\lambda_n^j} + \frac{  \left|x_n^j - x_n^k \right|^2 }{ \lambda_n^j \lambda_n^k} \right)  = \infty.
\end{align}

\end{theorem}

The profiles can live in different limiting geometries, which is one of the principal differences relative to previous analyses. We now turn to the detailed proof. Throughout this section, we write
\begin{align*}
G_{ \mathcal{O}}(x,y;  \lambda ): =  \left( - \Delta_{\mathcal{O}} -  \lambda   \right)^{-1} (x,y )
\end{align*}
for the Green's function of the Dirichlet Laplacian in a general open set $\mathcal{O}$. This function is symmetric under the interchange of $x$ and $y$.

\begin{definition}[\cite{KVZ1}] \label{de3.1v65}

Given a sequence $\{\mathcal{O}_n\}_n$ of open subsets of $\mathbb{R}^N $, we define
\begin{align*}
\widetilde{\lim} \,  \mathcal{O}_n :  =  \left\{ x \in \mathbb{R}^N : \liminf\limits_{n \to \infty } \dist  \left(x, \mathcal{O}_n^c  \right) > 0  \right\}.
\end{align*}
Writing $\tilde{O} = \widetilde{\lim} \, \mathcal{O}_n$, we say $\mathcal{O}_n \to \mathcal{O}$ if the following two conditions hold: $\mathcal{O} \triangle \tilde{O}$ is a finite set and
\begin{align}\label{eq3.1v65}
G_{\mathcal{O}_n}(x,y;  \lambda ) \to G_{\mathcal{O}} (x,y ;  \lambda )
\end{align}
for all $ \lambda  \in (-2 , - 1)$, all $x \in \mathcal{O}$, and uniformly for $y$ in compact subsets of $\mathcal{O} \setminus \{x \}$.

\end{definition}

Given sequences of scaling and translation parameters $N_n \in 2^{\mathbb{Z}}$ and $x_n \in \Omega$, we would like to consider the domains $N_n \left( \Omega -  \left\{x_n \right\} \right)$.
In \cite{KVZ1}, R. Killip, M. Visan, and X. Zhang prove that in Cases 2 and 3, $\Omega_n \to \mathbb{R}^N $, and in Case 4, $\Omega_n \to \mathbb{R}_+^N $. In fact, they also have proved
\begin{proposition}\label{pr3.6v65}
Assume $\Omega_n \to \Omega_\infty$ as $n \to \infty$ in the sense of Definition \ref{de3.1v65} and let $\Theta \in C_c^\infty ((0, \infty))$.
Then
\begin{align}\label{eq3.11v65}
\left\|  \left( \Theta  \left( - \Delta_{\Omega_n}  \right) - \Theta  \left( - \Delta_{\Omega_\infty}  \right) \right) \delta_y  \right\|_{\dot{H}^{-1} ( \mathbb{R}^N )} \to 0 \text{ as } n\to \infty,
\end{align}
uniformly for $y$ in compact subsets of $\widetilde{\lim } \, \Omega_n$.
\end{proposition}

As a preliminary for our proof, we use the following lemma to establish an important analytic tool, called "inverse" Sobolev inequality.

\begin{lemma}\label{le5.1v65}
Let $f \in \dot{H}_D^1 ( \Omega)$. Then we have
\begin{align}\label{eq4.16v65}
\|f\|_{L_x^\frac{2N }{ N -2}} \lesssim \|f \|_{\dot{H}_D^1 ( \Omega) }^\frac{ N -2}N   \cdot \sup\limits_{ L  \in 2^{\mathbb{Z}}}  \left\|P^\Omega_L f   \right\|_{L_x^\frac{2N }{N -2}}^\frac2N .
\end{align}

\end{lemma}

\begin{proof}

If $N=3$, by Lemma \ref{lea.7v65} and Lemma \ref{lea.6v65}, we have
\begin{align*}
\| f\|_{L_x^6}^6
& \lesssim \int  \left(\sum\limits_K  \left| P_K^\Omega f  \right|^2 \right)  \left( \sum\limits_M  \left| P_M^\Omega f  \right|^2 \right)  \left( \sum\limits_L
\left| P_L^\Omega f  \right|^2  \right) \,\mathrm{d}x
\\
& \lesssim \sum\limits_{K \le M \le L }  \left\| P_K^\Omega f \right\|_{L_x^6}  \left\| P_K^\Omega f  \right\|_{L_x^\infty}
\left\| P_M^\Omega f  \right\|_{L_x^6}^2  \left\| P_L^\Omega f  \right\|_{L_x^3}  \left\| P_L^\Omega f  \right\|_{L_x^6} \\
& \lesssim  \left( \sup\limits_{L \in 2^{\mathbb{Z}}}
\left\| P_L^\Omega f  \right\|_{L_x^6}^4 \right) \sum\limits_{K \le M \le L } K^\frac32 L^\frac12  \left\| P_K^\Omega f  \right\|_{L_x^2}
\left\| P_L^\Omega f  \right\|_{L_x^2} \\
& \lesssim  \left( \sup\limits_{L \in 2^{\mathbb{Z}}}  \left\| P_L^\Omega f  \right\|_{L_x^6}^4 \right) \sum\limits_{K \le M \le L } K^\frac12 L^{- \frac12}
\left\| \nabla P_K^\Omega f  \right\|_{L_x^2}  \left\|\nabla P_L^\Omega f  \right\|_{L_x^2},
\end{align*}
which leads to \eqref{eq4.16v65} via Schur's test and other elementary considerations.

If $N  \ge 4$, then one may modify the argument as follows: %By Lemma \ref{lea.7v65} and Lemma \ref{lea.6v65}, we have
\begin{align*}
\|f \|_{L_x^\frac{2 N }{ N -2}}^\frac{2 N }{ N -2}
& \lesssim \int  \left( \sum\limits_M  \left|P^\Omega _M f \right|^2 \right)^\frac{ N }{2( N -2)}   \left( \sum\limits_L   \left| P^\Omega_L  f \right|^2 \right)^\frac{ N }{2( N -2)} \,\mathrm{d}x
\\
& \lesssim \sum\limits_{M \le L } \int  \left|P^\Omega_M f \right|^\frac{ N }{ N -2}  \left|P^\Omega_L  f \right|^\frac{ N }{ N -2} \,\mathrm{d}x
\\
%& \lesssim  \left( \sup\limits_{K \in 2^{\mathbb{Z}}}  \left\|P_K^\Omega f   \right\|_{L^\frac{2 N }{ N -2}}  \right)^\frac4{ N -2} \sum\limits_{M \le L }
%\left\| P_M^\Omega f \right\|_{L_x^\frac{2 N }{ N -4}}  \left\| P_L^\Omega f \right\|_{L^2} \\
& \lesssim  \left( \sup\limits_{K\in 2^{\mathbb{Z}}}  \left\| P_K^\Omega f \right\|_{L_x^\frac{2 N }{ N -2}} \right)^\frac4{ N -2} \sum\limits_{M \le  L } M^{-1}
L^{-1} \left\|\nabla  P_M^\Omega f  \right\|_{L_x^\frac{2 N }{ N -4}}  \left\|\nabla P_L^\Omega f  \right\|_{L^2_x} \\
%& \lesssim  \left( \sup\limits_{K \in 2^{\mathbb{Z}}}  \left\| P_K^\Omega f  \right\|_{L_x^\frac{2 N }{ N -2}} \right)^\frac4{ N -2} \sum\limits_{M \le L }
%M L^{-1}  \left\|\nabla P_M^\Omega f  \right\|_{L_x^2}  \left\|\nabla P_L^\Omega f  \right\|_{L_x^2} \\
& \lesssim  \left( \sup\limits_{K \in 2^{\mathbb{Z}}}   \left\| P_K^\Omega f \right\|_{L_x^\frac{2 N }{ N -2}}  \right)^\frac4{N-2}  \left( \sum\limits_{ K \in 2^{\mathbb{Z}}}  \left\|\nabla P_K^\Omega f \right\|_{L_x^2}^2  \right).
\end{align*}

\end{proof}

\begin{proposition}[Inverse Sobolev inequality] \label{pr5.2v65}
Let $\{f_n \} \subseteq \dot{H}_D^1( \Omega)$. Suppose that
\begin{align*}
\lim\limits_{n \to \infty} \|f_n \|_{\dot{H}_D^1 ( \Omega)} = A < \infty\quad
\text{ and }\quad
\lim\limits_{n \to \infty} \|f_n \|_{L_x^\frac{2 N }{ N -2} ( \Omega)} = \epsilon > 0.
\end{align*}
Then there exist a subsequence in $n$, $\{\phi_n \} \subseteq \dot{H}_D^1 ( \Omega)$, $ \{N_n\} \subseteq 2^{\mathbb{Z}}$, $\{x_n \} \subseteq \Omega$ conforming to one of the following four cases:
\begin{enumerate}
	\item
	Case 1. $N_n \equiv N_\infty \in 2^{\mathbb{Z}}$ and $x_n \to x_\infty \in \Omega$ as $n \to \infty$. In this case, we choose $\phi \in \dot{H}_D^1( \Omega)$ and the subsequence so that $f_n  \rightharpoonup \phi $ in $\dot{H}_D^1( \Omega)$ as $n \to \infty$, and then set $\phi_n : = \phi.$
	
	\item Case 2. $N_n \to 0$ and $-N_n x_n \to x_\infty \in \mathbb{R}^N $ as $n \to \infty$.
	In this case, we choose $\tilde{\phi} \in \dot{H}^1( \mathbb{R}^N )$ and the subsequence so that
	$g_n(x) : = N_n^{- \frac{N-2}2} f_n  \left( N_n^{-1} x + x_n \right)  \rightharpoonup \tilde{\phi} (x) $ in $\dot{H}^1( \mathbb{R}^N )$ as $n \to \infty$,
	and then set
	\begin{align*}
		\phi_n (x) : = N_n^\frac{N-2}2  \left( \chi_n \tilde{\phi}  \right) (N_n(x- x_n)),
	\end{align*}
	where $\chi_n (x) = \Theta \left( \frac{d \left( N_n^{-1} x+ x_n  \right) }{ \diam  \left( \Omega^c \right)}  \right)$.
	
	\item Case 3. $N_n d(x_n) \to \infty$ as $n \to \infty$. In this case, we choose $\tilde{\phi} \in \dot{H}^1( \mathbb{R}^N )$ and the subsequence so that
	\begin{align*}
		g_n(x) : = N_n^{- \frac{N-2}2} f_n \left(N_n^{-1} x + x_n \right) \rightharpoonup \tilde{\phi} (x) \text{ in } \dot{H}^1(\mathbb{R}^N ) \text{ as } n \to \infty,
	\end{align*}
	and then set
	\begin{align*}
		\phi_n(x) : = N_n^\frac{N-2}2  \left( \chi_n \tilde{\phi} \right) ( N_n ( x- x_n)),
		\text{ where }
		\chi_n(x) = 1 - \Theta  \left( \frac{|x|}{ N_n d(x_n)}  \right).
	\end{align*}

	\item Case 4. $N_n \to \infty$ and $N_n d(x_n) \to d_\infty > 0$ as $n \to \infty$. In this case, we choose $\tilde{ \phi} \in \dot{H}_D^1( \mathbb{R}_+^N )$ and the subsequence so that
	\begin{align*}
		g_n(x): = N_n^{- \frac{N-2}2} f_n \left(N_n^{-1} R_n x + x_n^* \right)  \rightharpoonup \tilde{\phi}(x) \text{ in } \dot{H}^1( \mathbb{R}^N )
		\text{ as } n \to \infty,
	\end{align*}
	and then set
	\begin{align*}
		\phi_n(x) : = N_n^\frac{N-2}2 \tilde{ \phi} \left( N_n R_n^{-1}  \left(x- x_n^* \right) \right),
	\end{align*}
	where $R_n \in SO(N )$ satisfies
	\begin{align*}
		R_n e_N  = \frac{x_n - x_n^*}{  \left|x_n- x_n^*  \right|}
	\end{align*}
	and $x_n^* \in \partial \Omega$ is chosen so that $d(x_n) =  \left|x_n - x_n^* \right|$.
\end{enumerate}
Furtheremore, the following estimates hold:
\begin{align}
\liminf\limits_{n \to \infty} \|\phi_n \|_{\dot{H}_D^1}^2 &\gtrsim A^2  \left( \frac{\epsilon}A \right)^\frac{ N^2}2, \label{eq5.1v65}
\\
\liminf\limits_{n \to \infty}  \left( \|f_n \|_{\dot{H}_D^1( \Omega)}^2 - \|f_n - \phi_n \|_{\dot{H}_D^1}^2  \right) &\gtrsim A^2  \left( \frac\epsilon{A} \right)^\frac{ N^2}2, \label{eq5.2v65}\\
\limsup\limits_{n\to \infty} \|f_n(x) - \phi_n(x) \|_{L_x^\frac{2 N }{ N -2}}^\frac{2 N }{ N -2} &\le \epsilon^\frac{2 N}{ N -2}  \left( 1 - c  \left( \frac\epsilon{A} \right)^\frac{ N ( N +2)}2  \right). \label{eq5.3v65}
\end{align}
\end{proposition}

\begin{proof}
By Lemma \ref{le5.1v65} and the assumptions on $f_n$, for each $n$, there exists $N_n \in 2^{\mathbb{Z}}$ such that
\begin{align*}
\liminf\limits_{n \to \infty} \left\|P_{N_n}^\Omega  f_n  \right\|_{L_x^\frac{2 N }{ N -2}} \gtrsim \epsilon^\frac{ N }2 A^{- \frac{ N -2}2}.
\end{align*}
Set $\lambda_n = N_n^{-1}$. By the H\"older's inequality, we have
\begin{align*}
\epsilon^\frac{ N }2 A^{- \frac{ N -2}2}
& \lesssim \liminf\limits_{n \to \infty}  \left\|P_{N_n}^\Omega f_n  \right\|_{L_x^\frac{2 N }{ N -2}}
\lesssim \liminf\limits_{n \to \infty }  \left\|P_{N_n}^\Omega f_n  \right\|_{L_x^2}^\frac{ N -2} N   \left\|P_{N_n}^\Omega f_n  \right\|_{L_x^\infty}^\frac2N
\\
& \lesssim \liminf\limits_{n \to \infty}  \left(A N_n^{-1}  \right)^\frac{ N -2}N   \left\|P_{N_n}^\Omega f_n  \right\|_{L_x^\infty}^\frac2N .
\end{align*}
Thus we may find $x_n \in \mathbb{R}^N $ such that
\begin{align*}
\liminf\limits_{n \to \infty}  \left|  \left( P_{N_n}^\Omega f_n \right)( x_n) \right| N_n^{- \frac{N -2}2} \gtrsim \epsilon^{\frac{N^2}4} A^{1- \frac{N^2}4} .
\end{align*}
This leads to
\begin{align}\label{eq5.4v65}
 \left| \left(P_{N_n}^\Omega f_n \right)(x_n) \right| \gtrsim N_n^\frac{N -2}2 \epsilon^\frac{N^2}4 A^{1- \frac{N^2}4}.
\end{align}
Note that the different cases in Proposition \ref{pr5.2v65} are completely determined by the behavior of $x_n$ and $N_n$. We next observe that
\begin{align}\label{eq5.5v65}
N_n d(x_n) \gtrsim  \left( \frac\epsilon{A} \right)^\frac{N^2}4 \quad\text{ whenever } N_n \gtrsim 1.
\end{align}
Indeed, Theorem \ref{th2.2} implies that whenever $N_n \gtrsim 1$, it holds
\begin{align*}
\int_\Omega  \left| p_{N_n^{-2}}(x_n, y)  \right|^2 \,\mathrm{d}y
& \lesssim N_n^{2N}  \int_\Omega  \left| (N_n d(x_n))   \left(N_n d(x_n) + N_n  |x_n - y|  \right) e^{- c N_n^2 |x_n - y |^2 }  \right|^2 \,\mathrm{d} y
\\
& \lesssim  (N_n d(x_n ))^2  \left(N_n d(x_n) + 1 \right)^2 N_n^N .
\end{align*}
Write
\begin{align*}
\left(P_{N_n}^\Omega f_n  \right)(x_n) = \int_\Omega p_{N_n^{-2}}(x_n, y)
 \left(  P_{\le 2 N_n}^\Omega e^{- \frac1{N_n^2} \Delta_\Omega} P_{N_n}^\Omega f_n  \right)(y) \,\mathrm{d}y.
\end{align*}
Then by \eqref{eq5.4v65} and the Cauchy-Schwarz inequality, we have
\begin{align*}
  N_n^\frac{N-2}2 \epsilon^\frac{N^2}4 A^{1 - \frac{N^2}4}
& \lesssim (N_n d(x_n)) (N_n d(x_n) + 1) N_n^\frac{N}2  \left\|P_{ \le 2N_n}^\Omega e^{- \frac1{N_n^2} \Delta_\Omega} P_{N_n}^\Omega f_n  \right\|_{L_x^2}
\\
& \lesssim (N_n d(x_n)) (N_n d(x_n) + 1) N_n^\frac{N-2}2 \|f_n \|_{\dot{H}_D^1( \Omega)},
\end{align*}
from which \eqref{eq5.5v65} follows.

Due to \eqref{eq5.5v65}, upon taking a subsequence, we only need to consider the following four cases:
\begin{enumerate}
\item $N_n \sim 1$ and $N_n d(x_n) \sim 1$,

\item $N_n \to 0$ and $N_n d(x_n) \lesssim 1$,

\item $N_n d(x_n) \to \infty$ as $n \to \infty$,

\item $N_n \to \infty$ and $N_n d(x_n) \sim 1$.

\end{enumerate}

We next deal with these four cases in order.

\textbf{Case 1}. Up to a subsequence, we may assume
\begin{align*}
N_n \equiv N_\infty \in 2^{\mathbb{Z}} \text{ and } x_n \to x_\infty \in \Omega, \text{ as } n \to \infty.
\end{align*}
For the proof of later cases, we introduce
\begin{align*}
g_n (x ): = N_n^{- \frac{N-2} 2} f_n  \left(N_n^{-1} x + x_n \right).
\end{align*}
As $f_n$ is supported in $\Omega$, $g_n$ is supported in $\Omega_n : = N_n ( \Omega - \{x_n\})$. Moreover, we have
\begin{align*}
\|g_n \|_{\dot{H}_D^1( \Omega_n)} = \|f_n \|_{\dot{H}_D^1(  \Omega)} \lesssim A.
\end{align*}
Passing to a further subsequence, we may find $\tilde{\phi}$ so that $g_n \rightharpoonup \tilde{\phi}$ in $\dot{H}^1( \mathbb{R}^N )$ as $n \to \infty$.
Rescaling this weak convergence leads to
\begin{align}\label{eq5.6v65}
f_n(x) \rightharpoonup \phi(x) : = N_\infty^\frac{N-2} 2 \tilde{\phi} (N_\infty (x- x_\infty) ) \text{ in } \dot{H}_D^1(\Omega).
\end{align}
Since $\dot{H}_D^1( \Omega)$ is a weakly closed subset of $\dot{H}^1( \mathbb{R}^N )$, $\phi \in \dot{H}_D^1( \Omega)$.

We next show that $\phi$ is non-trivial. To this end, let $h: = P_{N_\infty}^\Omega \delta_{x_\infty}$. Then we obtain from the Bernstein inequality
\begin{align}\label{eq5.7v65}
\left\|   \left(- \Delta_\Omega \right)^{- \frac12} h  \right\|_{L^2(\Omega)}  =  \left\|  \left( - \Delta_\Omega \right)^{- \frac12} P_{N_\infty}^\Omega \delta_{x_\infty}  \right\|_{L^2( \Omega)}
\lesssim N_\infty^\frac{N-2}2,
\end{align}
which implies that $h \in \dot{H}_D^{-1} (\Omega)$. On the other hand, we note
\begin{align}\label{eq5.8v65}
\langle \phi, h \rangle = \lim\limits_{n \to \infty} \langle f_n , h \rangle
= \lim\limits_{n \to \infty}  \left\langle f_n, P_{N_\infty}^\Omega \delta_{x_\infty}  \right\rangle
= \lim\limits_{n \to \infty}  \left(P_{N_n}^\Omega f_n \right)(x_n) + \lim\limits_{n \to \infty}
\left\langle f_n, P_{N_\infty}^\Omega  \left( \delta_{x_\infty} - \delta_{x_n } \right)  \right\rangle.
\end{align}
Applying the basic elliptic theory
\begin{align}\label{eq5.9v65}
\| \nabla v \|_{L^\infty (|x| \le R)}
\lesssim R^{-1} \|v\|_{L^\infty ( |x| \le 2 R )} + R \| \Delta v \|_{L^\infty ( |x| \le 2 R)},
\end{align}
to $v(x) =  \left(P_{N_\infty}^\Omega f_n \right)(x+x_n)$ with $R= \frac12 d(x_n)$, we see that the second limit in \eqref{eq5.8v65} vanishes.

By the Bernstein inequalities and the fact that $d(x_n) \sim 1$, we have
\begin{align*}
 \left\| P_{N_\infty}^\Omega f_n  \right\|_{L_x^\infty} \lesssim N_\infty^\frac{N-2} 2 A\quad
\text{ and }\quad
 \left\|\Delta P_{N_\infty}^\Omega f_n  \right\|_{L_x^\infty} \lesssim N_\infty^\frac{N+2}2 A.
\end{align*}
Thus for $n$ sufficiently large, it follows from the fundamental theorem of calculus and \eqref{eq5.9v65} that
\begin{align}\label{eq5.10v65}
 \left| \left\langle f_n, P_{N_\infty}^\Omega  \left(  \delta_{x_\infty} - \delta_{x_n}  \right)  \right\rangle  \right|
\lesssim |x_\infty - x_n |  \left\| \nabla P_{N_\infty}^\Omega f_n  \right\|_{L^\infty ( |x| \le R)}
\lesssim A  \left( \frac{ N_\infty^\frac{N-2}2 }{d(x_n)} + N_\infty^\frac{N+2}2 d(x_n) \right) |x_\infty - x_n |,
\end{align}
which converges to zero as $n \to \infty$.

Therefore, we obtain from \eqref{eq5.4v65}, \eqref{eq5.7v65}, \eqref{eq5.8v65}, and \eqref{eq5.10v65} that
\begin{align}\label{eq5.11v65}
N_\infty^\frac{N-2}2 \epsilon^\frac{N^2}4 A^{1 - \frac{N^2}4}  \lesssim |\langle \phi, h \rangle |
\lesssim \|\phi \|_{\dot{H}_D^1( \Omega)} \|h \|_{\dot{H}_D^{-1} ( \Omega)}
\lesssim N_\infty^\frac{N-2}2 \|\phi \|_{\dot{H}_D^1( \Omega)},
\end{align}
which gives \eqref{eq5.1v65}.

Since $\dot{H}_D^1( \Omega)$ is a Hilbert space, \eqref{eq5.2v65} follows immediately from \eqref{eq5.1v65} and \eqref{eq5.6v65}.

It remains to prove decoupling for the $L_x^\frac{2 N }{ N -2}$ norm. Note first that since $f_n$ is bounded in $\dot{H}_D^1( \Omega)$, we may pass to a further subsequence so that
$f_n \to \phi$ in $L^2$-sense on any compact set via the Rellich-Kondrashov Theorem. Passing to yet another subsequence, we may assume furtheremore that $f_n \to \phi$ almost everywhere in $\Omega$.
Using the refined Fatou Lemma, and a change of variables, we find
\begin{align*}
\lim\limits_{n \to \infty}  \left( \|f_n \|_{L_x^\frac{2 N }{ N -2}}^\frac{2 N }{ N -2} - \|f_n - \phi_n \|_{L_x^\frac{2 N }{ N -2}}^\frac{2 N }{ N -2}  \right) = \|\phi \|_{L_x^\frac{2 N }{ N -2}( \Omega)}^\frac{2 N }{ N -2}.
\end{align*}
Then \eqref{eq5.3v65} would follow if we were able to prove
\begin{align}\label{eq5.12v65}
\| \phi \|_{L_x^\frac{2 N }{ N -2} ( \Omega)} \gtrsim \epsilon^\frac{N^2}4 A^{1- \frac{N^2}4}.
\end{align}
By the Mikhlin multiplier theorem, we have
\begin{align*}
\left\| P_{\le 2 N_\infty}^\Omega  \right\|_{L_x^\frac{2 N }{ N +2} \to L_x^\frac{2 N }{ N +2}} \lesssim 1.
\end{align*}
Using this and the Bernstein's inequality, we get
\begin{align*}
\| h \|_{L_x^\frac{2 N }{ N  +2}} \lesssim   \left\|P_{\le 2 N_\infty}^\Omega  \right\|_{L_x^\frac{2 N }{ N +2} \to L_x^\frac{2 N }{ N +2}}  \left\|P_{N_\infty}^\Omega  \right\|_{L_x^1 \to L_x^\frac{2 N }{ N +2}} \| \delta_{x_\infty} \|_{L_x^1}
\lesssim N_\infty^\frac{N-2}{2}.
\end{align*}
This together with \eqref{eq5.11v65} yields
\begin{align*}
N_\infty^\frac{N-2}2 \epsilon^\frac{N^2}4 A^{1- \frac{N^2}4}  \lesssim  |\langle \phi, h \rangle |
\lesssim N_\infty^\frac{N-2}{2} \|\phi \|_{L_x^\frac{2 N }{ N -2}},
\end{align*}
which justifies \eqref{eq5.12v65}.
\medskip

\textbf{Case 2}. As $N_n \to 0$, the condition $N_n d(x_n) \lesssim 1$ indicates that $\{ N_n x_n \}_{n \ge 1}$ is bounded. Thus, up to a subsequence, we may assume $-N_n x_n \xrightarrow{n \to \infty} x_\infty\in \mathbb{R}^N $.
As in Case 1, we set $\Omega_n : = N_n ( \Omega - \{x_n \})$. Note that the rescaled obstacles $\Omega_n^c$ shrink to $x_\infty$ as $n \to \infty$.

Since $f_n$ is bounded in $\dot{H}_D^1( \Omega)$, the sequence $g_n$ is bounded in $\dot{H}_D^1( \Omega_n ) \subseteq \dot{H}^1( \mathbb{R}^N )$.
Thus, up to a subsequence, we may find $\tilde{\phi}$ so that $g_n \rightharpoonup \tilde{\phi}$ in $\dot{H}^1(\mathbb{R}^N )$ as $n \to \infty$.

Next we show
\begin{align}\label{eq5.13v65}
\chi_n \tilde{\phi} \to \tilde{\phi},
\text{ or equivalently, }
 \left(1 - \chi \left( N_n^{-1} x + x_n  \right) \right) \tilde{\phi}(x) \to 0 \text{ in } \dot{H}^1( \mathbb{R}^N ).
\end{align}
Towards the proof of \eqref{eq5.13v65}, we first set
\begin{align*}
B_n : =  \left\{ x \in \mathbb{R}^N : \dist  \left(x, \Omega_n^c  \right) \le \diam  \left( \Omega_n^c \right)  \right\},
\end{align*}
which contains $\supp \, ( 1 - \chi_n)$ and $\supp \, ( \nabla \chi_n)$.
Since $N_n\to 0$, the measure of $B_n$ shrinks to zero as $n \to \infty$.
By the H\"older and Sobolev inequalities, we have
\begin{align*}
 \left\|  \left( 1 - \chi \left(N_n^{-1} x + x_n  \right) \right) \tilde{\phi} (x)  \right\|_{\dot{H}^1( \mathbb{R}^N )}
 \lesssim  \left\|\nabla \tilde{\phi}  \right\|_{L^2( B_n)} +  \left\|\tilde{\phi}  \right\|_{L^\frac{2N}{N-2} ( B_n)},
\end{align*}
which converges to zero by the dominated convergence theorem.

With \eqref{eq5.13v65} at hand the proofs of \eqref{eq5.1v65} and \eqref{eq5.2v65} now follow closely their Case 1 counterparts.
Indeed, set $h : = P_1^{\mathbb{R}^N } \delta_0$. Then
\begin{align*}
\left\langle \tilde{\phi}, h  \right\rangle = \lim\limits_{n \to \infty} \langle g_n , h \rangle
= \lim\limits_{n \to \infty}  \left\langle g_n , P_1^{\Omega_n} \delta_0  \right\rangle + \lim\limits_{n \to \infty} \left\langle g_n,  \left(P_1^{\mathbb{R}^N } - P_1^{\Omega_n}  \right) \delta_0  \right\rangle .
\end{align*}
By Proposition \ref{pr3.6v65} and the uniform boundedness of $\|g_n \|_{ \dot{H}^1( \mathbb{R}^N )}$, the second term vanishes. Therefore,
\begin{align}\label{eq5.14v65}
 \left|  \left\langle \tilde{\phi}, h  \right\rangle  \right|
& = \left| \lim\limits_{n \to \infty}  \left\langle g_n, P_1^{\Omega_n } \delta_0  \right\rangle \right|
 =  \left| \lim\limits_{n \to \infty}  \left\langle f_n, N_n^\frac{N+2}2  \left(P_1^{\Omega_n } \delta_0  \right) (N_n (x- x_n ))  \right\rangle  \right| \notag
\\
& =  \left|\lim\limits_{n \to \infty}  \left\langle f_n, N_n^{- \frac{N-2}2} P_{N_n}^\Omega \delta_{x_n}  \right\rangle  \right|
\gtrsim \epsilon^\frac{N^2}4 A^{1 - \frac{N^2}4},
\end{align}
where the last inequality follows from \eqref{eq5.4v65}. Thus, as in \eqref{eq5.11v65}, we get
\begin{align*}
\left\| \tilde{\phi}  \right\|_{\dot{H}^1( \mathbb{R}^N )} \gtrsim \epsilon^\frac{N^2}4 A^{1 - \frac{N^2}4}.
\end{align*}
Combining this with \eqref{eq5.13v65}, for $n$ sufficiently large, we obtain
\begin{align*}
\| \phi_n \|_{\dot{H}_D^1( \Omega)} =  \left\| \chi_n \tilde{\phi}  \right\|_{\dot{H}_D^1 ( \Omega_n )} \gtrsim
\epsilon^\frac{N^2}4 A^{1 - \frac{N^2}4},
\end{align*}
which gives \eqref{eq5.1v65}.

To prove the decoupling in $\dot{H}_D^1 ( \Omega)$, we write
\begin{align*}
&\quad \| f_n \|_{\dot{H}_D^1( \Omega)}^2 - \|f_n - \phi_n \|_{\dot{H}_D^1( \Omega)}^2
\\
& = 2 \langle f_n, \phi_n \rangle_{\dot{H}_D^1 ( \Omega)} - \|\phi_n \|_{\dot{H}_D^1 ( \Omega)}^2
= 2  \left\langle N_n^{- \frac{N-2}2} f_n (N_n^{-1} x + x_n ), \tilde{\phi}(x) \chi_x(x)  \right\rangle_{\dot{H}_D^1( \Omega_n )} -  \left\|\chi_n \tilde{\phi} \right\|_{\dot{H}_D^1( \Omega_n )}^2
\\
& = 2  \left\langle g_n, \tilde{\phi}  \right\rangle_{\dot{H}^1( \mathbb{R}^N )} - 2  \left\langle g_n , \tilde{\phi} ( 1 - \chi_n )  \right\rangle_{\dot{H}^1( \mathbb{R}^N )} -  \left\| \chi_n \tilde{\phi}  \right\|_{\dot{H}_D^1 ( \Omega_n)}^2 .
\end{align*}
From the weak convergence of $g_n$ to $\tilde{\phi}$, \eqref{eq5.13v65}, and \eqref{eq5.1v65}, we conclude
\begin{align*}
\lim\limits_{n \to \infty}  \left( \|f_n \|_{\dot{H}_D^1 ( \Omega)}^2 - \|f_n - \phi_n \|_{\dot{H}_D^1 ( \Omega)}^2  \right)
=  \left\| \tilde{\phi}  \right\|_{\dot{H}^1( \mathbb{R}^N )}^2 \gtrsim \epsilon^\frac{N^2}2 A^{2 - \frac{N^2}2 } .
\end{align*}
This verifies \eqref{eq5.2v65}.

We now prove the decoupling for $L_x^\frac{2 N }{N -2} ( \Omega)$ norm by showing
\begin{align}\label{eq5.15v65}
\liminf\limits_{n \to \infty}  \left( \|f_n \|_{L_x^\frac{2N }{N -2}}^\frac{2N }{ N -2} - \|f_n - \phi_n \|_{L_x^\frac{2 N }{ N -2}}^\frac{2 N }{ N -2}  \right)
=  \left\|\tilde{\phi}  \right\|_{L_x^\frac{2 N }{ N -2}}^\frac{2 N }{ N -2}.
\end{align}
Notice that \eqref{eq5.3v65} then follows from the lower bound
\begin{align}\label{eq5.16v65}
 \left\| \tilde{\phi}  \right\|_{L_x^\frac{2 N }{ N -2}}^\frac{2 N }{ N -2} \gtrsim \left( \epsilon^\frac{N^2}4 A^{1 - \frac{N^2}4} \right)^\frac{2N}{N-2}.
\end{align}
We prove \eqref{eq5.16v65} in a similar way as in Case 1: by \eqref{eq5.14v65} and the Mikhlin multiplier theorem
\begin{align*}
\epsilon^\frac{N^2}4 A^{1 - \frac{N^2}4}  \lesssim  \left|  \left\langle \tilde{\phi} , h  \right\rangle  \right|
\lesssim  \left\| \tilde{\phi}  \right\|_{L_x^\frac{2 N }{ N -2}} \| h \|_{L_x^\frac{2 N }{ N +2} }
\sim  \left\|\tilde{\phi}  \right\|_{L_x^\frac{2 N }{ N -2}}  \left\|P_1^{\mathbb{R}^N } \delta_0  \right\|_{L_x^\frac{2N }{ N +2}}
\lesssim  \left\|\tilde{\phi}  \right\|_{L_x^\frac{2 N }{ N -2}}.
\end{align*}
Plugging this into \eqref{eq5.15v65} leads to \eqref{eq5.3v65} in Case 2.

To establish \eqref{eq5.15v65}, we need two simple observations: the first one is
\begin{align}\label{eq5.17v65}
g_n - \chi_n \tilde{\phi} \to 0 \quad\text{ almost everywhere in } \mathbb{R}^N ,
\end{align}
while the second one is
\begin{align}\label{eq5.18v65}
 \left\| \chi_n \tilde{\phi} - \tilde{\phi}  \right\|_{L_x^\frac{2 N }{ N -2} ( \mathbb{R}^N )} \to 0.
\end{align}
For \eqref{eq5.17v65},  using the definition of $\tilde{\phi}$ together with \eqref{eq5.13v65}, we deduce
\begin{align*}
g_n - \chi_n \tilde{\phi} \rightharpoonup 0 \quad\text{ in } \dot{H}^1( \mathbb{R}^N ).
\end{align*}
Thus, by  a similar argument as in Case 1, and passing to a subsequence, we obtain \eqref{eq5.17v65}. For \eqref{eq5.18v65}, we simply use \eqref{eq5.13v65} and the Sobolev inequality.

Combining \eqref{eq5.17v65} with \eqref{eq5.18v65}, and passing to a subsequence if necessary, we obtain
\begin{align*}
g_n - \tilde{\phi} \to 0 \quad\text{ almost everywhere
 in } \mathbb{R}^N,
\end{align*}
which yields by the refined Fatou lemma
\begin{align*}
\liminf\limits_{n \to \infty}  \left( \|g_n \|_{L_x^\frac{2 N }{ N -2}}^\frac{2 N }{ N -2} -  \left\|g_n - \tilde{\phi}  \right\|_{L_x^\frac{2 N }{ N -2}}^\frac{2 N }{ N -2}  \right)
=  \left\|\tilde{\phi}  \right\|_{L_x^\frac{2 N }{N -2}}^\frac{2 N }{ N -2}.
\end{align*}
Combining this with \eqref{eq5.18v65} and a standard rescaling argument gives \eqref{eq5.15v65}.
\medskip

\textbf{Case 3}.
The proof of this case is in parallel to that of Case 2. The differing geometry of the two cases enters only in the use of Proposition \ref{pr3.6v65} and the analogue of estimate \eqref{eq5.13v65}.
As these first two inputs have already been proven in all cases, we only need to show
\begin{align}\label{eq5.19v65}
\chi_n \tilde{\phi} \to \tilde{\phi}, \text{ or equivalently, }
\Theta  \left( \frac{ |x|}{ \dist  \left( 0, \Omega_n^c \right)}  \right) \tilde{\phi}(x) \to 0 \text{ in } \dot{H}^1( \mathbb{R}^N ).
\end{align}
To this end, let
\begin{align*}
B_n : =  \left\{ x \in \mathbb{R}^N : |x| \ge \frac14 \dist \left(0, \Omega_n^c \right)  \right\}.
\end{align*}
Then by H\"older's inequality, we have
\begin{align*}
\left\| \Theta  \left( \frac{|x|}{ \dist \left(0, \Omega_n^c \right)}  \right) \tilde{\phi} (x)  \right\|_{\dot{H}^1( \mathbb{R}^N )}
\lesssim  \left\|\nabla \tilde{\phi}  \right\|_{L^2( B_n)} +  \left\| \tilde{\phi}  \right\|_{L^\frac{2 N }{N -2} (B_n)} .
\end{align*}
which converges to 0 by the dominated convergence theorem (as $1_{B_n} \to 0$ almost everywhere).
\medskip

\textbf{Case 4}. Passing to a subsequence, we may assume $N_n d(x_n) \to d_\infty > 0$.
By weak sequential compactness of balls in $\dot{H}^1( \mathbb{R}^N )$, we may find a subsequence and a $\tilde{\phi}\in \dot{H}^1( \mathbb{R}^N )$ so that $g_n \rightharpoonup \tilde{\phi} $ in this space.
Using the following useful characterization of Sobolev spaces,
\begin{align*}
\dot{H}_D^1( \mathbb{R}_+^N )  =  \left\{g \in \dot{H}^1( \mathbb{R}^N ): \int_{\mathbb{R}^N } g(x) \psi(x) \,\mathrm{d}x = 0 \text{ for all } \psi \in C_c^\infty (-\mathbb{R}_+^N )  \right\},
\end{align*}
we find $\tilde{\phi} \in \dot{H}_D^1(  \mathbb{R}_+^N )$ by noticing that for any compact set $K$ in the halfspace, $K \subseteq \Omega_n^c$ for $n$ sufficiently large.
Here $\Omega_n: = N_n R_n^{-1}  \left( \Omega -  \left\{ x_n^*  \right\}  \right)\supset \supp(g_n)$.
Since $\tilde{\phi} \in \dot{H}_D^1 ( \mathbb{R}_+^N )$, we have $\phi_n \in \dot{H}_D^1( \Omega)$, as is easily seen from
\begin{align*}
x \in \mathbb{R}_+^N \Longleftrightarrow N_n^{-1} R_n x + x_n^* \in \mathbb{R}_{+,n}^N : =  \left\{ y:  \left( x_n - x_n^* \right)  \left(y - x_n^* \right) > 0 \right\} \subseteq \Omega.
\end{align*}
This inclusion further shows that
\begin{align}\label{eq5.20v65}
\left\| \tilde{\phi}   \right\|_{\dot{H}_D^1(  \mathbb{R}_+^N )} = \|\phi_n \|_{\dot{H}_D^1(\mathbb{R}_{+,n}^N ) } = \|\phi_n \|_{\dot{H}_D^1 ( \Omega)}.
\end{align}

To prove \eqref{eq5.1v65}, it is thus sufficient to find a lower bound on $ \left\|\tilde{\phi}  \right\|_{\dot{H}_D^1 ( \mathbb{R}_+^N )}$. For this, let $h: = P_1^{ \mathbb{R}_+^N } \delta_{d_\infty e_N }$.
Then it follows from the Bernstein inequality that
\begin{align}\label{eq5.21v65}
 \left\|  \left( - \Delta_{\mathbb{R}_+^N} \right)^{- \frac12} h  \right\|_{L^2( \Omega)} \lesssim 1,
\end{align}
which also implies $h \in \dot{H}_D^{-1} ( \mathbb{R}_+^N )$.

Set $\tilde{x}_n : = N_n R_n^{-1}  \left( x_n - x_n^* \right)$. Then, $\tilde{x}_n \to d_\infty e_N $ and we obtain from  Proposition \ref{pr3.6v65} that
\begin{align*}
 \left\langle \tilde{\phi}, h  \right\rangle
& = \lim\limits_{n \to \infty}  \left(  \left\langle g_n , P_1^{\Omega_n} \delta_{\tilde{x}_n}  \right\rangle +  \left\langle g_n,  \left(P_1^{ \mathbb{R}_+^N } - P_1^{\Omega_n }  \right) \delta_{d_\infty e_N } \right\rangle
+  \left\langle g_n, P_1^{\Omega_n }  \left( \delta_{d_\infty e_N  } - \delta_{\tilde{x}_n }  \right)  \right\rangle  \right)
\\
&= \lim\limits_{n \to \infty}  \left( N_n^{- \frac{N-2}2}  \left(P_{N_n}^\Omega f_n  \right)(x_n ) +  \left\langle g_n, P_1^{\Omega_n }  \left( \delta_{d_\infty e_N  } - \delta_{\tilde{x}_n }  \right)  \right\rangle  \right).
\end{align*}
Arguing as in the proof of \eqref{eq5.10v65} and applying \eqref{eq5.9v65} to $v(x) =  \left(P_1^{\Omega_n } g_n  \right) \left( x + \tilde{x}_n  \right)$ with
$ R = \frac12 N_n d(x_n)$,  we obtain
\begin{align*}
\left|  \left\langle g_n, P_1^{\Omega_n}  \left( \delta_{d_\infty e_N } - \delta_{\tilde{x}_n}  \right)  \right\rangle  \right|
\lesssim A  \left(d_\infty^{-1} + d_\infty  \right)  \left|d_\infty e_N  - \tilde{x}_n  \right| \xrightarrow{n  \to \infty} 0.
\end{align*}
Therefore, we conclude
\begin{align*}
 \left|  \left\langle \tilde{\phi}, h \right\rangle  \right| \gtrsim \epsilon^{\frac{N^2}4} A^{1- \frac{N^2}4},
\end{align*}
which together with \eqref{eq5.20v65} and \eqref{eq5.21v65} gives \eqref{eq5.1v65}.

For the proof of \eqref{eq5.2v65}, we simply note that
\begin{align*}
\| f_n \|_{\dot{H}_D^1 ( \Omega)}^2 - \|f_n - \phi_n \|_{\dot{H}_D^1 ( \Omega)}^2
& = 2  \left\langle f_n , \phi_n  \right\rangle_{\dot{H}_D^1 ( \Omega)} - \| \phi_n \|_{\dot{H}_D^1( \Omega)}^2
\\
& = 2  \left\langle g_n, \tilde{\phi}  \right\rangle_{\dot{H}_D^1 ( \Omega_n )} -  \left\|\tilde{\phi}  \right\|_{\dot{H}_D^1 (  \mathbb{R}_+^N )}^2
\xrightarrow{n \to \infty}  \left\|\tilde{\phi}  \right\|_{\dot{H}_D^1( \mathbb{R}_+^N )}^2.
\end{align*}
The proof of \eqref{eq5.3v65} is very similar to the cases treated previously: one uses the Rellich-Kondrashov compactness theorem
 to show $g_n \to \tilde{\phi}$ almost everywhere and then the refined Fatou lemma to see that
\begin{align*}
\text{LHS of } \ \eqref{eq5.3v65} =  \left\|\tilde{\phi}  \right\|_{L^\frac{2N }{ N -2}( \mathbb{R}_+^N )}^\frac{2N }{ N -2}.
\end{align*}
Pairing with $h$ as in Cases 1 and 2, the lower bound on this quantity immediately follows.

\end{proof}

Finally, we are ready to prove Theorem \ref{th5.6v65}, based on an inductive application of Proposition \ref{pr5.2v65}.

\begin{proof}[Proof of Theorem \ref{th5.6v65}]

We will proceed inductively and extract one bubble at each time.
In the initial step, we set $w_n^0 : = f_n$.
Suppose we have a decomposition up to level $J \ge 0$ satisfying \eqref{eq5.26v65new} -- \eqref{eq5.28v65new}.
Passing to a subsequence if necessary, we may assume
$A_J : = \lim\limits_{n \to \infty}   \left\|w_n^J  \right\|_{\dot{H}_D^1 ( \Omega)}$
and $\epsilon_J : = \lim\limits_{n \to \infty}  \left\|w_n^J  \right\|_{L_x^\frac{2 N }{ N -2}( \Omega)}$.

If $\epsilon_J = 0$, then we just stop and set $J^* = J$.
If not, then we apply Proposition \ref{pr5.2v65} to $w_n^J$.
Passing to a subsequence in $n$, we find
$ \left\{ \phi_n^{J+1}  \right\} \subseteq \dot{H}_D^1( \Omega)$, $\left\{ \lambda_n^{J+1}  \right\} \subseteq 2^{\mathbb{Z}}$, and $ \left\{ x_n^{J+1}  \right\} \subseteq \Omega$, which conform to one of the four cases listed in the theorem. Then we relabel the parameters given by Proposition \ref{pr5.2v65} as follows:
$\lambda_n^{J+1} : = N_n^{-1}$ and the profiles $\tilde{\phi}^{J+1}$ are defined as weak limits:
\begin{align*}
\left(g_n^{J+1}  \right)^{-1} w_n^J  \rightharpoonup \tilde{\phi}^{J+1}, \text{ as } n \to \infty,
\end{align*}
where $g_n^{J+1}$ is given as in the theorem.
In Cases 2, 3, 4, we define $\phi^{J+1} : = \tilde{\phi}^{J+1}$, while in Case 1,
\begin{align*}
\phi^{J+1} (x) : = g_\infty^{J+1} \tilde{\phi}^{J+1} (x) : =  \left( \lambda_\infty^{J+1}  \right)^{- \frac{N-2}2 } \tilde{\phi}^{J+1}  \left( \frac{ x - x_\infty^{J+1} }{ \lambda_\infty^{J+1} }  \right).
\end{align*}
Finally, $\phi_n^{J+1}$ is defined as in the theorem.

In Case 1, we can rewrite it as
$\phi_n^{J+1} = g_\infty^{J+1} \tilde{\phi}^{J+1}$,
where $\Omega_\infty^{J+1} : =  \left( \lambda_\infty^{J+1}  \right)^{-1}  \left( \Omega -  \left\{ x_\infty^{J+1}  \right\}  \right)$ and that in all four cases, we have
\begin{align}\label{eq5.30v65}
\lim\limits_{n \to \infty }  \left\|  \left(g_n^{J+1}  \right)^{- 1} \phi_n^{J+1} - \tilde{\phi}^{J+ 1}  \right\|_{\dot{H}^1( \mathbb{R}^N ) } = 0;
\end{align}
see also \eqref{eq5.13v65} and \eqref{eq5.19v65} for Cases 2 and 3.

Next, we define $w_n^{J+1} : = w_n^J - \phi_n^{J+1}$.
By \eqref{eq5.30v65} and the construction of $\tilde{\phi}^{J+1}$ in each case, we have
\begin{align*}
\left( g_n^{J+1}  \right)^{-1} w_n^{J+ 1} \rightharpoonup 0 \quad\text{ in } \dot{H}^1( \mathbb{R}^N ), \text{ as } n \to \infty,
\end{align*}
which proves \eqref{eq5.28v65new} at the level $J+1$. Moreover, we also infer from Proposition \ref{pr5.2v65} that
\begin{align*}
\lim\limits_{n \to \infty}  \left(  \left\|w_n^J  \right\|_{\dot{H}_D^1( \Omega)}^2 -  \left\|\phi_n^{J+1} \right\|_{\dot{H}_D^1 ( \Omega)}^2 -  \left\|w_n^{J+ 1} \right\|_{\dot{H}_D^1( \Omega)}^2  \right) = 0.
\end{align*}
This together with the inductive hypothesis gives \eqref{eq5.26v65new} at the level $J+1$.
By a similar argument, we can prove \eqref{eq5.27v65new} at level $J+1$.

From Proposition \ref{pr5.2v65}, passing to a further subsequence, we obtain
\begin{align}\label{eq5.31v65}
\begin{split}
A_{J+1}^2 = \lim\limits_{n \to \infty}  \left\|w_n^{J+1}  \right\|_{\dot{H}_D^1( \Omega)}^2 \le A_J^2
\left( 1 - c  \left( \frac{\epsilon_J}{A_J}  \right)^\frac{N^2 }2   \right) \le A_J^2, \\
\epsilon_{J+1}^{\frac{2 N }{ N -2}} = \lim\limits_{n \to \infty}  \left\|w_n^{J+1}  \right\|_{L_x^\frac{2 N }{ N -2} ( \Omega)}^\frac{2 N }{ N -2}
 \le \epsilon_J^{ \frac{2N}{N-2} }
\left( 1- c \left( \frac{\epsilon_J}{A_J}  \right)^\frac{N(N+2)}2  \right).
\end{split}
\end{align}

If $\epsilon_{J+ 1} = 0$, then we stop and set $J^* = J+ 1$; in this case, %we have
\begin{align}\label{eq5.25v}
\lim\limits_{J \to J^*} \limsup\limits_{n \to \infty} \|w_n^J \|_{L_x^\frac{2N}{N-2}( \Omega)} = 0
\end{align}
is automatic.
If $\epsilon_{J+ 1} > 0$, then we continue the induction process. If this process does not terminate in finitely many steps, then we set $J^* = \infty$; in this case, \eqref{eq5.31v65} implies $\epsilon_J \xrightarrow{J \to \infty} 0$ and so \eqref{eq5.25v} follows.

The desired equation
\begin{align*}
\lim\limits_{J \to \infty} \limsup\limits_{n \to \infty}  \left\| e^{tz \Delta_\Omega} w_n^J  \right\|_{L_t^\infty L_x^\frac{2 N }{ N -2} (\mathbb{R}_+ \times \Omega )} = 0.
\end{align*}
 is an immediate consequence of \eqref{eq5.25v} by noticing the well-known fact that for any $h  \in L^\frac{2 N }{ N -2}( \Omega  )$,
\begin{align*}
\left\|e^{tz \Delta_\Omega}  h   \right\|_{L_t^\infty L_x^\frac{2N }{ N -2} ( \mathbb{R}_+ \times  \Omega  )} \le \| h  \|_{L^\frac{2 N }{ N -2}_x ( \Omega )}.
\end{align*}
The estimate \eqref{eq5.25v65new} follows as
\begin{align*}
\|e^{tz \Delta_\Omega} w_n^J \|_{L_{t,x}^\frac{2(N+2)}{N-2} }
\le
\|e^{tz \Delta_\Omega} w_n^J \|_{L_t^\infty L_x^\frac{2N}{N-2}}^\frac4{N+2} \|e^{tz \Delta_\Omega} w_n^J \|_{L_t^2 L_x^\frac{2N}{N-4}}^\frac{N-2}{N+2}.
\end{align*}
Next we verify the asymptotic orthogonality condition \eqref{eq5.29v65new} by a contradiction argument. Assume \eqref{eq5.29v65new} fails to be true for some pair $(j, k)$. Without loss of generality, we may assume $j < k$ and \eqref{eq5.29v65new} holds for all pairs $(j,l)$ with $j < l < k$.
Passing to a subsequence, we may assume
\begin{align*}
\frac{\lambda_n^j}{ \lambda_n^k} \to \lambda_0 \in (0, \infty), \quad \frac{x_n^j - x_n^k}{ \sqrt{\lambda_n^j \lambda_n^k} } \to x_0, \text{ as } n \to \infty.
\end{align*}
From the inductive relation
\begin{align*}
w_n^{k-1} = w_n^j - \sum\limits_{l = j+1}^{k - 1} \phi_n^l
\end{align*}
and the definition of $\tilde{\phi}^k$, we obtain
\begin{align*}
\tilde{\phi}^k = \underset{ n \to \infty} {w-\lim }  \left(g_n^k \right)^{-1} w_n^{k - 1}\notag  =  \underset{n \to \infty } { w-\lim } \left( g_n^k \right)^{-1} w_n^j -  \sum\limits_{l = j+1}^{k-1}  \underset{n \to \infty } { w-\lim }
\left(g_n^k \right)^{-1} \phi_n^l=:A_1+A_2.
\end{align*}
Next, we prove that these weak limits are zero, which would be a contradiction to the non-triviality of $\tilde{\phi}^k$. Rewrite $A_1$ as follows:
\begin{align*}
\left( g_n^k \right)^{-1} w_n^j =  \left(g_n^k \right)^{-1} g_n^j  \left(g_n^j \right)^{-1} w_n^j.
\end{align*}
By \eqref{eq5.28v65new}, and the fact that adjoint of the unitary operators $\left(g_n^k  \right)^{-1} g_n^j$ converge strongly, we obtain that the weak limit in $A_1$ is $0$.

To complete the proof of \eqref{eq5.29v65new}, it remains to show $A_2 = 0$. For all $j < l < k$, we write
\begin{align*}
\left( g_n^k \right)^{-1} \phi_n^l =  \left(g_n^k \right)^{-1} g_n^j  \left(g_n^j \right)^{-1} \phi_n^l.
\end{align*}
Arguing as for $A_1$, it suffices to show
\begin{align*}
\left( g_n^j \right)^{-1} \phi_n^l \rightharpoonup 0 \text{ in } \dot{H}^1( \mathbb{R}^N ) \text{ as } n \to \infty .
\end{align*}
Using a density argument, this reduces to the statement that for all $\phi \in C_c^\infty  \left( \widetilde{\lim} \,  \Omega_n^l \right)$
\begin{align}\label{eq5.35v65}
I_n : =  \left( g_n^j \right)^{-1} g_n^l \phi \rightharpoonup 0 \text{ in } \dot{H}^1( \mathbb{R}^N ), \text{ as } n \to \infty.
\end{align}
Depending on which cases $j$ and $l$ fall into, we can rewrite $I_n$ as follows:
\begin{itemize}
\item[(a).]
If both $j$ and $l$ conform to Case 1, 2, or 3, then
\begin{align*}
I_n =  \left( \frac{\lambda_n^j}{ \lambda_n^l}  \right)^\frac{N-2}2 \phi \left( \frac{ \lambda_n^j x + x_n^j - x_n^l }{ \lambda_n^l}  \right).
\end{align*}

\item[(b).]
If $j$ conforms to Case 1, 2, or 3 and $l$ to Case 4, then
\begin{align*}
I_n =  \left( \frac{\lambda_n^j}{ \lambda_n^l} \right)^\frac{N-2}2 \phi  \left( \frac{  \left(R_n^l \right)^{-1}  \left( \lambda_n^j x + x_n^j -  \left(x_n^l \right)^*  \right)}{ \lambda_n^l} \right).
\end{align*}

\item[(c).]
If $j$ conforms to Case 4 and $l$ to Case 1, 2, or 3, then
\begin{align*}
I_n =  \left( \frac{\lambda_n^j}{ \lambda_n^l} \right)^\frac{N-2}2 \phi  \left( \frac{ R_n^j \lambda_n^j x +  \left(x_n^j \right)^* - x_n^l}{\lambda_n^l}  \right).
\end{align*}

\item[(d).]
If both $j$ and $l$ conform to Case 4, then
\begin{align*}
I_n =  \left( \frac{\lambda_n^j}{ \lambda_n^l}  \right)^\frac{N-2}2 \phi  \left( \frac{  \left(R_n^l \right)^{-1}  \left(R_n^j \lambda_n^j x +  \left(x_n^j \right)^* - \left(x_n^l \right)^*  \right)}{ \lambda_n^l}  \right).
\end{align*}

\end{itemize}

We first prove \eqref{eq5.35v65} when the scaling parameters are not comparable, that is,
\begin{align*}
\lim\limits_{n \to \infty}  \left( \frac{\lambda_n^j}{\lambda_n^l} + \frac{ \lambda_n^l}{ \lambda_n^j} \right)  = \infty.
\end{align*}
In this case, we shall treat all the four cases simultaneously. By the Cauchy-Schwarz inequality, for all $\psi \in C_c^\infty( \mathbb{R}^N )$, we have
\begin{align*}
\left| \langle I_n, \psi \rangle_{\dot{H}^1( \mathbb{R}^N )}  \right|
& \lesssim \min  \left( \|\Delta I_n \|_{L^2( \mathbb{R}^N )} \|\psi \|_{L^2( \mathbb{R}^N )} , \|I_n \|_{L^2( \mathbb{R}^N ) } \|\Delta \psi \|_{L^2( \mathbb{R}^N )}  \right)
\\
& \lesssim \min  \left( \frac{ \lambda_n^j}{ \lambda_n^l} \|\Delta \phi \|_{L^2( \mathbb{R}^N )} \|\psi \|_{L^2( \mathbb{R}^N )}, \frac{\lambda_n^l}{ \lambda_n^j} \|\phi \|_{L^2( \mathbb{R}^N )} \|\Delta \psi \|_{L^2( \mathbb{R}^N )}  \right),
\end{align*}
which converges to zero as $n \to \infty$.
Thus, in this case $A_2 = 0$ and we get the desired contradiction.

It remains to deal with the situation when
\begin{align*}
\frac{\lambda_n^j}{ \lambda_n^l} \to \lambda_0\in (0, \infty), \ \frac{ \left|x_n^j - x_n^l  \right|^2}{ \lambda_n^j \lambda_n^l} \to \infty, \text{ as } n \to \infty .
\end{align*}
In this case, we need the following simple observation: assume either $\Omega_n \equiv \Omega$ or $\{\Omega_n \}$ conforms to one of the three scenarios considered in Proposition \ref{pr5.2v65}.
Let $f \in C_c^\infty  \left( \widetilde{\lim} \, \Omega_n \right)$ and let $\{x_n \}_{n \ge 1} \subseteq \mathbb{R}^N $.
Then
\begin{align*}
	f (x+ x_n ) \rightharpoonup 0 \quad\text{ in } \dot{H}^1( \mathbb{R}^N  ) \text{ as } n \to \infty,
\end{align*}
whenever $|x_n | \to \infty$.

Case (a) follows directly from the above observation as it implies
\begin{align*}
\lambda_0^\frac{N-2}2 \phi \left( \lambda_0 x + y_n  \right) \rightharpoonup 0 \text{ in } \dot{H}^1( \mathbb{R}^N ), \text{ as } n \to \infty,
\end{align*}
where $y_n : = \frac{ x_n^j - x_n^l}{ \lambda_n^l} = \frac{ x_n^j - x_n^l}{ \sqrt{\lambda_n^l \lambda_n^j} } \sqrt{ \frac{ \lambda_n^j}{ \lambda_n^l} } \to \infty$, as $n \to \infty$.

For Case (b), as $SO(N )$ is a compact group, we can proceed similarly provided that
\begin{align*}
\frac{  \left|x_n^j -  \left(x_n^l \right)^*  \right|}{ \lambda_n^l} \to \infty \text{ as } n \to \infty.
\end{align*}
But this is immediate consequence of the triangle inequality:
for $n$ sufficiently large,
\begin{align*}
\frac{  \left| x_n^j -  \left(x_n^l \right)^*  \right|}{ \lambda_n^l} \ge \frac{  \left|x_n^j - x_n^l \right|}{ \lambda_n^l} - \frac{  \left|x_n^l -  \left(x_n^l \right)^*  \right|}{ \lambda_n^l} \ge \frac{  \left|x_n^j - x_n^l  \right|}{ \lambda_n^l} - 2 d_\infty^l \xrightarrow{n \to \infty} \infty.
\end{align*}
The proof of Case (c) is entirely similar. For Case (d), note that for $n$ sufficiently large, we have
\begin{align*}
\frac{  \left| \left(x_n^j \right)^* -  \left(x_n^l \right)^*  \right|}{ \lambda_n^l} & \ge \frac{  \left|x_n^j - x_n^l \right|}{ \lambda_n^l}
- \frac{  \left|x_n^j -  \left(x_n^j \right)^*  \right|}{ \lambda_n^l} - \frac{  \left|x_n^l -  \left(x_n^l \right)^*  \right|}{ \lambda_n^l} \\
& \ge \frac{  \left|x_n^j - x_n^l \right|}{ \sqrt{\lambda_n^j \lambda_n^l}} \sqrt{ \frac{\lambda_n^j} {\lambda_n^l} } - \frac{d \left(x_n^j \right) }{ \lambda_n^j} \frac{ \lambda_n^j}{ \lambda_n^l} - \frac{ d \left(x_n^l \right)}{ \lambda_n^l}\\
&\ge \frac12 \sqrt{ \lambda_0 } \frac{  \left|x_n^j - x_n^l \right|}{ \sqrt{ \lambda_n^j \lambda_n^l}} - 2 \lambda_0 d_\infty^j - 2  d_\infty^l \xrightarrow{n\to\infty} \infty.
\end{align*}
The desired weak convergence follows again from the observation.
This completes the proof.
\end{proof}

\subsection{Approximation of profiles}\label{subse3.2}

In this subsection, we show the uniform boundedness of the spacetime norm of nonlinear profiles.
Since embedding these nonlinear profiles corresponds to different limiting geometry (back to $\Omega$), we need to consider several different cases.

We first treat the case when the rescaled domains $\Omega_n$ expand to fill $\mathbb{R}^N$, where $\Omega_n = \lambda_n^{-1} ( \Omega- \{x_n \} )$. This case corresponds to Case 2 of Theorem \ref{th5.6v65}.

\begin{proposition}\label{th6.1}
Suppose that $\{\lambda_n \} \subseteq 2^{\mathbb{Z}}$ satisfies $\lambda_n \to \infty$ as $n \to \infty$, and that $\{ x_n\} \subseteq \Omega$ satisfies $- \frac{x_n}{ \lambda_n} \to x_\infty \in \mathbb{R}^N  $ as $n \to \infty$. Let $\phi \in \dot{H}^1( \mathbb{R}^N )$ and define
\begin{align*}
\phi_n(x) =  \frac1{ \lambda_n^{\frac{N -2}2}} ( \chi_n \phi) \left( \frac{x-x_n}{ \lambda_n} \right),
\end{align*}
where $\chi_n(x) = \Theta  \left( \frac{d  \left( \lambda_n x + x_n \right) }{ \diam  \left( \Omega^c \right)} \right).$
Then for $n$ sufficiently large there exists a global solution $v_n$ to $\CGL_\Omega$ with initial data $v_n(0) = \phi_n$ satisfying
$$\|v_n \|_{L_{t,x}^{ \frac{2( N +2)}{ N -2} }( \mathbb{R}_+ \times \Omega)} \lesssim 1,$$
with the implicit constant depending only on $\|\phi\|_{\dot{H}^1}$.

\end{proposition}

\begin{proof}

To prove this result, we shall use the stability Theorem \ref{th2.9}. For this, we need to obtain approximate solutions to $\CGL_\Omega$ that satisfy the necessary estimates required in Theorem \ref{th2.9}; see \eqref{eq6.11}, \eqref{eq6.12}, and \eqref{eq6.14} below.

The first step towards such approximate solutions to $\CGL_\Omega$ is to  construct global solutions to $\CGL_{\mathbb{R}^N }$. More precisely, fix $0<\theta\ll 1$ and let $w_n$ and $w_\infty$ be solutions to $\CGL_{\mathbb{R}^N  }$ with initial data $w_n(0) =
P _{\le \lambda_n^\theta}^{\mathbb{R}^N } \phi $ and $w_\infty (0) = \phi$.
By Theorem \ref{th1.1v37}, we know that $w_n$ and $w_\infty$ are global solutions satisfying
\begin{align}\label{eq6.3}
\|w_n \|_{\dot{S}^1( \mathbb{R}_+ \times \mathbb{R}^N  )} + \|w_\infty \|_{\dot{S}^1( \mathbb{R}_+ \times \mathbb{R}^N )} \lesssim 1,
\end{align}
with the implicit constant depending only on $\|\phi \|_{\dot{H}^1}$. Moreover, by Theorem \ref{th2.9}, we have
\begin{align*}
\lim\limits_{n \to \infty}  \left\|w_n - w_\infty  \right\|_{\dot{S}^1( \mathbb{R}_+ \times \mathbb{R}^N  )}  = 0.
\end{align*}
It follows from Bernstein's inequality that
\begin{align*}
 \left\| P _{\le \lambda_n^\theta}^{\mathbb{R}^N } \phi   \right\|_{\dot{H}^s} \lesssim \lambda_n^{\theta( s- 1)} \quad \text{ for any } s \ge 1,
\end{align*}
and so Lemma \ref{le2.10} gives
\begin{align*}
 \left\| |\nabla |^s w_n  \right\|_{\dot{S}^1( \mathbb{R}_+ \times \mathbb{R}^N  ) } \lesssim \lambda_n^{\theta s } \quad \text{ for any } 0 \le s \le 2,
\end{align*}
with the implicit constant depending only on $\|\phi \|_{\dot{H}^1}$. Combining this with the Gagliardo-Nirenberg inequality, we obtain for all $0 \le s \le 2 $,
\begin{align}\label{eq6.6}
 \left\| |\nabla |^s w_n  \right\|_{L_{t,x}^\infty } \lesssim \lambda_n^{ \theta  \left( s + \frac12  \right)}.
\end{align}
Finally, using the $\CGL_{\mathbb{R}^N }$ equation we get
\begin{align}\label{eq6.7}
 \left\| \partial_t w_n \right\|_{L_{t,x}^\infty} \lesssim \|\Delta w_n \|_{L_{t,x}^\infty } + \|w_n \|_{L_{t,x}^\infty}^\frac{N +2}{ N -2}  \lesssim \lambda_n^{\frac{5}{2}  \theta}.
\end{align}
To construct the approximate solution to $\CGL_\Omega$, we add the following term $z_n$ to $\chi_n w_n(t)$:
\begin{align*}
z_n(t) : =  \bar{z}  \int_0^t e^{  (t - s)  {z} \Delta_{\Omega_n } }
\left( \Delta_{\Omega_n } \chi_n \right) w_n  \left( s, - \lambda_n^{- 1} x_n \right) \,\mathrm{d}s,
\end{align*}
which comes from the high frequency reflections off the obstacle.

We next show that for any $T > 0$ and for all  $0 \le s < \frac{N }2$, there holds
\begin{equation}\label{eq6.8}
 \limsup\limits_{n \to \infty} \|z_n \|_{\dot{X}^1 ([0, T] \times \Omega_n )} = 0
\end{equation}
and
\begin{equation}\label{eq6.9}
	\left\|  \left( - \Delta_{\Omega_n}  \right)^\frac{s}2 z_n \right\|_{L_t^\infty L_x^2 ([0, T] \times \Omega_n )} \lesssim \lambda_n^{s - \frac{N}2 + \frac52\theta}
	\left( T+ \lambda_n^{- 2 \theta} \right).
\end{equation}
Indeed, integration by parts, we have
\begin{align*}
z_n(t) & = - \int_0^t  \left( e^{t  {z} \Delta_{\Omega_n}} \partial_s e^{- s  {z} \Delta_{\Omega_n}} \chi_n  \right) w_n \left( s, - \lambda_n^{- 1} x_n \right) \,\mathrm{d}s
\\
& = - \chi_n w_n  \left( t, - \lambda_n^{- 1} x_n  \right) + e^{t  {z} \Delta_{\Omega_n }}  \left( \chi_n w_n(0, - \lambda_n^{- 1} x_n) \right)
+ \int_0^t  \left( e^{(t- s)  {z} \Delta_{\Omega_n} } \chi_n  \right) \partial_s w_n  \left(s, - \lambda_n^{- 1} x_n  \right) \,\mathrm{d}s.
\end{align*}
We first estimate the $L_t^\frac{ N +2}{ N -2}  \dot{W}_x^{1, \frac{2 N ( N +2)}{ N^2 - 2 N + 8 }}$ norm of $z_n$. Using the Strichartz inequality, Theorem \ref{th2.3}, \eqref{eq6.6} and \eqref{eq6.7}, we get
\begin{align*}
& \quad \| z_n \|_{L_t^\frac{ N +2}{ N -2}  \dot{W}_x^{1, \frac{2 N ( N +2) }{ N^2 - 2 N  + 8 }} }
\\
& \lesssim
\left\|\nabla \chi_n (x) w_n \left(t, - \lambda_n^{- 1} x_n \right)  \right\|_{L_t^\frac{ N +2}{ N -2} L_x^\frac{2 N ( N +2)}{ N^2 - 2 N + 8 }}
+  \left\|\nabla \chi_n (x) w_n \left(0, - \lambda_n^{- 1} x_n \right)  \right\|_{L_x^2}
+  \left\|\nabla \chi_n(x) \partial_t w_n  \left(t, - \lambda_n^{- 1} x_n \right)  \right\|_{L_t^1 L_x^2}
\\
& \lesssim T^\frac{ N -2}{ N +2}  \left\|\nabla \chi_n  \right\|_{L_x^\frac{2 N ( N +2)}{ N^2 - 2 N + 8 }}
\|w_n \|_{L_{t,x}^\infty}
+  \left\|\nabla \chi_n  \right\|_{L_x^2} \|w_n \|_{L_{t,x}^\infty}
+ T  \left\|\nabla \chi_n  \right\|_{L_x^2} \| \partial_t w_n \|_{L_{t,x}^\infty}
\\
& \lesssim T^\frac{ N -2}{ N +2}  \lambda_n^{- \frac{( N -2)^2}{2( N +2) } + \frac\theta2}
+ \lambda_n^{- \frac{ N -2}2  + \frac\theta2} + T \lambda_n^{- \frac{ N -2}2  + \frac52 \theta} \to 0, \text{ as } n \to \infty.
\end{align*}
Similarly, using the Sobolev embedding, we obtain
\begin{align}\label{eq6.10}
& \quad \| z_n \|_{L_{t,x}^{ \frac{2( N +2)}{ N  - 2} }}
\\
& \lesssim  \left\|  \left( - \Delta_{\Omega_n} \right)^\frac12 z_n  \right\|_{L_t^{ \frac{2( N +2)}{N -2} } L_x^\frac{2 N ( N +2)}{ N^2 + 4 }} \notag
\\
& \lesssim  \left\| \nabla \chi_n (x) w_n \left(t, - \lambda_n^{- 1} x_n  \right)   \right\|_{L_t^{ \frac{2( N +2)}{ N -2} } L_x^\frac{2  N ( N +2) }{ N^2 + 4 }}
+  \left\| \nabla \chi_n (x) w_n \left(0, - \lambda_n^{- 1} x_n \right)  \right\|_{L_x^2}
+  \left\|\nabla \chi_n (x) \partial_t w_n  \left(t, - \lambda_n^{- 1} x  \right)  \right\|_{L_t^1 L_x^2} \notag
\\
& \lesssim T^\frac{ N -2}{2( N +2)}  \|\nabla \chi_n \|_{L_x^\frac{ 2 N ( N +2) }{ N^2 +4 }} \| w_n \|_{L_{t,x}^\infty}
+ \|\nabla \chi_n \|_{L_x^2} \|w_n \|_{L_{t,x}^\infty}
+ T \|\nabla \chi_n \|_{L_x^2} \|\partial_t w_n \|_{L_{t,x}^\infty } \notag
\\
& \lesssim T^\frac{ N -2}{2( N +2) } \lambda_n^{ - \frac{ N ( N -2)}{2( N +2)} + \frac\theta2} + \lambda_n^{- \frac{ N -2}2 + \frac\theta2 }
+ T \lambda_n^{ - \frac{ N -2}2  + \frac52 \theta} \to 0, \text{ as } n \to \infty. \notag
\end{align}
This proves \eqref{eq6.8}.

To establish \eqref{eq6.9}, we argue as before and estimate
\begin{align*}
& \quad \left\|  \left( - \Delta_{\Omega_n} \right)^\frac{s}2 z_n  \right\|_{L_t^\infty L_x^2} \\
&\lesssim    \left\|  \left( - \Delta \right)^\frac{s}2 \chi_n w_n \left(t, - \lambda_n^{- 1} x_n  \right)  \right\|_{L_t^\infty L_x^2}
+  \left\| \left( - \Delta \right)^\frac{s}2 \chi_n w_n \left(0, - \lambda_n^{- 1} x_n  \right)  \right\|_{L_x^2}
+  \left\|  \left( - \Delta  \right)^\frac{s}2 \chi_n \partial_t w_n \left( t, - \lambda_n^{- 1} x_n  \right)  \right\|_{L_t^1 L_x^2} \\
&\lesssim
   \left\|   \left( - \Delta  \right)^\frac{s}2 \chi_n  \right\|_{L_x^2} \|w_n \|_{L_{t,x}^\infty}
+ T  \left\|  \left( - \Delta \right)^\frac{s}2  \chi_n  \right\|_{L_x^2} \|\partial_t w_n \|_{L_{t,x}^\infty} \\
&\lesssim \lambda_n^{s - \frac{ N }2 + \frac\theta2} + T \lambda_n^{ s- \frac{ N }2 + \frac52 \theta}
\lesssim \lambda_n^{s- \frac{ N }2 + \frac52 \theta}  \left( T + \lambda_n^{- 2 \theta}  \right).
\end{align*}

We are now in a position to introduce the approximate solution
\begin{align*}
\tilde{v}_n (t,x) =
\begin{cases}
 \frac1{ \lambda_n^{ \frac{ N -2}2} }   \left( \chi_n w_n + z_n  \right)  \left( \frac{t}{  \lambda_n^{ 2} },  \frac{ x- x_n}{ \lambda_n }  \right),   \ 0 \le t \le \lambda_n^2 T, \\
 \\
 e^{  \left( t - \lambda_n^2 T \right) {z} \Delta_\Omega } \tilde{v}_n  \left( \lambda_n^2 T , x  \right),          \qquad \qquad  t > \lambda_n^2 T,
\end{cases}
\end{align*}
where $T > 0$ is a parameter to be chosen later.
Note that $\tilde{v}_n$ has finite spacetime bound. Indeed, using a change of variables, the Strichartz inequality, \eqref{eq6.3}, \eqref{eq6.9} and \eqref{eq6.10}, we get
\begin{align}\label{eq6.11}
 \left\|\tilde{v}_n  \right\|_{L_{t,x}^{ \frac{2( N +2)}{ N - 2} } ( \mathbb{R}_+ \times \Omega)}
& \lesssim  \| \chi_n w_n + z_n \|_{L_{t,x}^{ \frac{2( N +2)}{ N -2} }([0, T] \times \Omega_n)} +  \left\| ( \chi_n w_n + z_n )(T)  \right\|_{\dot{H}_D^1 ( \Omega_n )}
\\
& \lesssim \|w_n \|_{L_{t,x}^{ \frac{2( N +2)}{ N -2} } ( \mathbb{R}_+ \times \mathbb{R}^N )} + \| z_n  \|_{L_{t,x}^{ \frac{2( N +2)}{ N -2} } ( [0, T] \times \Omega_n ) }
+ \| \chi_n \|_{L_x^\infty } \|\nabla w_n \|_{L_t^\infty L_x^2 ( \mathbb{R}_+ \times \mathbb{R}^N )} \notag
\\
& \quad  + \|\nabla \chi_n \|_{L_x^N   } \|w_n \|_{L_t^\infty L_x^\frac{2 N }{ N -2}  ( \mathbb{R}_+ \times \mathbb{R}^N  )}
+  \left\|  \left( - \Delta_{\Omega_n } \right)^\frac12 z_n  \right\|_{L_t^\infty L_x^2 ( [0, T] \times \Omega_n )} \notag
\\
& \lesssim 1 + T^\frac{ N -2}{2( N +2)} \lambda_n^{ - \frac{ N ( N -2)}{2( N +2)}  + \frac\theta2} + \lambda_n^{- \frac{ N -2}2  + \frac\theta2}
+ T \lambda_n^{- \frac{ N -2}2 + \frac52 \theta}. \notag
\end{align}
Next, observe that
\begin{align}\label{eq6.12}
\lim\limits_{T \to \infty } \limsup\limits_{n \to \infty }  \left\|  \left( - \Delta_\Omega \right)^\frac12 e^{t  {z} \Delta_\Omega}  \left( \tilde{v}_n (0) - \phi_n  \right)  \right\|_{L_t^{ \frac{2( N +2)}{ N -2}} L_x^\frac{2 N ( N +2)}{ N^2 + 4 } } = 0.
\end{align}
Indeed, using the Strichartz inequality, a change of variables, and H\"older's inequality, we estimate
\begin{align*}
 \left\| \left( - \Delta_\Omega \right)^\frac12 e^{t  {z} \Delta_\Omega }  \left( \tilde{v}_n (0) - \phi_n  \right)  \right\|_{L_t^{ \frac{2( N +2)}{ N -2} } L_x^\frac{ 2 N ( N +2)}{ N^2 + 4 }}
& \lesssim  \left\|  \left( - \Delta_\Omega \right)^\frac12  \left( \tilde{v}_n (0) - \phi_n  \right)  \right\|_{L_x^2}
\\
& \lesssim  \left\| \nabla  \left( \chi_n  P _{\le \lambda_n^\theta}^{\mathbb{R}^N } \phi  - \chi_n  \phi  \right)  \right\|_{L_x^2}
\\
& \lesssim \| \nabla \chi_n \|_{L_x^N  }  \left\| P _{\le \lambda_n^\theta}^{\mathbb{R}^N } \phi  - \phi  \right\|_{L_x^\frac{2N }{ N -2} }
+ \| \chi_n \|_{L_x^\infty }  \left\| \nabla \left( P _{\le \lambda_n^\theta}^{\mathbb{R}^N }  \phi  - \phi  \right) \right\|_{L_x^2},
\end{align*}
which converges to zero as $n\to \infty$.

Now, we claim that $\tilde{v}_n$ is an approximate solution to $\CGL_\Omega$ in the sense that
\begin{align}\label{eq6.14}
  \partial_t \tilde{v}_n - z  \Delta_\Omega \tilde{v}_n  =
  -  z \left|\tilde{v}_n  \right|^\frac4{ N -2}  \tilde{v}_n + e_n , \notag
\intertext{ with }
\lim\limits_{T \to \infty } \limsup\limits_{n \to \infty } \| e_n \|_{ \dot{N}^1( \mathbb{R}_+ \times \Omega   ) } = 0.
\end{align}
We start by verifying \eqref{eq6.14} on the large time interval $ \left( \lambda_n^2 T, \infty \right)$.
By the definition of $\tilde{v}_n$, in this regime we have $e_n =  z  \left|\tilde{v}_n \right|^\frac4{ N -2} \tilde{v}_n$.
Using Theorem \ref{th2.3}, the Strichartz inequality, and \eqref{eq6.9}, we estimate
\begin{align*}
\|  e_n \|_{\dot{N}^1   \left( \left( \lambda_n^2 T, \infty  \right)  \times \Omega \right) }
& \lesssim  \left\|  \left( - \Delta_\Omega \right)^\frac12  \left(  \left|\tilde{v}_n  \right|^\frac4{ N -2}  \tilde{v}_n  \right)  \right\|_{L_t^\frac{ N +2}{ N }  L_x^\frac{2 N ( N +2)}{ N^2 + 2 N   + 8 } \left(  \left( \lambda_n^2 T, \infty  \right) \times \Omega \right)}
\\
& \lesssim   \left\|  \left( - \Delta_\Omega \right)^\frac12 \tilde{v}_n  \right\|_{L_t^\frac{ N +2}{ N-2}  L_x^\frac{ 2 N ( N +2) }{ N^2 - 2 N  + 8 }
\left(
 \left( \lambda_n^2 T, \infty  \right) \times \Omega \right)}  \left\|\tilde{v}_n  \right\|_{L_{t,x}^{\frac{2( N +2)}{ N -2} }  \left(   \left( \lambda_n^2 T, \infty  \right) \times \Omega  \right) }^\frac4{ N -2}
\\
& \lesssim  \left\| \left( - \Delta_\Omega \right)^\frac12 \left( \chi_n w_n (T) + z_n (T)  \right)  \right\|_{L_x^2 }  \left\|\tilde{v}_n  \right\|_{L_{t,x}^{ \frac{2( N +2)}{ N -2} } \left(   \left( \lambda_n^2 T, \infty  \right) \times \Omega  \right)}^\frac4{ N -2}
\\
& \lesssim  \left( 1+ \lambda_n^{- \frac{ N -2}2 + \frac52 \theta }  \left(T + \lambda_n^{- 2 \theta }  \right)  \right)
\left\|\tilde{v}_n \right\|_{L_{t,x}^{ \frac{2( N +2)}{ N -2} }
\left(   \left( \lambda_n^2 T, \infty  \right) \times \Omega \right)}^\frac4{ N -2} .
\end{align*}
Thus, to establish \eqref{eq6.14} on $ \left( \lambda_n^2 T, \infty \right)$, it suffices to notice
\begin{equation}\label{eq6.15}
\begin{aligned}
&\quad\left\| e^{  \left( t- \lambda_n^2 T  \right) {z} \Delta_\Omega } \tilde{v}_n  \left( \lambda_n^2 T \right) \right\|_{L_{t,x}^{ \frac{2( N +2)}{ N -2}} (   \left( \lambda_n^2 T, \infty  \right) \times \Omega ) }\\
&\lesssim  \left\|\tilde{v}_n \left(\lambda_n^2 T \right) \right\|_{\dot{H}^1}\lesssim \|w_n(T)\|_{\dot{H}^1}
+ \|z_n(T)\|_{\dot{H}^1}\xrightarrow{n,T\to \infty}0.
\end{aligned}
\end{equation}
Thus we are left to prove \eqref{eq6.14} on the time interval $ \left[0, \lambda_n^2 T \right] $. For these values of time, we compute
\begin{align*}
e_n (t,x )
&  =  \left(  \left(  \partial_t  - z  \Delta_\Omega  \right) \tilde{v}_n
+  z \left|\tilde{v}_n  \right|^\frac4{ N -2}  \tilde{v}_n  \right)(t,x)
\\
& = - \lambda_n^{- \frac{ N +2}2}  z \left( \Delta \chi_n  \right) \left( \lambda_n^{- 1} ( x  - x_n )  \right) w_n  \left( \lambda_n^{- 2 } t, - \lambda_n^{- 1} x_n \right)
+ \lambda_n^{- \frac{ N +2}2 }  z \left( \Delta \chi_n  w_n  \right)  \left( \lambda_n^{- 2 } t , \lambda_n^{- 1} ( x - x_n ) \right)
\\
& \quad + 2 \lambda_n^{- \frac{ N +2}2 } z   \left( \nabla \chi_n \cdot \nabla w_n \right)  \left( \lambda_n^{- 2} t , \lambda_n^{- 1} ( x - x_n )  \right)\\
& \quad  + \lambda_n^{- \frac{ N +2}2 }  z \left( - \chi_n
 \left| w_n \right|^\frac4{ N -2} w_n
+   \left| \chi_n w_n + z_n  \right|^\frac4{ N -2}   \left( \chi_n w_n + z_n  \right)  \right) \left(\lambda_n^{- 2} t, \lambda_n^{- 1} ( x - x_n ) \right).
\end{align*}
Thus, using a change of variables and Theorem \ref{th2.3}, we estimate
\begin{align}
\|  e_n \|_{\dot{N}^1 \left( \left[0, \lambda_n^2 T \right]  \times \Omega \right) }
& \lesssim \left\|  \left( - \Delta_\Omega \right)^\frac12 e_n  \right\|_{L_{t,x}^\frac{2( N +2)}{ N +4}  \left(   \left[0, \lambda_n^2 T \right] \times \Omega  \right)}  \notag   \\
& \lesssim \left\| \nabla \left( \Delta \chi_n \left( w_n (t,x) - w_n  \left(t, - \lambda_n^{- 1} x_n \right) \right) \right)  \right\|_{L_{t,x}^\frac{2( N +2)}{ N +4}
([0, T] \times \Omega_n ) } \notag \\
& +  \left\| \nabla \left( \nabla \chi_n \cdot \nabla w_n  \right)  \right\|_{L_{t,x}^\frac{2( N +2)}{ N +4}  ([0, T] \times \Omega_n ) }  \notag \\
&   +  \left\| \nabla  \left( \chi_n | w_n |^\frac4{ N -2} w_n
 -  \left| \chi_n w_n + z_n \right|^\frac4{ N -2} \left( \chi_n w_n + z_n \right)  \right) \right\|_{L_{t,x}^\frac{2( N +2)}{ N +4}   ([0, T] \times \Omega_n ) } \notag \\
& := A_{1, n}  + A_{2, n}  + A_{3,n}. \notag
\end{align}
Using H\"older's inequality, the fundamental theorem of calculus, and \eqref{eq6.6}, we estimate
\begin{align*}
A_{1,n } & \lesssim T^\frac{ N +4}{2( N +2) } \| \Delta \chi_n \|_{L_x^\frac{2( N +2)}{ N +4}  } \| \nabla w_n \|_{L_{t,x}^\infty}
+ T^\frac{ N +4}{2( N +2) } \|\nabla \Delta  \chi_n \|_{L_x^\frac{2( N +2)}{ N +4}  }
\left\|w_n (t,x) - w_n  \left(t, - \lambda_n^{- 1} x_n \right)  \right\|_{L_{t,x}^\infty ( \mathbb{R}_+ \times \supp ( \Delta \chi_n ) )}
\\
& \lesssim  T^\frac{ N +4}{2( N +2) } \lambda_n^{- \frac{ N^2 - 8 }{2( N +2) } + \frac32 \theta}
+ T^\frac{ N +4}{2( N +2) } \lambda_n^{- \frac{ N^2  - 8 }{2( N +2) }} \lambda_n^{- 1} \|\nabla w_n \|_{L_{t,x}^\infty}
\lesssim T^\frac{ N +4}{2( N +2) } \lambda_n^{- \frac{ N^2 - 8 }{2( N +2) }  + \frac32 \theta } \xrightarrow{n\to \infty} 0.
\end{align*}
Notice that the cancellations induced by the introduction of $z_n$ were essential in order to control this term. Next,
\begin{align*}
A_{2, n}  &
\le T^\frac{ N +4}{2( N +2) }  \left( \|\Delta \chi_n \|_{L_x^\frac{2( N +2)}{ N +4}  } \| \nabla w_n \|_{L_{t,x}^\infty} + \|\nabla \chi_n \|_{L_x^\frac{2( N +2)}{ N +4} } \| \Delta w_n \|_{L_{t,x}^\infty} \right)\\
& \le T^\frac{ N +4}{2( N +2) } \left( \lambda_n^{- \frac{ N^2 - 8}{2( N +2) }+ \frac32 \theta } + \lambda_n^{- \frac{ N^2 + 2 N  - 4 }{ 2( N +2) }+ \frac52 \theta} \right) \xrightarrow{n\to \infty} 0.
\end{align*}
Finally, we turn our attention to $A_{3,n } $. A simple algebraic computation yields
\begin{align*}
A_{3, n}  \lesssim &  T^\frac{ N +4}{2( N +2) } %\left(
 \left\|\nabla \left( \left( \chi_n - \chi_n^\frac{ N +2}{ N -2}  \right) |w_n|^\frac{4}{ N -2}w_n   \right) \right\|_{L_t^\infty L_x^\frac{2( N +2)}{ N +4}  }
\\
& +  T^\frac{ N +4}{2( N +2) } \left(   \left\|z_n^\frac4{ N -2}  \nabla z_n  \right\|_{L_t^\infty L_x^\frac{2( N +2)}{ N +4} } +
 \left\| |z_n|^\frac4{N-2} \nabla ( \chi_n w_n )  \right\|_{L_t^\infty L_x^\frac{2(N+2)}{N+4} } +  \left\|  |\chi_n w_n|^\frac4{N-2} \nabla z_n  \right\|_{L_t^\infty L_x^\frac{2(N+2)}{N+4}} \right)
 \\
& +  T^\frac{ N +4}{2( N +2) }  \left\| \nabla ( \chi_n w_n) |\chi_n w_n|^\frac{6- N }{ N -2} z_n  \right\|_{L_t^\infty L_x^\frac{2(N+2)}{N+4} },
\end{align*}
where all spacetime norms are over $[0, T ] \times \Omega_n$. Using H\"older's inequality and \eqref{eq6.6}, we estimate
\begin{align*}
\left\| \nabla  \left(  \left( \chi_n - \chi_n^\frac{ N +2}{ N -2}   \right) w_n^\frac{ N +2}{ N -2}  \right) \right\|_{L_t^\infty L_x^\frac{2( N +2)}{ N +4}  }
& \lesssim  \| \nabla \chi_n \|_{L_x^\frac{2( N +2)}{ N +4}  } \|w_n \|_{L_{t,x}^\infty}^\frac{ N +2}{ N -2}
+  \left\|\chi_n - \chi_n^\frac{ N +2}{ N -2}   \right\|_{L_x^\frac{2( N +2) }{ N +4}  } \|w_n \|_{L_{t,x}^\infty}^\frac4{ N -2}  \|\nabla w_n \|_{L_{t,x}^\infty }
\\
& \lesssim \lambda_n^{- \frac{ N^2 + 2 N - 4 }{2( N +2) } + \frac{ N +2}{2( N -2)} \theta } +
\lambda_n^{- \frac{N(N+4)}{2(N+2)} + \frac{3N - 2}{2(N-2)} \theta }.
\end{align*}
Using \eqref{eq6.9}, the Sobolev embedding, and Theorem \ref{th2.3}, we obtain
\begin{align*}
\left\| | z_n|^\frac4{N-2}  \nabla z_n \right\|_{L_t^\infty L_x^\frac{2(N+2)}{N+4}}
& \lesssim \|\nabla z_n \|_{L_t^\infty L_x^2} \|z_n \|_{L_t^\infty L_x^\frac{4(N+2)}{N-2}}^\frac4{N-2}
\\
& \lesssim \| \nabla z_n \|_{L_t^\infty L_x^2}  \left\| |\nabla |^\frac{N(N+6)}{4(N+2)} z_n  \right\|_{L_t^\infty L_x^2}^\frac4{N-2}
\\
& \lesssim \lambda_n^{1- \frac{N}2 + \frac52\theta }  \left(T+ \lambda_n^{-2 \theta}  \right)  \left( \lambda_n^{\frac{N(N+6)}{4(N+2)} - \frac{N}2 + \frac52 \theta }
\left(T+ \lambda_n^{- 2 \theta}  \right) \right)^\frac4{N- 2}
\\
& \sim  \left(T+ \lambda_n^{- 2 \theta}  \right)^\frac{N+2}{N-2} \lambda_n^{1- \frac{N}2 + \frac52 \theta } \lambda_n^{\frac4{N-2}  \left( \frac{N(N+6)}{4(N+2)} - \frac{N}2 + \frac52 \theta  \right)}
\sim  \left(T+ \lambda_n^{- 2 \theta}  \right)^\frac{N+2}{N-2} \lambda_n^{ - \frac{N^2 + 2 N- 4}{2(N+2)} + \frac{5(N+2)}{2(N-2)} \theta }.
\end{align*}
Similarly, using also \eqref{eq6.9}, the Sobolev embedding, and Theorem \ref{th2.3}, we obtain
\begin{align*}
& \left\| |z_n|^\frac4{N-2} \nabla ( \chi_n w_n )  \right\|_{L_t^\infty L_x^\frac{2(N+2)}{N+4} }
\\
& \lesssim \|\nabla \chi_n \|_{L_x^N} \|w_n \|_{L_t^\infty L_x^\frac{ 2 N(N+2 )^2}{ (N-2)(N^2 + 2 N - 4)}} \|z_n \|_{L_t^\infty L_x^\frac{2N(N+2)^2}{(N-2)(N^2 + 2 N - 4)}}^\frac4{N-2}
+ \|\nabla w_n \|_{L_t^\infty L_x^2} \|z_n \|_{L_t^\infty L_x^\frac{4(N+2)}{N-2}}^\frac4{N-2}
\\
& \lesssim \|w_n \|_{L_t^\infty L_x^\frac{2N(N+2)^2}{ (N-2)(N^2 + 2N- 4)}}  \|z_n \|_{L_t^\infty L_x^\frac{2N(N+2)^2}{ (N-2)(N^2 + 2 N - 4)}}^\frac4{N- 2}
+ \|z_n \|_{L_t^\infty L_x^\frac{4(N+2)}{ N- 2}}^\frac4{N-2} \\
& \lesssim   \left\| |\nabla |^\frac{2(N^2 + 3 N- 2)}{ (N+2)^2}  w_n  \right\|_{L_t^\infty L_x^2}
\left\|  |\nabla |^\frac{2(N^2 + 3 N - 2)}{ (N+2)^2} z_n  \right\|_{L_t^\infty L_x^2}^\frac4{N-2} +  \left\| |\nabla |^\frac{N(N+6)}{4(N+2)} z_n  \right\|_{L_t^\infty L_x^2}^\frac4{N-2} \\
& \lesssim \lambda_n^{ \theta  \left( \frac{2(N^2 + 3 N - 2)}{(N+2)^2} - 1  \right) }   \left( \lambda_n^{ \frac{2(N^2 + 3N- 2) }{(N+2)^2} - \frac{N}2 + \frac52 \theta } \left(T+ \lambda_n^{- 2\theta}  \right) \right)^\frac4{N - 2}
+  \left( \lambda_n^{ \frac{N(N+6)}{ 4(N+2)} - \frac{N}2 + \frac52 \theta }  \left(T+ \lambda_n^{- 2\theta }  \right)  \right)^\frac4{N-2}  \\
& \lesssim \lambda_n^{- \frac{4(N^3 - 8 N + 8 )}{ 2(N+2)^2 (N-2)} }  \left(T+ \lambda_n^{- 2\theta}  \right)^\frac4{N-2} \lambda_n^{\theta  ( \frac{N^2 + 2 N- 8 }{(N+2)^2} + \frac{10}{N - 2} )}
+  \left(T+ \lambda_n^{- 2\theta}  \right)^\frac4{N-2 } \lambda_n^{\theta \frac{10}{N-2} } \lambda_n^{- \frac{N}{ N+2} },
\end{align*}
and
\begin{align*}
&  \left\| \nabla ( \chi_n w_n) |\chi_n w_n|^\frac{6- N}{ N- 2} z_n  \right\|_{L_t^\infty L_x^\frac{2(N+2)}{N+4} }
\\
& \lesssim  \left\|\nabla \chi_n  \right\|_{L_x^N} \|w_n \|_{L_t^\infty L_x^\frac{ 2N(N+2)^2}{ (N - 2) (N^2 + 2N- 4)}}^\frac4{N-2} \|z_n \|_{L_t^\infty L_x^\frac{2N(N+2)^2}{(N-2) (N^2 + 2 N - 4)}}
+  \|\nabla w_n \|_{L_t^\infty L_x^2} \|w_n \|_{L_t^\infty L_x^\frac{4(N+2)}{N-2}}^\frac{6-N}{N-2} \|z_n \|_{L_t^\infty L_x^\frac{4(N+2)}{N-2}} \\
& \lesssim
\left\| |\nabla |^\frac{2(N^2 + 3N - 2)}{ (N+2)^2}  w_n  \right\|_{L_t^\infty L_x^2}^\frac4{N-2}
\left\| |\nabla |^\frac{2(N^2 + 3N-2)}{ (N+2)^2}  z_n  \right\|_{L_t^\infty L_x^2}
+  \left\| |\nabla |^\frac{N(N+6)}{ 4(N+2)} w_n  \right\|_{L_t^\infty L_x^2}^\frac{6 - N}{ N- 2}
\left\|  |\nabla |^\frac{N(N+6)}{ 4(N+2)} z_n  \right\|_{L_t^\infty L_x^2} \\
& \lesssim  \left( \lambda_n^{ \theta  \left( \frac{2(N^2 + 3 N - 2)}{(N+2)^2} - 1 \right) } \right)^\frac4{N- 2} \lambda_n^{ \frac{2(N^2 + 3N- 2)}{(N+2)^2} - \frac{N}2 + \frac52 \theta }
\left(T+ \lambda_n^{- 2\theta}  \right)
+  \left( \lambda_n^{\theta  \left( \frac{N(N+6)}{4(N+2)} - 1 \right)}  \right)^\frac{6- N}{N-2} \lambda_n^{ \frac{N(N+6)}{ 4(N+2)} - \frac{N}2 + \frac52 \theta}
\left(T+ \lambda_n^{- 2 \theta}  \right) \\
& \sim \lambda_n^{- \frac{N^3 - 8 N + 8 }{ 2(N+2)^2} }  \left(T+ \lambda_n^{- 2\theta}  \right) \lambda_n^{\theta  \left( \frac{2(N^2 + 3N - 2)}{ (N+2)^2} - 1 \right)\frac4{N-2}  + \frac52 \theta}
+ \lambda_n^{\theta \left( \frac{N(N+6)}{4(N+2)} - 1 \right) \frac{6-N}{N-2} + \frac52 \theta }  \left(T+ \lambda_n^{- 2\theta}  \right) \lambda_n^{- \frac{N(N-2)}{ 4(N+2)} }.
\end{align*}
In addition, we have
\begin{align*}
& \left\| |\chi_n w_n|^\frac4{N-2} \nabla z_n  \right\|_{L_t^\infty L_x^\frac{2(N+2)}{ N+4} } \\
& \lesssim \|\nabla z_n \|_{L_t^\infty L_x^2}  \|w_n \|_{L_t^\infty L_x^\frac{4(N+2)}{ N-2}}^\frac4{N-2}
\\
& \lesssim \|\nabla z_n \|_{L_t^\infty L_x^2}  \left\| |\nabla |^\frac{N(N+6)}{4(N+2)} w_n  \right\|_{L_t^\infty L_x^2}^\frac4{N-2}
\\
& \lesssim \lambda_n^{1- \frac{N}2 + \frac52 \theta}  \left(T+ \lambda_n^{- 2 \theta}  \right) \lambda_n^{\theta  \left( \frac{N(N+6)}{ 4 (N+2)} - 1 \right)  \frac4{N-2}}
 \sim  \left(T+ \lambda_n^{- 2\theta} \right) \lambda_n^{- \frac{N-2}2 } \lambda_n^{  \left( \frac52 + \frac4{N- 2}  \left( \frac{N(N+6)}{4(N+2)}- 1 \right) \right) \theta }.
\end{align*}
Putting everything together, we derive
\begin{align*}
A_{3, n}   \to 0, \text{ as } n \to \infty.
\end{align*}
Therefore,
\begin{align*}
\lim\limits_{T\to \infty} \limsup\limits_{n \to \infty} \|e_n \|_{\dot{N}^1 \left( \left[0, \lambda_n^2 T \right]  \times \Omega \right)} = 0,
\end{align*}
which together with \eqref{eq6.15} gives \eqref{eq6.14}.

Using \eqref{eq6.11}, \eqref{eq6.12}, and \eqref{eq6.14}, and invoking Theorem \ref{th2.9}, for $n$ and $T$ sufficiently large, we obtain a global solution $v_n$ to $\CGL_\Omega$ with initial data $v_n(0) = \phi_n$ and
\begin{align*}
\|v_n \|_{L_{t,x}^{\frac{2( N +2)}{ N -2}}   \left( \mathbb{R}_+ \times \Omega \right)} \lesssim 1.
\end{align*}

\end{proof}

Our next result concerns the scenario when the rescaled obstacles $\Omega_n^c$
 (where $\Omega_n = \lambda_n^{- 1}  \left( \Omega - \left\{ x_n \right\} \right)$) are retreating to infinity, which corresponds to Case 3 of Theorem \ref{th5.6v65}.

\begin{proposition}\label{th6.3}

Let $\{x_n \} \subseteq \Omega$ and $\{ \lambda_n \} \subseteq 2^{\mathbb{Z}}$ be such that $\beta_n: =
\frac{d(x_n)}{ \lambda_n} \to \infty$ as $n\to \infty$. Let $\phi \in \dot{H}^1( \mathbb{R}^N  )$ and define
\begin{align*}
\phi_n(x) = \frac1{  \lambda_n^{ \frac{ N -2}2}}  ( \chi_n \phi) \left( \frac{x- x_n}{\lambda_n } \right),
\end{align*}
where $\chi_n (x) = 1- \Theta \left( \frac{ |x|}{ \beta_n } \right)$.
 Then for $n$ sufficiently large, there exists a global solution $v_n$ to $\CGL_\Omega$ with initial data $v_n(0) = \phi_n$ which satisfies
\begin{align*}
\| v_n \|_{L_{t,x}^{ \frac{2( N +2)}{ N -2} } ( \mathbb{R}_+ \times \Omega)} \lesssim 1,
\end{align*}
with the implicit constant depending only on $\|\phi \|_{\dot{H}^1}$.

\end{proposition}

\begin{proof}

As before, we fix $0<\theta\ll 1$ and let $w_n$ and $w_\infty$ be solutions to $\CGL_{\mathbb{R}^N  }$ with initial data
$w_n (0) = P_{ \le \beta_n^\theta }^{\mathbb{R}^N } \phi $, and $w_\infty (0) = \phi$. By Theorem \ref{th1.1v37}, we know that
$w_n$ and $w_\infty$ are global solutions obeying the global spacetime norms. Moreover, arguing as in the proof of Proposition \ref{th6.1} and invoking
 Theorem \ref{th2.9} and Lemma \ref{le2.10}, we see that $w_n$ and $w_\infty$ satisfy the following
\begin{align}\label{eq6.24}
& \|w_n \|_{\dot{S}^1( \mathbb{R}_+ \times \mathbb{R}^N  )} + \|w_\infty \|_{\dot{S}^1(\mathbb{R}_+ \times \mathbb{R}^N  )} \lesssim 1, \\
& \lim\limits_{n \to \infty} \|w_n - w_\infty \|_{\dot{S}^1( \mathbb{R}_+ \times \mathbb{R}^N  )} = 0, \notag\\
& \left\|  |\nabla |^s w_n  \right\|_{\dot{S}^1( \mathbb{R}_+ \times \mathbb{R}^N  )} \lesssim
\beta_n^{s\theta}\qquad \text{ for all } 0 \le s \le 2.\notag
\end{align}
Let $T > 0$ be fixed and determined later. Set
\begin{align*}
\tilde{v}_n(t,x) =
\begin{cases}
\frac1{ \lambda_n^{\frac{ N -2}2} } ( \chi_n w_n) \left( \frac{t}{ \lambda_n^{ 2} } ,  \frac{x- x_n}{ \lambda_n}   \right), 0 \le t \le \lambda_n^2 T,
\\
\\
e^{ \left(t - \lambda_n^2 T \right)  {z} \Delta_\Omega } \tilde{v}_n \left( \lambda_n^2 T, x  \right), t > \lambda_n^2 T.
\end{cases}
\end{align*}
Then $\tilde{v}_n$ has finite space-time norm. Indeed, using a change of variables, the Strichartz inequality, the H\"older inequality, the Sobolev embedding, and \eqref{eq6.24}, we get
\begin{align}\label{eq6.25}
\left\| \tilde{v}_n  \right\|_{L_{t,x}^{ \frac{2( N +2)}{ N -2} } ( \mathbb{R}_+ \times \Omega)}
& \lesssim \|  \chi_n w_n \|_{L_{t,x}^{ \frac{2( N +2)}{ N -2} } ( [0, T] \times \Omega_n )}
+ \| \chi_n w_n (T) \|_{ \dot{H}_D^1( \Omega_n )}  \\
& \lesssim \| w_n \|_{\dot{S}^1( \mathbb{R}_+ \times \mathbb{R}^N  )} + \|\nabla \chi_n \|_{L_x^N } \|w_n (T) \|_{L_x^\frac{2 N }{ N -2}  } + \| \chi_n \|_{L_x^\infty } \|\nabla w_n(T) \|_{L_x^2} \lesssim 1.
\notag
\end{align}
Observe that
\begin{align}\label{eq6.26}
\lim\limits_{T \to \infty} \limsup\limits_{n \to \infty}  \left\| \left( - \Delta_\Omega \right)^\frac12 e^{t  {z} \Delta_\Omega}  \left( \tilde{v}_n(0) - \phi_n  \right) \right\|_{L_t^{ \frac{2( N +2)}{ N -2} } L_x^\frac{2 N ( N +2)}{ N^2 + 4 }} = 0.
\end{align}
In fact, by the Strichartz estimate and a change of variables, we have
\begin{align*}
&  \left\|  \left( - \Delta_\Omega \right)^\frac12 e^{t  {z} \Delta_\Omega }  \left( \tilde{v}_n(0) - \phi_n \right) \right\|_{L_t^{ \frac{2( N +2)}{ N -2} } L_x^\frac{2 N ( N +2) }{ N^2 + 4 } ( \mathbb{R}_+ \times \Omega)}   \\
& \lesssim \left\|  \left( - \Delta_{\Omega_n }  \right)^\frac12  \left( \chi_n P^{\mathbb{R}^N }_{> \beta_n^\theta } \phi  \right) \right\|_{L_x^2 ( \Omega_n )}
\\
& \lesssim \| \nabla \chi_n \|_{L_x^N }  \left\|P_{>
\beta_n^\theta }^{\mathbb{R}^N } \phi   \right\|_{L_x^\frac{2 N }{ N -2}  }
+ \| \chi_n \|_{L_x^\infty}  \left\|\nabla  P_{>  \beta_n^\theta }^{\mathbb{R}^N } \phi   \right\|_{L_x^2} \xrightarrow{n\to \infty} 0.
\end{align*}
Now, we claim that $\tilde{v}_n$ is an approximate solution to $\CGL_\Omega$ in the sense that
\begin{align}\label{eq6.29}
\lim\limits_{T \to \infty} \limsup\limits_{n \to \infty}   \left\|  \left( - \Delta_\Omega \right)^\frac12  \left( \left(   \partial_t - z  \Delta_\Omega \right)\tilde{v}_n +  z   \left|\tilde{v}_n  \right|^\frac4{ N -2}  \tilde{v}_n  \right)  \right\|_{ {N}^0 ( \mathbb{R}_+ \times \Omega ) }  = 0.
\end{align}
We first verify \eqref{eq6.29} for $t \in \left(  \lambda_n^2 T, \infty \right)$.
 In this case \eqref{eq6.29} follows as
\begin{equation*}
\left\| e^{ \left(t - \lambda_n^2 T \right)  {z} \Delta_\Omega } \tilde{v}_n  \left( \lambda_n^2 T \right)  \right\|_{L_{t,x}^{ \frac{2( N +2)}{ N -2} }
 \left( \left( \lambda_n^2 T, \infty \right) \times \Omega \right) }\lesssim  \left\|\tilde{v}_n \left(\lambda_n^2 T \right) \right\|_{\dot{H}_x^1}\lesssim \|w_n(T)\|_{\dot{H}_x^1}\xrightarrow{T\to \infty}0.
\end{equation*}
Next we show \eqref{eq6.29} on the time interval $ \left[0, \lambda_n^2 T \right]$. Direct computation gives
\begin{align*}
&  \left(  \left(   \partial_t - z  \Delta_\Omega \right) \tilde{v}_n
+  z \left|\tilde{v}_n \right|^\frac4{ N -2}  \tilde{v}_n  \right)(t,x)
\\
& =  - \lambda_n^{- \frac{ N +2}2} z \left(  \left( \chi_n - \chi_n^\frac{ N +2}{ N -2}   \right) |w_n|^\frac4{ N -2}  w_n  \right) \left( \lambda_n^{- 2 } t , \lambda_n^{- 1}(x- x_n ) \right)\\
&\quad + 2 \lambda_n^{- \frac{ N +2}2} z    \left( \nabla \chi_n \cdot \nabla w_n \right) \left( \lambda_n^{- 2} t, \lambda_n^{- 1}(x - x_n) \right)
 + \lambda_n^{- \frac{ N +2}2}  z \left( \Delta \chi_n w_n \right) \left( \lambda_n^{- 2}t , \lambda_n^{- 1}( x- x_n ) \right).
\end{align*}
Thus, using a change of variables and Theorem \ref{th2.3}, we obtain
\begin{equation}\label{eq3.28v26}
\begin{aligned}
&  \left\| \left( - \Delta_\Omega \right)^\frac12 \left(  \left(   \partial_t - z  \Delta  \right) \tilde{v}_n %-  \mu
 + z  \left|\tilde{v}_n \right|^\frac4{ N -2} \tilde{v}_n  \right)  \right\|_{  {N}^0 \left( \left[0, \lambda_n^2 T \right] \times \Omega  \right)}  \\
& \lesssim   \left\| \nabla  \left(  \left( \chi_n - \chi_n^\frac{ N +2}{ N -2}   \right) |w_n|^\frac4{ N -2}  w_n  \right) \right\|_{ {N}^0([0, T]\times \Omega_n )}
\\
& + \left\| \nabla  \left( \nabla \chi_n \cdot \nabla w_n  \right)  \right\|_{  {N}^0 ([0, T] \times \Omega_n )}
+ \left\|\nabla  \left( \Delta  \chi_n \cdot w_n  \right) \right\|_{ {N}^0([0,T] \times \Omega_n )}  \\
& : = B_{1,n }  + B_{2, n} ,
\end{aligned}
\end{equation}
where $B_{1,n }$ is the first term on the right hand side of \eqref{eq3.28v26}, and $B_{2,n } $ denotes the second and third term on the right hand side of \eqref{eq3.28v26}.
Using H\"older's inequality, we can estimate $B_{1,n } $ as follows:
\begin{align*}
B_{1,n }  & \lesssim \left\| \left( \chi_n - \chi_n^\frac{ N +2}{ N -2}  \right) |w_n |^\frac4{ N -2}  \nabla w_n  \right\|_{L_{t,x}^\frac{2( N +2)}{ N +4} }
+  \left\|\nabla \chi_n  \left( 1 - \frac{ N +2}{ N -2}  \chi_n^\frac4{ N -2}  \right) |w_n|^\frac{ N +2}{N -2}  \right\|_{L_t^\frac{ N +2}N   L_x^\frac{2 N ( N +2)}{ N^2 + 2N  + 8 }}   \\
& \lesssim \| \nabla w_n \|_{L_{t,x}^\frac{2( N +2) } N  }  \left( \|w_n - w_\infty \|_{L_{t,x}^{ \frac{2( N +2)}{ N -2} }}^\frac4{ N -2}  +  \left\| 1_{|x| \sim
\beta_n} w_\infty \right\|_{L_{t,x}^{ \frac{2( N +2)}{ N -2} }}^\frac4{N -2}   \right)  \\
& \ + \|w_n \|_{L_t^\frac{ N +2}{ N - 2}  L_x^{ \frac{2 N ( N +2)}{( N -2)^2 }}} \| \nabla \chi_n \|_{L_{x}^N  }  \left( \|w_n - w_\infty \|_{L_{t,x}^{ \frac{2( N +2)}{ N -2} }}^\frac4{ N -2} +  \left\| 1_{|x| \sim \beta_n} w_\infty \right\|_{L_{t,x}^{ \frac{2( N +2)}{ N -2} }}^\frac4{ N -2}  \right) \to 0,
\end{align*}
by the dominated convergence theorem and \eqref{eq6.24}. Similarly,
\begin{align*}
B_{2,n } & \lesssim T \left(  \| \Delta \chi_n \|_{L_x^\infty} \| \nabla w_n \|_{L_t^\infty L_x^2}
+ \|\nabla \chi_n \|_{L_x^\infty } \|\Delta w_n \|_{L_t^\infty L_x^2}
+ \|\nabla \Delta \chi_n \|_{L_x^N  } \|w_n \|_{L_t^\infty L_x^\frac{2N }{N -2} } \right)\\
& \lesssim T \left( \beta_n^{- 2} +  \beta_n^{\theta- 1}  \right) \xrightarrow{n\to \infty} 0.
\end{align*}
This completes the proof of \eqref{eq6.29}.

Using \eqref{eq6.25}, \eqref{eq6.26} and \eqref{eq6.29}, and invoking Theorem \ref{th2.9}, for $n$ sufficiently large, we obtain a global solution $v_n$ to $\CGL_\Omega$ with initial data $v_n(0) = \phi_n$ which satisfies
\begin{align*}
\| v_n \|_{L_{t,x}^{ \frac{2( N +2)}{ N -2} } ( \mathbb{R}_+ \times \Omega)} \lesssim 1
\quad\text{ and }\quad
\lim\limits_{T\to \infty } \limsup\limits_{n \to \infty}  \left\|v_n(t) - \tilde{v}_n (t)  \right\|_{\dot{S}^1( \mathbb{R}_+ \times \Omega )} = 0.
\end{align*}
\end{proof}

Our final result treats the situation when the rescaled domains $\Omega_n$ expand to fill a half space, where $\Omega_n : = \lambda_n^{- 1} R_n^{- 1}  \left( \Omega -  \left\{  x_n^* \right\}  \right)$. Here, $x_n^\ast \in \partial \Omega$ be such that $ \left|x_n - x_n^* \right| = d(x_n )$ and $R_n \in SO( N )$ satisfies
$R_n e_N  = \frac{x_n - x_n^*}{ \left|x_n - x_n^* \right|}$.
This case corresponds to Case 4 of Theorem \ref{th5.6v65}.

\begin{proposition}\label{th6.4}
Let $\{\lambda_n \} \subseteq 2^{\mathbb{Z}}$ and $\{ x_n \} \subseteq \Omega$ be such that $\lambda_n \to 0$ and $\beta_n : =
\frac{d(x_n)}{ \lambda_n } \to \beta_\infty > 0$ as $n \to \infty$. Let $\phi \in \dot{H}_D^1( \mathbb{R}_+^N  )$ and define
\begin{align*}
\phi_n(x) =  \frac1{ \lambda_n^{  \frac{N -2}2}}  \phi \left( \frac{R_n^{- 1}  \left(x - x_n^* \right)}{ \lambda_n } \right).
\end{align*}
Then for $n$ sufficiently large there exists a global solution $v_n$ to $\CGL_\Omega$ with $v_n(0) = \phi_n$ satisfying
\begin{align*}
\|v_n \|_{L_{t,x}^{ \frac{2(N +2)}{N -2} } ( \mathbb{R}_+ \times \Omega)} \lesssim 1,
\end{align*}
with the implicit constant depending only on $\|\phi \|_{\dot{H}^1}$.

\end{proposition}

\begin{proof}
As in the proof of previous cases, we first construct global solutions to $\CGL_{\mathbb{R}_+^N }$.

Fix $0<\theta\ll 1$ and let $w_n$ and $w_\infty$ be solutions to $\CGL_{\mathbb{R}_+^N  }$ with initial data $w_n (0) = P^{\mathbb{R}^N_+}_{ \le \lambda_n^{- \theta}}
\phi$ and $w_\infty (0) = \phi$. Then $w_n$ and $w_\infty$ are global solutions and obey
\begin{align}\label{eq4.28v37}
\|w_n \|_{\dot{S}^1( \mathbb{R}_+ \times \mathbb{R}_+^N  )} + \|w_\infty \|_{\dot{S}^1 \left( \mathbb{R}_+ \times \mathbb{R}_+^N  \right)} \lesssim 1,
\end{align}
with the implicit constant depending only on $\| \phi \|_{ \dot{H}^1( \mathbb{R}_+^N  )}$. Indeed, we can interpret such solutions as solutions to $\CGL_{\mathbb{R}^N }$  that are odd under reflection in $\partial \mathbb{R}_+^N  $. Then by Theorem \ref{th1.1v37} these solutions are global and obey \eqref{eq4.28v37}. Moreover, arguing as in the proof of Theorems \ref{th6.1} and \ref{th6.3} and using Theorem \ref{th2.9} and Lemma \ref{le2.10}, we have
\begin{align}\label{eq6.36}
\lim\limits_{n \to \infty} \|w_n - w_\infty \|_{\dot{S}^1( \mathbb{R}_+ \times \mathbb{R}_+^N  )} = 0, \text{ and }
\left\| \left( - \Delta_{\mathbb{R}_+^N  } \right)^\frac{s}2 w_n  \right\|_{L_t^\infty L_x^2 \left( \mathbb{R}_+ \times \mathbb{R}_+^N   \right)} \lesssim \lambda_n^{- \theta ( s- 1)} \text{ for any } 0 \le s \le 2.
\end{align}
Fix $T> 0$ to be chosen later. On the time interval $ \left[0, \lambda_n^2 T \right]$, we embed $w_n$ using a boundary straightening diffeomorphism $\Psi_n$ of a neighborhood of zero in $\Omega_n$ of size $L_n : = \lambda_n^{- 2 \theta}$ into a corresponding neighborhood in $\mathbb{R}_+^N  $. Indeed, it follows from  \cite[Section 2]{DFV-2022} that we can find $\Psi_n: \Omega_n \cap \left\{  \left|x^\perp \right| \le L_n \right\} \to \mathbb{R}_+^N $ of the form
\begin{align*}
	\Psi_n (x)  : =   \left( x^\perp, x_N   + \psi_n \left( x^\perp \right) \right),
\end{align*}
where $\psi_n:  \left\{x\in \Omega_n: \left|x^\perp \right|\leq L_n \right\}\to \mathbb{R}$ is a smooth function satisfying
\begin{align}\label{eq6.37}
	\psi_n (0) = 0, \ \nabla \psi_n (0) = 0,\   \left| \nabla \psi_n  \left(x^\perp  \right)  \right| \lesssim \lambda_n^{1 - 2 \theta},\
	\left|\partial^\alpha \psi_n  \left(x^\perp \right) \right| \lesssim \lambda_n^{ |\alpha | - 1},\quad \text{ for all } |\alpha | \ge 2.
\end{align}
Consequently, $\Psi_n$ satisfies
\begin{align}\label{eq6.38}
	\left| \det \left( \partial \Psi_n \right) \right| \sim 1 \quad\text{ and }\quad   \left|\partial \Psi_n \right| \lesssim 1.
\end{align}
Now, we are ready to define the approximate solution. Define a cut-off $\chi_n : \mathbb{R}^N   \to [0, 1]$ via
\begin{align*}
\chi_n (x) : = 1 - \Theta \left( \frac{x}{L_n } \right).
\end{align*}
Let $\tilde{w}_n : = \chi_n w_n$ and set
\begin{align*}
\tilde{v}_n(t,x) : =
\begin{cases}
\frac1{ \lambda_n^{ \frac{N -2}2 } }
\left( \tilde{w}_n  \left( \lambda_n^{-2} t  \right) \circ \Psi_n  \right)
\left( \frac{  R_n^{- 1}  \left( x- x_n^* \right) }{ \lambda_n }\right),
\ 0 \le t \le \lambda_n^2 T, \\
\\
e^{ \left( t- \lambda_n^2 T \right) {z} \Delta_\Omega} \tilde{v}_n  \left( \lambda_n^2 T, x \right),  \qquad\qquad\qquad   t > \lambda_n^2 T.
\end{cases}
\end{align*}
Observe that $\tilde{v}_n$ has finite spacetime norm. Indeed, by the Strichartz inequality, a change of variables, and \eqref{eq6.38}, we have
\begin{align}\label{eq6.39}
 \left\|\tilde{v}_n \right\|_{L_{t,x}^{ \frac{2(N +2)}{N -2} }  \left( \mathbb{R}_+ \times \Omega \right) }
& \lesssim \left\|\tilde{w}_n \circ \Psi_n \right\|_{L_{t,x}^{ \frac{2(N +2)}{N -2} }( \mathbb{R}_+ \times \Omega_n )}
+ \left\|\tilde{w}_n (T) \circ \Psi_n \right\|_{\dot{H}_D^1 ( \Omega_n )}
\\
& \lesssim  \left\| \tilde{w}_n \right\|_{L_{t,x}^{ \frac{2( N +2)}{N -2} } ( \mathbb{R}_+ \times \mathbb{R}_+^N  )}
+  \left\|\tilde{w}_n(T) \right\|_{\dot{H}^1( \mathbb{R}_+^N  )} \lesssim 1. \notag
\end{align}
Next, we prove the following asymptotic agreement of initial data:
\begin{align}\label{eq6.40}
\lim\limits_{T\to \infty} \limsup\limits_{n \to \infty}  \left\| \left( - \Delta_\Omega \right)^\frac12 e^{t  {z} \Delta_\Omega}  \left( \tilde{v}_n(0) - \phi_n \right) \right\|_{L_t^{ \frac{2( N +2)}{ N -2} } L_x^\frac{2N (N +2) }{ N^2 + 4 }
\left( \mathbb{R}_+ \times \Omega \right) } = 0.
\end{align}
By the Strichartz inequality and a change of variables, we have
\begin{align*}
& \left\|  \left(- \Delta_\Omega \right)^\frac12 e^{t  {z} \Delta_\Omega} \left( \tilde{v}_n(0) - \phi_n \right) \right\|_{L_t^{ \frac{2( N +2)}{N -2} } L_x^{\frac{ 2 N ( N +2)}{ N^2 + 4 }} ( \mathbb{R}_+ \times \Omega)}
\\
& \lesssim  \left\|  \left( \chi_n P^{\mathbb{R}^N_+}_{ \le \lambda_n^{- \theta} } \phi \right) \circ \Psi_n - \phi \right\|_{\dot{H}_D^1( \Omega_n)}
\\
& \lesssim  \left\|\nabla  \left(  \left( \chi_n P^{\mathbb{R}^N_+}_{> \lambda_n^{- \theta} }  \phi \right) \circ \Psi_n \right) \right\|_{L_x^2 } +  \left\|\nabla \left( \left( \chi_n \phi \right) \circ \Psi_n - \chi_n \phi \right) \right\|_{L_x^2}
+ \left\|\nabla  \left(  \left( 1 - \chi_n  \right) \phi  \right) \right\|_{L_x^2}.
\end{align*}
As $\lambda_n \to 0$ we have $\left\|\nabla P^{\mathbb{R}^N_+}_{> \lambda_n^{- \theta} } \phi \right\|_{L_x^2} \to 0$ as $n \to \infty$. Thus, using \eqref{eq6.38} we see that the first term converges to 0.
For the second term, we note that $\Psi_n (x) \to x$ in $C^1$. Approximating $\phi$ by $C_c^\infty \left( \mathbb{R}_+^N  \right)$ functions we see that the second term converges to 0.
Finally, by the dominated convergence theorem and the fact that $L_n = \lambda_n^{- 2 \theta} \to \infty$, the last term converges to 0 as well.

Now, we claim that $\tilde{v}_n$ is an approximate solution to $\CGL_\Omega$ in the sense that
\begin{align}\label{eq6.43}
\lim\limits_{T\to \infty } \limsup\limits_{n \to \infty }  \left\|  \left( - \Delta_\Omega \right)^\frac12  \left( \left(  \partial_t  - z  \Delta_\Omega \right) \tilde{v}_n
+  z  \left|\tilde{v}_n \right|^\frac4{ N -2} \tilde{v}_n  \right) \right\|_{ {N}^0 ( \mathbb{R}_+  \times \Omega )} = 0.
\end{align}
We first control the contribution of $ \big[\lambda_n^2 T, \infty \big)$. This reduces to prove
\begin{align*}
\lim\limits_{T\to \infty} \limsup\limits_{n \to \infty}  \left\| e^{ \left( t- \lambda_n^2 T \right)  {z} \Delta_\Omega} \tilde{v}_n \left( \lambda_n^2 T \right)  \right\|_{L_{t,x}^{ \frac{2( N +2)}{N -2} }  \left( \left( \lambda_n^2 T, \infty \right) \times \Omega \right)} = 0,
\end{align*}
which follows since
\[
 \left\| e^{ \left( t- \lambda_n^2 T \right)  {z} \Delta_\Omega} \tilde{v}_n \left( \lambda_n^2 T \right)  \right\|_{L_{t,x}^{ \frac{2( N +2)}{ N -2} }  \left( \left( \lambda_n^2 T, \infty \right) \times \Omega \right)}\lesssim  \left\|\tilde{w}_n(T) \right\|_{\dot{H}_x^1}\lesssim \|w_n(T)\|_{\dot{H}_x^1}\xrightarrow{n,T\to \infty}0.
\]
Next we control the contribution of the time interval $ \left[0, \lambda_n^2 T \right]$ to \eqref{eq6.43}.
Direct computation gives
\begin{align*}
\Delta  \left( \tilde{w}_n \circ \Psi_n \right) =  \left( \partial_k \tilde{w}_n  \circ \Psi_n  \right) \Delta \Psi_n^k
+  \left( \partial_{kl} \tilde{w}_n \circ \Psi_n \right) \partial_j \Psi_n^l \partial_j \Psi_n^k,
\end{align*}
where $\Psi_n^k$ denotes the $k-$th component of $\Psi_n$ and repeated indices are summed.

As $\Psi_n (x) = x +  \left(0, \psi_n  \left( x^\perp  \right) \right)$, we have
\begin{align*}
\Delta \Psi_n^k =  \mathcal{O} \left( \partial^2 \psi_n \right), \partial_j \Psi_n^l = \delta_{jl} + \mathcal{O} \left( \partial \psi_n \right),
\partial_j \Psi_n^l \partial_j \Psi_n^k = \delta_{jl} \delta_{jk} + \mathcal{O} \left( \partial \psi_n \right) + \mathcal{O} \left( \left( \partial \psi_n \right)^2 \right),
\end{align*}
where we use $\mathcal{O}$ to denote a collection of similar terms.
For example, $\mathcal{O}  \left( \partial \psi_n \right)$ contains terms of the form $c_j \partial_{x_j} \psi_n$ for some constants $c_j \in \mathbb{R}$, which may depend on the indices $k$ and $l$ appearing on the left hand side.
Therefore,
\begin{align*}
& \left( \partial_k \tilde{w}_n \circ \Psi_n \right) \Delta \Psi_n^k = \mathcal{O}  \left( \left( \partial \tilde{w}_n \circ \Psi_n  \right) \left( \partial^2 \psi_n \right) \right), \\
&  \left( \partial_{kl} \tilde{w}_n \circ \Psi_n \right) \partial_j \Psi_n^l \partial_j \Psi_n^k
= \Delta \tilde{w}_n \circ \Psi_n + \mathcal{O}  \left( \left( \partial^2 \tilde{w}_n \circ \Psi_n  \right) \left( \partial \psi_n +  \left( \partial \psi_n \right)^2  \right) \right)
\end{align*}
and so
\begin{align*}
& \left(   \partial_t - z  \Delta_{\Omega_n } \right)  \left( \tilde{w}_n \circ \Psi_n  \right)
 + z   \left(  \left|\tilde{w}_n \right|^\frac4{ N -2}  \tilde{w}_n  \right) \circ \Psi_n
\\
& =  \left( \left(   \partial_t - z  \Delta_{\mathbb{R}_+^N  }  \right) \tilde{w}_n
+  z \left|\tilde{w}_n  \right|^\frac4{ N -2}  \tilde{w}_n  \right) \circ \Psi_n
+ \mathcal{O} \left( \left( \partial \tilde{w}_n \circ \Psi_n  \right) \left( \partial^2 \psi_n \right) \right) + \mathcal{O}  \left( \left( \partial^2 \tilde{w}_n \circ \Psi_n \right) \left( \partial \psi_n +  \left( \partial \psi_n  \right)^2  \right) \right).
\end{align*}
By a change of variables and \eqref{eq6.38}, we get
\begin{align}
& \left\|  \left( - \Delta_\Omega \right)^\frac12  \left( \left(   \partial_t - z  \Delta_\Omega  \right) \tilde{v}_n
 + z \left|\tilde{v}_n \right|^\frac4{ N -2}  \tilde{v}_n  \right)  \right\|_{L_t^1 L_x^2  \left(  \left[0, \lambda_n^2 T \right] \times \Omega \right)}
\notag \\
& =  \left\|  \left( - \Delta_{\Omega_n }  \right)^\frac12  \left( \left(   \partial_t - z \Delta_{\Omega_n } \right)  \left( \tilde{w}_n \circ \Psi_n  \right)
 +  z    \left(  \left|\tilde{w}_n  \right|^\frac4{ N -2}
 \tilde{w}_n \right) \circ \Psi_n  \right)  \right\|_{L_t^1 L_x^2 ( [0, T] \times \Omega_n )} \notag
\\
& \lesssim \left\|  \left( - \Delta_{\Omega_n}  \right)^\frac12  \left(  \left(  \left(   \partial_t - z  \Delta_{\mathbb{R}_+^N  }  \right) \tilde{w}_n
+  z  \left| \tilde{w}_n \right|^\frac4{ N -2}  \tilde{w}_n  \right) \circ \Psi_n \right) \right\|_{L_t^1 L_x^2 ( [0, T] \times \Omega_n )} \notag\\
& \quad +  \left\|  \left( - \Delta_{\Omega_n } \right)^\frac12  \left(  \left( \partial \tilde{w}_n \circ \Psi_n  \right) \partial^2 \psi_n \right)  \right\|_{L_t^1 L_x^2 ( [0, T] \times \Omega_n )}
\notag \\
& \quad  +  \left\|  \left( - \Delta_{\Omega_n } \right)^\frac12  \left( \left( \partial^2 \tilde{w}_n \circ \Psi_n  \right)  \left( \partial \psi_n +  \left( \partial \psi_n \right)^2  \right) \right) \right\|_{L^1_t L_x^2 ([0, T ] \times \Omega_n )}
\notag \\
& \lesssim  \left\| \nabla  \left( \left(   \partial_t - z \Delta_{\mathbb{R}_+^N  } \right) \tilde{w}_n %- \mu
 + z   \left|\tilde{w}_n \right|^\frac4{ N -2}  \tilde{w}_n \right)  \right\|_{L_t^1 L_x^2  \left([0, T] \times \mathbb{R}_+^N   \right)} \notag   \\
& \quad  +  \left\| \nabla  \left( \left( \partial \tilde{w}_n \circ \Psi_n  \right) \partial^2 \psi_n \right)  \right\|_{L_t^1 L_x^2 ([0, T] \times \Omega_n )} \notag \\
& \quad  +  \left\|\nabla  \left( \left( \partial^2 \tilde{w}_n \circ \Psi_n \right)  \left( \partial \psi_n +  \left( \partial \psi_n \right)^2 \right) \right) \right\|_{L_t^1 L_x^2 ([0, T] \times \Omega_n)}
\notag \\
& : = C_{1,n}  + C_{2,n }  +C_{3,n } .  \notag
\end{align}
Using \eqref{eq6.36}, \eqref{eq6.37} and \eqref{eq6.38}, we can estimate the last two terms as follows:
\begin{align*}
C_{2, n}  & \lesssim
\left\|  \left( \partial \tilde{w}_n \circ \Psi_n \right) \partial^3 \psi_n  \right\|_{L_t^1 L_x^2 ( [0, T] \times \Omega_n )}
+ \left\| \left( \partial^2 \tilde{w}_n \circ \Psi_n  \right) \partial \Psi_n \partial^2 \psi_n  \right\|_{L_t^1 L_x^2 ([0, T] \times \Omega_n )}
\\
& \lesssim T \lambda_n^2  \left\| \nabla \tilde{w}_n  \right\|_{L_t^\infty L_x^2} + T \lambda_n  \left\| \partial^2 \tilde{w}_n \right\|_{L_t^\infty L_x^2}
\\
& \lesssim T \lambda_n^2  \left( \| \nabla \chi_n \|_{L_x^N } \|w_n \|_{L_t^\infty L_x^\frac{2 N }{ N -2} }  + \|\nabla w_n \|_{L_t^\infty L_x^2} \right)
\\
& \quad + T\lambda_n  \left(  \left\|\partial^2 \chi_n \right\|_{L_x^N  } \|w_n \|_{L_t^\infty L_x^\frac{2 N }{ N -2} } + \| \nabla \chi_n \|_{L_x^\infty}
\| \nabla w_n \|_{L_t^\infty L_x^2} +  \left\|\partial^2 w_n \right\|_{L_t^\infty L_x^2} \right)
\\
& \lesssim T \lambda_n^2 + T \lambda_n  \left( L_n^{-1} + \lambda_n^{- \theta} \right) \xrightarrow{n\to\infty} 0,
\end{align*}
and similarly,
\begin{align*}
C_{3,n }   & \lesssim
\left\|  \left( \partial^2 \tilde{w}_n \circ \Psi_n  \right)  \left( \partial^2 \psi_n   + \partial \psi_n \partial^2 \psi_n \right) \right\|_{L_t^1 L_x^2 ([0, T] \times \Omega_n )}\\
&
\quad +  \left\|  \left( \partial^3 \tilde{w}_n \circ \Psi_n  \right) \left( \partial \Psi_n  \left( \partial \psi_n +  \left( \partial \psi_n \right)^2 \right) \right) \right\|_{L_t^1 L_x^2 ([0, T] \times \Omega_n )} \\
& \lesssim T  \left( \lambda_n + \lambda_n^{2 - 2 \theta} \right)  \left\| \partial^2 \tilde{w}_n \right\|_{L_t^\infty L_x^2}
+ T  \left( \lambda_n^{1 - 2\theta} + \lambda_n^{2 - 4 \theta}  \right)  \left\|\partial^3 \tilde{w}_n  \right\|_{L_t^\infty L_x^2 }  \\
& \lesssim T\lambda_n  \left(L_n^{- 1} + \lambda_n^{ - \theta}  \right) \\
& \quad + T \lambda_n^{1 - 2\theta}  \left(  \left\| \partial^3 \chi_n  \right\|_{L_x^N } \|w_n \|_{L_t^\infty L_x^\frac{2 N }{ N -2} }
+  \left\|\partial^2 \chi_n  \right\|_{L_x^\infty }  \left\|\nabla w_n  \right\|_{L_x^2}
+ \|\nabla \chi_n  \|_{L_x^\infty }  \left\|\partial^2 w_n  \right\|_{L_t^\infty L_x^2} +  \left\|\partial^3 w_n  \right\|_{L_t^\infty L_x^2}  \right)  \\
& \lesssim T \lambda_n  \left( L_n^{- 1} + \lambda_n^{ - \theta}  \right)+ T \lambda_n^{ 1 - 2 \theta}  \left( L_n^{- 2 }+ L_n^{-1} \lambda_n^{- \theta} + \lambda_n^{- 2 \theta} \right) \xrightarrow{n\to\infty} 0.
\end{align*}
Finally, we consider $C_{1,n }$. A direct computation gives
\begin{align*}
 \left(   \partial_t - z  \Delta_{\mathbb{R}_+^N  } \right) \tilde{w}_n
 + z  \left|\tilde{w}_n  \right|^\frac4{ N -2} \tilde{w}_n
=  -  z \left( \chi_n - \chi_n^\frac{ N +2}{ N -2}   \right) |w_n |^\frac4{ N -2}  w_n + 2 z \nabla \chi_n \cdot \nabla w_n +  z \Delta \chi_n w_n.
\end{align*}
We then bound each term as follows.
\begin{align*}
 \left\| \nabla  \left( \Delta \chi_n w_n  \right) \right\|_{L_t^1 L_x^2  \left([0, T] \times \mathbb{R}_+^N   \right)}
& \lesssim T  \left(  \left\|\partial^3 \chi_n  \right\|_{L_x^N  } \|w_n \|_{L_t^\infty L_x^\frac{2 N }{ N -2}  } +  \left\|\partial^2 \chi_n  \right\|_{L_x^\infty }
\|\nabla w_n \|_{L_t^\infty L_x^2}  \right) \\
 & \lesssim T L_n^{- 2} \xrightarrow{n\to\infty} 0, \\
\left\| \nabla  \left( \nabla \chi_n \cdot \nabla w_n  \right)  \right\|_{L_t^1 L_x^2  \left([0, T] \times \mathbb{R}_+^N   \right)}
& \lesssim T  \left(  \left\|\partial^2 \chi_n  \right\|_{L_x^\infty }  \left\|\nabla w_n  \right\|_{L_t^\infty L_x^2}
+ \|\nabla \chi_n \|_{L_x^\infty }
\left\|\partial^2 w_n  \right\|_{L_t^\infty L_x^2} \right)
\\
& \lesssim T  \left( L_n^{- 2} + L_n^{- 1} \lambda_n^{- \theta}  \right) \xrightarrow{n\to\infty} 0.
\end{align*}
Finally, for the first term, we have
\begin{align*}
& \quad \left\| \nabla  \left( \left( \chi_n - \chi_n^\frac{ N +2}{N -2}   \right) |w_n|^\frac4{ N -2}  w_n  \right)  \right\|_{ {N}^0 \left([0, T] \times \mathbb{R}_+^N   \right)}\\
& \lesssim  \left\|  \left( \chi_n - \chi_n^\frac{ N +2}{ N -2} \right) |w_n|^\frac{4 }{ N -2}  \nabla w_n  \right\|_{L_{t,x}^\frac{2( N +2)}{ N +4}
 ([0, T] \times \mathbb{R}_+^N  )}
+  \left\| |w_n |^\frac{ N +2}{ N -2}  \nabla \chi_n  \right\|_{{L_t^2  L_x^\frac{2 N}{ N+2 } } \left([0, T] \times \mathbb{R}_+^N  \right)}   \\
& \lesssim  \left\|w_n 1_{ |x|\sim L_n }  \right\|_{L_{t,x}^{ \frac{2( N +2)}{ N -2} }}^\frac4{ N -2}  \|\nabla w_n \|_{L_{t,x}^\frac{2( N +2)} N }
+ \| \nabla \chi_n \|_{L_x^N  }  \left\| w_n 1_{|x| \sim L_n }  \right\|_{L_{t,x}^{ \frac{2( N +2)}{ N -2}}}^\frac4{ N -2}  \|\nabla w_n \|_{{L_t^\frac{2(N +2)}{ N -2}  L_x^\frac{2 N ( N +2)}{ N^2 +N + 2 }}}  \\
& \lesssim  \left\| 1_{|x| \sim L_n } w_\infty  \right\|_{L_{t,x}^{ \frac{2(N +2)}{ N -2} }}^\frac4{N -2}
+ \|w_\infty - w_n \|_{L_{t,x}^{ \frac{2( N +2)}{ N -2}}}^\frac4{N -2}  \xrightarrow{n\to\infty} 0.
\end{align*}
This completes the proof of \eqref{eq6.43}.

Using \eqref{eq6.39}, \eqref{eq6.40}, and \eqref{eq6.43}, and invoking Theorem \ref{th2.9}, for $n$ large enough we obtain a global solution $v_n$ to $\CGL_\Omega$ with initial data
$v_n(0) = \phi_n$ and $\|v_n \|_{L_{t,x}^{ \frac{2( N +2)}{ N -2} }  \left( \mathbb{R}_+ \times \Omega \right)} \lesssim 1. $
Moreover,
\begin{align*}
\lim\limits_{T\to \infty } \limsup\limits_{n \to \infty}  \left\|v_n(t) - \tilde{v}_n (t)  \right\|_{\dot{S}^1( \mathbb{R}_+ \times \Omega )} = 0.
\end{align*}

\end{proof}

\subsection{Existence of a critical element}

In this subsection, we show the existence of a critical element.
Using the induction on energy argument together with \eqref{eq7.2} and  Theorem \ref{th2.9}, we will prove a compactness result for optimizing sequences of
critical elements. As a preliminary, we need a technical decoupling lemma.
 Before the statement, we first define operators $G_n^j$ on general functions of spacetime. These act on linear solutions in a manner corresponding to the action of $g_n^j $ on initial data in Theorem \ref{th5.6v65}. The exact definition depends on the case to which the index $j$ conforms.
In Cases 1, 2, and 3, we set
\begin{align*}
\left(G_n^j f \right)(t,x) : =  \frac1{ \left( \lambda_n^j \right)^{ \frac{ N -2}2} } f   \left(  \frac{t}{  \left( \lambda_n^j \right)^{2} } ,   \frac{x - x_n^j}{  \lambda_n^j }\right),
\end{align*}
while in Case 4, we define
\begin{align*}
\left(G_n^j f \right)(t,x): =  \frac1{ \left( \lambda_n^j \right)^{ \frac{ N -2}2} }  f \left(  \frac{t}{ \left( \lambda_n^j \right)^{ 2} } ,   \frac{
  \left( R_n^j \right)^{-1}  \left( x -  \left(x_n^j \right)^* \right)} {   \lambda_n^j  } \right).
\end{align*}
Here, the parameters $\lambda_n^j, x_n^j,  \left(x_n^j \right)^*$, and $R_n^j$ are defined as in Theorem \ref{th5.6v65}.

The asymptotic orthogonality condition \eqref{eq5.29v65new} gives rise to asymptotic decoupling of the nonlinear profiles.
\begin{lemma}[Decoupling of nonlinear profiles]\label{le7.3}
For $j \ne k$, we have
\begin{align*}
\lim\limits_{n \to \infty}  \left(  \left\|v_n^j v_n^k  \right\|_{L_{t,x}^\frac{ N +2}{ N -2}  ( \mathbb{R}_+ \times \Omega)} +   \left\|v_n^j \nabla v_n^k  \right\|_{L_{t,x}^\frac{ N +2}{2( N -2)}  ( \mathbb{R}_+ \times \Omega)} +  \left\|\nabla v_n^j \nabla v_n^k  \right\|_{{L_t^\frac{ N +2}{ N -2}  L_x^\frac{ N ( N +2)}{ N^2 + N  + 2 }}( \mathbb{R}_+ \times \Omega)} \right) = 0.
\end{align*}
\end{lemma}

\begin{proof}
Observe that as a consequence of the asymptotic orthogonality condition \eqref{eq5.29v65new}, one has the following asymptotic decoupling (see for instance \cite{Ke}):
\begin{equation}\label{le7.1}
\begin{aligned}
	\left\| G_n^j \psi^j G_n^k \psi^k \right\|_{L_{t,x}^\frac{ N +2}{ N - 2} ( \mathbb{R}_+ \times \mathbb{R}^N  )}
	&+ \left\|G_n^j \psi^j \nabla  \left(G_n^k \psi^k \right) \right\|_{L_{t,x}^\frac{ N +2}{2( N -2)} ( \mathbb{R}_+ \times \mathbb{R}^N  )}\\
	&+ \left\|\nabla  \left(G_n^j \psi^j \right) \nabla \left(G_n^k \psi^k \right)  \right\|_{L_{t,x}^\frac{ N +2}{ N }  ( \mathbb{R}_+ \times \mathbb{R}^N  )}
	\xrightarrow{n\to\infty} 0.
\end{aligned}
\end{equation}
For any $\epsilon > 0$, there exist $N_\epsilon \in \mathbb{N}$ and $\psi_\epsilon^j, \psi_\epsilon^k \in C_c^\infty ( \mathbb{R}_+ \times \mathbb{R}^N  )$, so that
\begin{align}\label{equ:3.60}
\left\| v_n^j - G_n^j \psi_\epsilon^j  \right\|_{\dot{X}^1( \mathbb{R}_+ \times \mathbb{R}^N  )}
+ \left\|v_n^k - G_n^k \psi_\epsilon^k \right\|_{\dot{X}^1( \mathbb{R}_+ \times \mathbb{R}^N  )} < \epsilon.
\end{align}
Thus, we obtain from \eqref{eq7.11} and \eqref{le7.1} that
\begin{align*}
 \left\| v_n^j v_n^k  \right\|_{L_{t,x}^\frac{ N +2}{ N -2 } }
& \le  \left\|v_n^j  \left( v_n^k - G_n^k \psi_\epsilon^k  \right)  \right\|_{L_{t,x}^\frac{ N +2}{ N  - 2} }
 +  \left\| \left( v_n^j - G_n^j \psi_\epsilon^j \right) G_n^k \psi_\epsilon^k \right\|_{L_{t,x}^\frac{ N +2}{ N -2} }
+  \left\| G_n^j \psi_\epsilon^j G_n^k \psi_\epsilon^k  \right\|_{L_{t,x}^\frac{ N +2}{ N -2} } \\
& \lesssim  \left\|v_n^j  \right\|_{\dot{X}^1}  \left\|v_n^k - G_n^k \psi_\epsilon^k  \right\|_{\dot{X}^1}
+  \left\| v_n^j - G_n^j \psi_\epsilon^j  \right\|_{\dot{X}^1}  \left\|\psi_\epsilon^k  \right\|_{\dot{X}^1}
+  \left\| G_n^j \psi_\epsilon^j G_n^k \psi_\epsilon^k  \right\|_{L_{t,x}^\frac{ N +2}{ N -2} }
\\
& \lesssim_{E_c, \delta} \epsilon + o(1), \quad\text{ as } n \to \infty.
\end{align*}
Since $\epsilon > 0$ was arbitrary, this proves the first asymptotic decoupling statement.

The second decoupling statement follows by a similar argument as above and thus we turn to the third assertion, where one needs to take care of the error term. By a similar argument as above, interpolation and \eqref{equ:3.60}, we estimate
\begin{align*}
& \left\| \nabla v_n^j \nabla v_n^k  \right\|_{{L_t^\frac{ N +2}{ N -2}  L_x^\frac{ N ( N +2)}{ N^2 + N  + 2 }}}\\
 &\le \left\|\nabla v_n^j  \left( \nabla v_n^k - \nabla \left( G_n^k \psi_\epsilon^k  \right) \right) \right\|_{{L_t^\frac{ N +2}{ N -2}  L_x^\frac{ N ( N +2)}{ N^2 + N + 2 }}}\\
&\quad + \left\| \left( \nabla v_n^j - \nabla \left(  G_n^j \psi_\epsilon^j  \right) \right) \nabla  \left( G_n^k \psi_\epsilon^k  \right)  \right\|_{{L_t^\frac{ N +2}{ N -2}  L_x^\frac{ N ( N +2) }{ N^2 +  N + 2 }}}+
\left\|\nabla  \left( G_n^j \psi_\epsilon^j  \right) \nabla  \left( G_n^k \psi_\epsilon^k  \right) \right\|_{{L_t^\frac{ N +2}{N -2}  L_x^\frac{ N ( N +2)}{ N^2 + N  + 2 }}} \\
& \lesssim_{E_c, \delta} \epsilon +  \left\|\nabla  \left( G_n^j \psi_\epsilon^j \right) \, \nabla \left(  G_n^k \psi_\epsilon^k \right)   \right\|_{L_{t,x}^\frac{ N +2} N }^{{\frac{ N -2}N}}
\left\|\nabla \psi_\epsilon^j  \right\|_{L_x^2}^{{\frac{3N+2}{N(N+2)}}}   \left\|\nabla \psi_\epsilon^k  \right\|_{L_x^2}^{{\frac{3N+2}{N(N+2)}}}\\
&\lesssim_{E_c, \delta} \epsilon + o(1), \quad\text{ as } n \to \infty,
\end{align*}
where we used \eqref{le7.1} in the last step. As $\epsilon > 0$ was arbitrary, the proof is completed.

\end{proof}

\begin{proposition}[Pre-compactness]\label{pr7.2}
Let $u_n: \big[0, T_n^* \big)  \times \Omega \to \mathbb{C}$ be a sequence of solutions with $E(u_n(0) ) \to E_c$ as $n \to \infty$,
and
\begin{align}\label{eq7.4}
\lim\limits_{n \to \infty} \|u_n \|_{ L_{t,x}^\frac{2( N +2)}{ N  - 2} \left( \big[0, T_n^*  \big) \times \Omega\right)  }  = \infty.
\end{align}
Then there is a sequence of times $t_n \in \big[0, T_n^*  \big) $, such that the sequence $u_n(t_n)$ has a subsequence that converges strongly in $\dot{H}_D^1(\Omega)$ as $n \to \infty$.
\end{proposition}

\begin{proof}	

Applying Theorem \ref{th5.6v65} to the bounded sequence $u_n(0)$ in $\dot{H}_D^1(\Omega)$ and passing to a subsequence if necessary, we obtain
\begin{align}\label{eq7.5}
u_n(0) = \sum\limits_{j =1 }^J \phi_n^j + w_n^J
\end{align}
with the properties stated in Theorem \ref{th5.6v65}. As a particular consequence, the following the energy decoupling holds: for any finite $0 \le J \le J^*$
\begin{align}\label{eq7.6}
\lim\limits_{n \to \infty} \left(E(u_n(0) ) - \sum\limits_{j = 1}^J E \left( \phi_n^j \right) - E \left(w_n^J \right) \right) = 0.
\end{align}

Now, we shall prove that $J^* = 1$.
This means that $w_n^1 \to 0$ in $\dot{H}_D^1( \Omega)$ and the only profile $\phi_n^1$ conforms to Case 1 of Theorem \ref{th5.6v65}. We divide our proofs below in two situations:
\medskip

\textbf{Situation 1.} $\sup\limits_j \limsup\limits_{n \to \infty} E \left( \phi_n^j \right) = E_c$.
\medskip

From the non-triviality of the profiles and the convergence $ \left\| \phi_n^j  \right\|_{\dot{H}_D^1}\to \left\| \phi^j  \right\|_{\dot{H}^1}$, we infer
\begin{align*}
\liminf\limits_{n \to \infty} E \left( \phi_n^j \right) > 0 \quad \text{ for every finite } 1 \le j \le J^*.
\end{align*}
Thus, up to a subsequence, \eqref{eq7.6} implies that there is a single profile in the decomposition \eqref{eq7.5} (namely, $J^* = 1$) and we can write
\begin{align}\label{eq7.7}
u_n(0) = \phi_n + w_n, \text{ with $\lim\limits_{n \to \infty} \|w_n \|_{\dot{H}_D^1} = 0$.}
\end{align}
If $\phi_n$ conforms to Cases 2, 3, or 4 of Theorem \ref{th5.6v65}, then by the Theorems \ref{th6.1}, \ref{th6.3} or \ref{th6.4}, there are global solutions $v_n$ to $\CGL_\Omega$ with data $v_n(0) = \phi_n$ that admit a uniform spacetime bound.

By Theorem \ref{th2.9}, this spacetime bound extends to the solutions $u_n$ if $n\gg 1$, which clearly contradicts \eqref{eq7.4}. Therefore, $\phi_n$ must conform to Case 1 and \eqref{eq7.7} becomes
$u_n(0) = \phi + w_n$, with $\lim\limits_{n \to \infty} \|w_n \|_{\dot{H}_D^1} = 0$. This gives the desired compactness.
\medskip

\textbf{Situation 2.} $\sup\limits_j \limsup\limits_{n \to \infty} E \left( \phi_n^j \right) \le E_c - 2 \delta$ for some $\delta > 0$.
\medskip

Observe that for each finite $J \le J^*$, we have $E \left( \phi_n^j \right) \le E_c - \delta$ for all $1 \le j \le J$ and $n$ sufficiently large.

If $j$ conforms to Case 1, then we let $v^j: I^j \times \Omega \to \mathbb{C}$ to be the maximal-lifespan solution to \eqref{eq:gl} with initial data $v^j(0) = \phi^j$. In this case, we simply take $v_n^j(t,x):= v^j(t,x)$. Then $v_n^j$ is also a solution to \eqref{eq:gl} on the time interval $I_n^j : = I^j$. Note that for all $n\in \mathbb{N}$, we have $0 \in I_n^j$, and
\begin{align}\label{eq7.9}
\left\| v_n^j (0) - \phi_n^j  \right\|_{\dot{H}_D^1} = \left\| v^j (0) - \phi^j  \right\|_{\dot{H}_D^1}=0.
\end{align}
Combining this with $E \left( \phi_n^j \right) \le E_c - \delta$ and \eqref{eq7.2}, we deduce that $v_n^j = v^j$ are global solutions that obey
\begin{align*}
\left\|v^j \right\|_{ L_{t,x}^\frac{2( N +2)}{ N  - 2}  \left( \mathbb{R}_+ \times \Omega\right) }
= \left\|v_n^j  \right\|_{ L_{t,x}^\frac{2( N +2)}{ N  - 2}  \left( \mathbb{R}_+ \times \Omega \right) }
 \le L (E_c - \delta) < \infty.
\end{align*}
This, together with the Strichartz inequality, shows that all Strichartz norms of $v_n^j$ are finite, and in particular so is the $\dot{X}^1$ norm.

The above property allows us to approximate $v_n^j$ in $\dot{X}^1(\mathbb{R}_+\times \Omega)$ by $C_c^\infty( \mathbb{R}_+ \times \mathbb{R}^N  )$ functions. More precisely, for any $\epsilon > 0$, there exist $N_\epsilon^j \in \mathbb{N}$ and $\psi_\epsilon^j \in C_c^\infty( \mathbb{R}_+ \times \mathbb{R}^N  )$ so that for $n \ge N_\epsilon^j$, we have
\begin{align}\label{eq7.10}
\left\|v_n^j - G_n^j \psi_\epsilon^j  \right\|_{\dot{X}^1} < \epsilon.
\end{align}
Specifically, choosing $\tilde{\psi}_\epsilon^j \in C_c^\infty( \mathbb{R}_+ \times \mathbb{R}^N  )$ such that
\begin{align*}
\left\|v^j - \tilde{\psi}_\epsilon^j  \right\|_{\dot{X}^1( \mathbb{R}_+ \times \mathbb{R}^N  )} < \frac\epsilon2,
\end{align*}
and then set
\begin{align*}
\psi_\epsilon^j (t,x) : =  \left( \lambda_\infty^j \right)^\frac{ N -2}2  \tilde{\psi}_\epsilon^j  \left(  \left( \lambda_\infty^j \right)^2 t , \lambda_\infty^j x + x_\infty^j  \right).
\end{align*}

When $j$ conforms to Cases 2, 3, or 4, we apply the nonlinear embedding theorems to construct the nonlinear profiles $v_n^j$.
More precisely, let $v_n^j$ be the global solutions to $\CGL_\Omega$ constructed in Theorems \ref{th6.1}, \ref{th6.3}, or \ref{th6.4}, as appropriate.
In particular, $v_n^j$ satisfies \eqref{eq7.10} and $\sup\limits_{n, j} \left\| v_n^j \right\|_{L_{t,x}^\frac{2( N +2)}{ N -2} ( \mathbb{R}_+ \times \Omega) } < \infty$.

In all cases, we can use \eqref{eq7.3} together with our bounds on the spacetime norms of $v_n^j$ and the finiteness of $E_c$ to deduce
\begin{align}\label{eq7.11}
 \left\| v_n^j  \right\|_{\dot{X}^1( \mathbb{R}_+ \times \Omega)} \lesssim_{E_c, \delta} E \left( \phi_n^j \right)^\frac12 \lesssim_{E_c, \delta} 1.
\end{align}
Combining this with \eqref{eq7.6}, we deduce
\begin{align}\label{eq7.12}
\limsup\limits_{n \to \infty } \sum\limits_{j = 1}^J  \left\|v_n^j  \right\|_{\dot{X}^1( \mathbb{R}_+ \times \Omega)}^2 \lesssim_{E_c, \delta} \limsup\limits_{n \to \infty} \sum\limits_{j = 1}^J E \left( \phi_n^j \right) \lesssim_{E_c, \delta} 1 ,
\end{align}
uniformly for finite $J \le J^*$.

By Lemma \ref{le7.3}, we can bound the sum of the nonlinear profiles in $\dot{X}^1$, as follows:
\begin{align}\label{eq7.13}
\limsup\limits_{n \to \infty}  \left\| \sum\limits_{j = 1}^J v_n^j  \right\|_{\dot{X}^1( \mathbb{R}_+ \times \Omega)} \lesssim_{E_c, \delta} 1
\end{align}
uniformly for finite $J\le J^*$. Indeed, by Young's inequality, \eqref{eq7.11}, \eqref{eq7.12} and Lemma \ref{le7.3},
\begin{align*}
\iint_{\mathbb{R}_+ \times \Omega} \Big| \sum\limits_{j = 1}^J v_n^j(t,x)  \Big|^\frac{2( N +2)}{ N -2}  \,\mathrm{d}x \mathrm{d}t
& \lesssim \sum\limits_{j = 1}^J    \iint_{\mathbb{R}_+\times \Omega}    \left| v_n^j (t,x)  \right|^\frac{2( N +2)}{ N -2} \,\mathrm{d}x \mathrm{d}t  + C_J
\sum\limits_{j \ne k } \iint_{\mathbb{R}_+ \times \Omega}  \left|v_n^j  \right|  \left|v_n^k \right|^\frac{ N +6}{ N -2}\,\mathrm{d}x \mathrm{d}t
\\
& \lesssim_{E_c, \delta} 1 +
C_J \sum\limits_{j \ne k }  \left\|v_n^j v_n^k \right\|_{L_{t,x}^\frac{ N +2}{ N -2} }  \left\|v_n^k \right\|_{L_{t,x}^{ \frac{2( N +2)}{ N -2} }}^\frac8{N -2}
\lesssim_{E_c, \delta} 1 + C_J  o(1), \quad \text{ as } n \to \infty.
\end{align*}
Similarly,
\begin{align*}
& \left\| \sum\limits_{j = 1}^J \nabla v_n^j \right\|_{{L_t^\frac{2(N +2)}{ N -2}  L_x^\frac{2 N  ( N +2) }{ N^2+ N  + 2 }}}^2
 =  \left\|  \left( \sum\limits_{j = 1}^J \nabla v_n^j \right)^2  \right\|_{{L_t^\frac{ N +2}{N -2}  L_x^\frac{ N ( N +2) }{N^2 + N  + 2 }}}
\\
&\qquad \lesssim \sum\limits_{j = 1}^J  \left\|\nabla v_n^j \right\|_{{L_t^\frac{2(N +2)}{ N -2}  L_x^\frac{2 N ( N +2) }{ N^2 + N + 2 }}}^2
+ \sum\limits_{j \ne k}  \left\|\nabla v^j_n \nabla v_n^k \right\|_{{L_t^\frac{ N +2}{ N -2}  L_x^\frac{ N ( N +2) }{ N^2 + N + 2 }}}
\lesssim_{E_c, \delta} 1  + o(1), \quad\text{ as } n \to \infty.
\end{align*}
This completes the proof of \eqref{eq7.13}.

The same argument combined with \eqref{eq7.6} shows that given $\eta > 0$, there exists $J' = J'( \eta)$ such that
\begin{align*}
\limsup\limits_{n \to \infty}  \left\|\sum\limits_{j = J'}^J v_n^j \right\|_{\dot{X}^1( \mathbb{R}_+ \times \Omega)} \le \eta,
\end{align*}
uniformly in $J \ge J'$.

Now we are ready to construct an approximate solution to $\CGL_\Omega$.
For each $n$ and $J$, we define
\begin{align*}
u_n^J: = \sum\limits_{j = 1}^J v_n^j + e^{t {z} \Delta_\Omega} w_n^J.
\end{align*}
It is clear that $u_n^J$ is defined globally in time. In order to apply Theorem \ref{th2.9}, it suffices to verify the following three properties of  $u_n^J$:
\begin{align}\label{eq3.37v54}
  \left\|u_n^J (0) - u_n (0) \right\|_{\dot{H}_D^1 } \to 0,  \quad\text{ as $n \to \infty$ for any $J$, }
\end{align}
\begin{align}\label{eq3.38v54}
\text{ $\limsup\limits_{n \to \infty}  \left\|u_n^J  \right\|_{\dot{X}^1( \mathbb{R}_+ \times \Omega)} \lesssim_{E_c, \delta} 1$, \quad uniformly in $J$, }
\end{align}
and
\begin{align}\label{eq3.39v54}
\text{ $\lim\limits_{J \to \infty} \limsup\limits_{n \to \infty}  \left\|  \left(  \partial_t -  z  \Delta_\Omega \right) u_n^J
 + z \left|u_n^J \right|^\frac4{ N -2}  u_n^J \right\|_{\dot{N}^1( \mathbb{R}_+ \times \Omega )} = 0.$}
\end{align}
These properties imply that for sufficiently large $n$ and $J$, $u_n^J$ is an approximate solution to \eqref{eq:gl} with finite spacetime norm, which asymptotically matches $u_n(0)$ at time $t = 0$. Using Theorem \ref{th2.9}, we infer that for $n$, $J$ sufficiently large, the solution $u_n$ inherits the spacetime bounds of $u_n^J$, thus contradicting \eqref{eq7.4}.

We now turn to the proof of \eqref{eq3.37v54}, \eqref{eq3.38v54}, and \eqref{eq3.39v54}. The first two properties are easy to verify. Indeed, \eqref{eq3.37v54} follows directly from \eqref{eq7.5} and \eqref{eq7.9}, while \eqref{eq3.38v54} follows from \eqref{eq7.13} and Theorem \ref{th2.8}:
\begin{align*}
\limsup\limits_{n \to \infty} \left\|u_n^J  \right\|_{ \dot{X}^1( \mathbb{R}_+ \times \Omega )}
\lesssim \limsup\limits_{n \to \infty}   \left\|\sum\limits_{j = 1}^J v_n^j \right\|_{\dot{X}^1( \mathbb{R}_+ \times \Omega )} + \limsup\limits_{n \to \infty}  \left\|w_n^J \right\|_{ \dot{H}_D^1 }
\lesssim_{E_c, \delta} 1.
\end{align*}

It remains to verify \eqref{eq3.39v54}. Writing $F( u ) =  -  | u |^\frac4{ N -2} u $, we obtain via a direct computation
\begin{align}\label{eq7.15}
\left(  \bar{z} \partial_t  - \Delta_\Omega \right) u_n^J - F \left(u_n^J \right)
& = \sum\limits_{j = 1}^J F \left(v_n^j \right) - F \left(u_n^J \right)
\\
& = \sum\limits_{j = 1}^J F \left(v_n^j \right) - F \left( \sum\limits_{j = 1}^J v_n^j \right) + F \left(u_n^J - e^{t {z} \Delta_\Omega} w_n^J \right) - F \left(u_n^J \right). \notag
\end{align}
Taking the derivative, we have
\begin{align*}
 \left| \nabla  \left( \sum\limits_{j = 1}^J F \left( v_n^j \right) - F \left( \sum\limits_{j = 1}^J v_n^j \right) \right)  \right|
\lesssim_J \sum\limits_{j \ne k}    \left|\nabla v_n^j  \right|  \left|v_n^k  \right|^\frac4{ N -2}
\end{align*}
and hence, using \eqref{eq7.11} and Lemma \ref{le7.3}, we conclude
\begin{align*}
&  \left\| \nabla \left( \sum\limits_{j = 1}^J F \left(v_n^j \right) - F \left( \sum\limits_{j = 1}^J v_n^j \right) \right) \right\|_{ \dot{N}^0 ( \mathbb{R}_+ \times \Omega )} \\
&
 \lesssim_J \sum\limits_{j \ne k }  \left\|  \left|\nabla v_n^j \right| \left|v_n^k \right|^\frac4{ N -2}  \right\|_{L_{t,x}^\frac{2( N +2)}{ N +4} }
 \lesssim_J
\sum\limits_{j \ne k}  \left\|\nabla v_n^j v_n^k \right\|_{L_{t,x}^\frac{ N +2}{2( N -2) }}
 \left\|v_n^k \right\|_{L_{t,x}^{ \frac{2( N +2)}{ N -2} }}^\frac{6- N }{ N -2}  \lesssim_{J, E_c, \delta} o(1),  \quad\text{ as } n \to \infty.
\end{align*}
Thus, we obtain from the above estimate and Theorem \ref{th2.3} that
\begin{align}\label{eq7.16}
\lim\limits_{J \to \infty} \limsup\limits_{n \to \infty}  \left\|\sum\limits_{j = 1}^J F(v_n^j) - F \left( \sum\limits_{j = 1}^J v_n^j \right)  \right\|_{\dot{N}^1( \mathbb{R}_+ \times \Omega )} = 0.
\end{align}
For the second difference in \eqref{eq7.15}, we will show
\begin{align}\label{eq7.17}
\lim\limits_{J \to \infty} \limsup\limits_{n \to \infty} \left\| \nabla  \left( F  \left(u_n^J - e^{t  {z} \Delta_{\Omega}} w_n^J \right) - F \left(u_n^J \right)  \right)  \right\|_{L_{t,x}^\frac{2( N +2)}{ N  + 4}} = 0.
\end{align}
Note that
\begin{align*}
& \quad  \left\| \nabla  \left( F \left(u_n^J - e^{tz \Delta_\Omega} w_n^J \right) - F \left(u_n^J \right)  \right)  \right\|_{L_{t,x}^\frac{2( N +2)}{ N +4} }
\\
& \  \lesssim  \left\|\nabla  \left( e^{tz \Delta_\Omega} w_n^J  \right)  \right\|_{L_{t,x}^\frac{2( N +2)}N  }  \left\|e^{tz \Delta_\Omega } w_n^J  \right\|_{L_{t,x}^\frac{2( N +2)}{ N -2}}^\frac4{ N -2}
+  \left\|\nabla u_n^J \right\|_{L_{t,x}^\frac{2( N +2)} N }  \left\|e^{tz \Delta_\Omega} w_n^J  \right\|_{L_{t,x}^\frac{2( N +2)}{N -2}}^\frac4{N -2}
\\
& \ +  \left\|  \left|u_n^J \right|^\frac4{N -2} \nabla \left(  e^{tz \Delta_\Omega} w_n^J  \right) \right\|_{L_{t,x}^\frac{2( N +2)}{ N +4} }
+  \left\|\nabla u_n^J  \right\|_{L_{t,x}^\frac{2( N +2)} N }  \left\|e^{tz \Delta_\Omega} w_n^J \right\|_{L_{t,x}^\frac{2( N +2)}{ N -2}}
\left\|u_n^J  \right\|_{L_{t,x}^\frac{2( N +2)}{ N -2} }^\frac{6- N }{ N -2}.
\end{align*}
Using \eqref{eq5.25v65new}, Theorem \ref{th2.3}, \eqref{eq3.38v54}, Strichartz inequality, and the fact that $w_n^J$ is bounded in $ \dot{H}_D^1$, we see that \eqref{eq7.17} follows once we are able to establish
\begin{align} \label{eq7.18}
\lim\limits_{J \to \infty} \limsup\limits_{n \to \infty}
 \left\|  \left|u_n^J \right|^\frac4{ N -2} \nabla \left(  e^{tz \Delta_\Omega } w_n^J \right)  \right\|_{L_{t,x}^\frac{2( N +2)}{ N +4}} = 0.
\end{align}
By H\"older's inequality, Theorem \ref{th2.3}, \eqref{eq3.38v54} and the Strichartz inequality, we have
\begin{align*}
& \quad \left\|  \left|u_n^J \right|^\frac4{N -2} \nabla  \left( e^{t z \Delta_\Omega } w_n^J \right)  \right\|_{L_{t,x}^\frac{2( N +2)}{N +4}}
\\
& \lesssim  \left\|u_n^J  \right\|_{L_{t,x}^\frac{2(N +2)}{N -2}}^\frac3{ N -2}  \left\| \nabla \left(  e^{tz \Delta_\Omega} w_n^J  \right)  \right\|_{L_{t,x}^\frac{2( N +2)}N }^\frac{N -3}{N -2} \left\| u_n^J \nabla \left(  e^{t z \Delta_\Omega} w_n^J \right)  \right\|_{L_{t,x}^\frac{N +2}{N -1}}^\frac1{N -2}
\\
& \lesssim  \left\| \sum\limits_{j = 1}^J v_n^j \nabla \left( e^{tz \Delta_\Omega} w_n^J  \right)  \right\|_{L_{t,x}^\frac{ N +2}{N -1}}^\frac1{N -2}
+  \left\|e^{tz \Delta_\Omega} w_n^J  \right\|_{L_{t,x}^\frac{2( N +2)}{ N -2}}^\frac1{ N -2}  \left\|\nabla \left(  e^{tz \Delta_\Omega} w_n^J \right) \right\|_{L_{t,x}^\frac{2( N +2)}{N }}^\frac1{N-2}
\\
& \lesssim  \left\| e^{tz \Delta_\Omega} w_n^J  \right\|_{L_{t,x}^\frac{2( N +2)}{ N -2}}^\frac1{ N -2}
+  \left\| \sum\limits_{j = 1}^J v_n^j \nabla \left( e^{tz \Delta_\Omega} w_n^J  \right)  \right\|_{L_{t,x}^\frac{ N +2}{ N -1}}^\frac1{ N -2}.
\end{align*}
By \eqref{eq5.25v65new}, the contribution of the first term to \eqref{eq7.18} is fine.
Thus, the proof of \eqref{eq7.18} reduces to show that for each $1 \le j \le J'$, it holds
\begin{align}\label{eq7.19}
\lim\limits_{J \to \infty} \limsup\limits_{n \to \infty}  \left\|v_n^j \nabla \left(  e^{t {z} \Delta_\Omega} w_n^J  \right) \right\|_{L_{t,x}^\frac{ N +2}{ N -1 } } = 0.
\end{align}
To this end, fix $ 1 \le j \le J'$. Let $\epsilon > 0$, $\psi_\epsilon^j \in C_c^\infty$ be as in \eqref{eq7.10}, and let $R, T> 0$ be such that $\psi_\epsilon^j$ is supported in the cylinder $[0, T] \times \{ |x| \le R\}$.
Then
\begin{align*}
& \supp  \left( T_n^j \psi_\epsilon^j  \right)
\subseteq \big[ 0, \left( \lambda_n^j \right)^2  T \big] \times  \left\{  \left| x- x_n^j  \right| \le \lambda_n^j R  \right\}
\intertext{ and }
&  \left\|T_n^j \psi_\epsilon^j  \right\|_{L_{t,x}^\infty} \lesssim  \left( \lambda_n^j \right)^{- \frac{ N -2}2}  \left\| \psi_\epsilon^j  \right\|_{L_{t,x}^\infty}.
\end{align*}
If $j$ conforms to Case 4, then we replace $x_n^j$ above by $ \left(x_n^j \right)^*$.
Then, we deduce from Proposition \ref{co2.14} that
\begin{align*}
\left\| T_n^j \psi_\epsilon^j \nabla e^{t  {z} \Delta_\Omega} w_n^J  \right\|_{L_{t,x}^\frac{ N +2}{ N -1 } }
\lesssim T^\frac{31}{180} R^\frac7{45}  \left\| \psi_\epsilon^j  \right\|_{L_{t,x}^\infty}  \left\| e^{t  {z} \Delta_\Omega} w_n^J  \right\|_{L_{t,x}^{\frac{2( N +2)}{ N -2} }}^\frac1{18} \left\|w_n^J  \right\|_{\dot{H}_D^1}^\frac{17}{18}
\lesssim_{\psi_\epsilon^j, E_c}  \left\|e^{t {z} \Delta_\Omega} w_n^J \right\|_{L_{t,x}^{\frac{2( N +2)}{ N -2} }}^\frac1{18}.
\end{align*}
Combining this with \eqref{eq7.10} and using Theorem \ref{th2.3} and the Strichartz inequality, we infer that
\begin{align*}
\left\| v_n^j \nabla e^{t {z} \Delta_\Omega} w_n^J \right\|_{L_{t,x}^\frac{ N +2}{ N -1}}
& \lesssim  \left\|v_n^j - T_n^j \psi_\epsilon^j  \right\|_{\dot{X}^1}  \left\| \nabla e^{t  {z} \Delta_\Omega} w_n^J \right\|_{{L_t^\frac{2(N +2)}{ N -2} L_x^\frac{2 N ( N +2)}{N^2 + N  + 2 }}}
+ C \left( \psi_\epsilon^j, E_c \right)  \left\|e^{t {z} \Delta_\Omega} w_n^J \right\|_{L_{t,x}^{ \frac{2( N +2)}{ N -2} }}^\frac1{18}
\\
& \lesssim \epsilon E_c + C \left( \psi_\epsilon^j, E_c \right)
\left\|e^{t  {z} \Delta_\Omega} w_n^J  \right\|_{L_{t,x}^{ \frac{2( N +2)}{ N -2} }}^\frac1{18}.
\end{align*}
Using \eqref{eq5.25v65new}, we thus obtain
\begin{align*}
\text{ LHS of }\eqref{eq7.19} \lesssim_{E_c} \epsilon.
\end{align*}
As $\epsilon > 0$ was arbitrary, \eqref{eq7.19} follows.

This completes the proof of \eqref{eq7.18} and thus also the proof of \eqref{eq7.17}. The final property \eqref{eq3.39v54} then follows by combining \eqref{eq7.16} and \eqref{eq7.17}. The proof is thus complete.

\end{proof}

As an immediate consequence of the Palais-Smale condition, we obtain that the failure of Theorem \ref{th1.3} implies the existence of a critical element.
\begin{theorem}[Existence of a critical element]\label{th7.4}
Suppose Theorem \ref{th1.3} fails to be true. Then there exist a critical energy $0 < E_c < \infty$ and a global solution $u_c $ to \eqref{eq:gl} with $E(u_c (0) ) = E_c$, which blows up in the sense that
$\|u_c \|_{ L_{t,x}^\frac{2( N +2)}{ N  - 2}  \left( \mathbb{R}_+ \times \Omega \right) }  = \infty$, and $ \left\{u_c (t): t \in \mathbb{R}_+ \right\}$ is pre-compact in $\dot{H}_D^1( \Omega) $. Thus, there exist $C:\mathbb{R}_+ \to \mathbb{R}_+$ such that for any $t \in \mathbb{R}_+$, and $\eta > 0$, we have
\begin{align*}
\int_{|x| \ge C( \eta)} |\nabla u_c (t,x)|^2 \,\mathrm{d}x  \le \eta.
\end{align*}

\end{theorem}

\begin{proof}
If Theorem \ref{th1.3} fails to be true, then there must exist a critical energy $0< E_c < \infty$ and a sequence of solutions $u_n : \big[0, T_n^* \big)  \times \Omega \to \mathbb{C}$ such that
\begin{align}
E(u_n(0) ) &\to E_c \notag
\intertext{ and }
\|u_n \|_{ L_{t,x}^\frac{2( N +2)}{ N - 2}  \left(  \big[0, T_n^* \big) \times \Omega\right)  }   &\to \infty, \quad\text{ as } n \to \infty. \label{eq7.21}
\end{align}
Applying Proposition \ref{pr7.2} and passing to a subsequence, we find $t_n \in  \big[0, T_n^* \big)$ and $\phi \in \dot{H}_D^1( \Omega)$ such that
\begin{align*}
u_n(t_n) \to \phi \quad\text{ in } \dot{H}_D^1( \Omega), \text{ as } n \to \infty .
\end{align*}
In particular, $E( \phi) = E_c$.
We take $u_c : [0, T^* ) \times \Omega \to \mathbb{C}$ to be the maximal-lifespan solution to \eqref{eq:gl} with initial data $u_c (0) = \phi$.
From Theorem \ref{th2.9} and \eqref{eq7.21}, we obtain
\begin{align}\label{eq7.22}
\|u_c \|_{ L_{t,x}^\frac{2( N +2)}{ N  - 2} \left(  \big[0, T^*  \big)  \times \Omega \right)  }  = \infty.
\end{align}
Next, we prove that the orbit of $u_c $ is pre-compact in $\dot{H}_D^1$. For any sequence $ \left\{t_n'  \right\} \subset  \left[0, T^*  \right) $, \eqref{eq7.22} implies $
\|u_c \|_{ L_{t,x}^\frac{2( N +2)}{ N  - 2}  \left(  \left[ t_n' , \infty  \right) \times \Omega \right) }  = \infty$.
Thus, by Proposition \ref{pr7.2}, we see that $u_c (t_n')$ admits a subsequence that converges strongly in $\dot{H}_D^1$.
Therefore, $ \left\{u_c (t): t \in  \big[0, T^* \big) \right\}$ is pre-compact in $\dot{H}_D^1(\Omega)$.

We now use a contradiction argument to show that the solution $u$ is global in time. Suppose $T^*  < \infty$ and let $t_n\to  T^* $. Invoking Proposition \ref{pr7.2} and passing to a subsequence, we find $\phi \in \dot{H}_D^1$ such that $u(t_n) \to \phi $ in $\dot{H}_D^1$.
From the local well-posedness theory, there exist $T = T ( \phi) > 0$ and a unique solution $v: [ 0, T] \times \Omega \to \mathbb{C}$ to \eqref{eq:gl} with initial data $v(0) = \phi$ such that $\|v \|_{ L_{t,x}^\frac{2( N +2)}{ N  - 2} \left( [0, T] \times \Omega \right) }  < \infty$.
By Theorem \ref{th2.9}, for $n$ sufficiently large we find a unique solution $\tilde{u}_n : [0, t_n + T] \times \Omega \to \mathbb{C}$ to \eqref{eq:gl} with data $\tilde{u}_n(t_n) = u_c (t_n)$ and $ \left\| \tilde{u}_n  \right\|_{L_{t,x}^\frac{2( N +2)}{ N - 2} \left( [0, t_n + T] \times \Omega \right) }    < \infty$.
From uniqueness of solutions to \eqref{eq:gl}, we have $\tilde{u}_n = u_c $.
Thus, taking $n$ sufficiently large, we see that $u_c $ can be extended beyond $T^* $, which contradicts the fact that $\big[0, T^* \big) $ is the maximal lifespan of $u_c $.
\end{proof}

\section{Nonexistence of the critical element}\label{se6v23}
In this section, we prove the critical element of $\CGL_\Omega $ does not exist. More precisely, we prove that
\begin{theorem}\label{thm:nonexistence}
The critical element $u_c$ in Theorem \ref{th7.4} does not exist.	
\end{theorem}

The key idea to show Theorem \ref{thm:nonexistence} is a Morawetz type inequality. We borrow the idea of Morawetz inequality of $\NLS_{\Omega}$ and rewrite the $\CGL_{\Omega}$ in the form of $\NLS_{\Omega}$:
\begin{align}\label{eq:gl_e}
    iu_t + \Im{z} \Delta u =  \Im{z} |u|^{\frac{4}{ N -2}}u + \mathcal{D},
\end{align}
where $\mathcal{D} = i \Re{z}\Delta u -  i\Re{z} |u|^{\frac{4}{ N -2}}u$ is the dissipation of $\CGL$.
We use \eqref{eq:gl_e} to write
    \begin{align}\label{eq:moment}
    \begin{aligned}
           \partial_t \Im { \left(\Bar{u} \nabla u \right)}
           & = -2\Im{z} \sum\limits_{j =1}^N \partial_j\Re{ \left( \nabla  u \partial_j \Bar{u} \right)} + \frac{\Im{z}}{2}\Delta  \nabla  \left(|u|^2 \right)
     + \frac{2( N -2)}{ N }\Im{z}  \nabla  \left(|u|^{2 N /( N -2)} \right) \\
           &\quad + \Re{z}\Im{ \left( \nabla u \Delta \Bar{u} \right)}
           - \Re{z}\Im{ \left( \nabla u |u|^{\frac{4}{ N -2}}\Bar{u} \right)} +\Re{z}\Im{
           \left(\Delta  \nabla u \Bar{u} \right)} - \Re{z}\Im{ \left(\Bar{u}  \nabla  \left(|u|^{\frac{4}{ N -2}}u \right) \right)}.
    \end{aligned}
    \end{align}
In below, we first establish a Morawetz inequality for $\CGL_{\Omega}$.
\begin{lemma}\label{lem:morawetz}
Let $u$ be a solution to $\CGL_{\Omega}$. Then for any $A\ge 1$ with $A\sqrt{|I|} \ge \diam(\Om^c)$,
\begin{align}\label{est:morawetz}
  \begin{aligned}
    \int_I \int_{\Om} \frac{|u(t, x)|^{\frac{2 N }{ N -2}}}{|x|} \,\mathrm{d}x\mathrm{d}t \lesssim A\sqrt{|I|},
  \end{aligned}
\end{align}
where the implicit constant depends only on $E(u(0) )$.
\end{lemma}

\begin{proof}
  Let $\phi \in C_0^\infty (\mathbb{R}^N )$ be a radial function such that $\phi(x) = 1$ for $|x| \le 1$ and $\phi(x) = 0$ when $|x| \ge 2$. Let $R \ge \text{diam} \left(\Om^c \right)$ and define $a(x) := |x| \phi \left(x/R \right)$. Observe that for $|x| > R$,
  \begin{align*}
    |\nabla a| \lesssim 1, \  \left| \nabla^2 a \right| \lesssim \frac1R, \  \left|\Delta^2 a \right| \lesssim \frac1{R^3},
  \end{align*}
  and for $0 < |x| \le R$,
  \begin{align*}
    \nabla a (x) = \frac{x}{|x|}, \  \nabla^2 a \ \text{ is positive definite}, \  \Delta^2 a < 0.
  \end{align*}
By \eqref{eq:moment}, we have
  \begin{align}\nonumber
      \partial_t \int_\Om\Im { \left(\Bar{u}  \nabla u \right)}\cdot  \nabla a \,\mathrm{d}x
      & = -2\int_\Om \left( \Im{z}\partial_j\Re{ \left( \nabla  u \partial_j \Bar{u} \right)}
      +\frac{\Im{z}}{2}\Delta  \nabla  \left(|u|^2 \right)
      + \frac{2( N -2)}{ N }\Im{z}  \nabla  \left(|u|^{\frac{2 N }{ N -2}} \right) \right) \cdot  \nabla a \,\mathrm{d}x \\\nonumber
           &\quad + \int_\Om \left(\Re{z}\Im{ \left( \nabla u \Delta \Bar{u} \right)}
           - \Re{z}\Im{ \left( \nabla u |u|^{\frac{4}{ N -2}}\Bar{u} \right)} +\Re{z}\Im{ \left(\Delta  \nabla u \Bar{u} \right)} \right) \cdot  \nabla a \,\mathrm{d}x \\\label{eq:momentum_1}
           &\quad - \int_\Om  \Re{z}\Im{ \left(\Bar{u}  \nabla  \left(|u|^{\frac{4}{ N -2}}u \right) \right)} \cdot  \nabla a \,\mathrm{d}x .
  \end{align}
  By applying the Sobolev embedding and usual trace theory, the first term on the right-hand side of \eqref{eq:momentum_1} is bounded below by
  \begin{align}\label{est:dispersive}
  \begin{aligned}
    &\quad-2\int_\Om \left( \Im{z}\partial_j\Re{ \left( \nabla  u \partial_j \Bar{u} \right)} +\frac{\Im{z}}{2}\Delta  \nabla  \left(|u|^2 \right)
    + \frac{2( N -2)}{ N }\Im{z}  \nabla  \left(|u|^{\frac{2 N }{ N -2}} \right) \right) \cdot  \nabla a  \,\mathrm{d}x \\
    &\gtrsim \int_{|x| \le R} \frac{|u|^{2 N /( N -2)}}{|x|} \,\mathrm{d}x  - \frac1R  \left(\|u\|_{\dot{H}^1_x}^2 + \|u\|_{L^{\frac{2 N }{ N -2}}_x}^2 \right).
  \end{aligned}
  \end{align}
For the second term on the right-hand side of \eqref{eq:momentum_1}, we may use integration by parts to obtain
  \begin{align}\label{eq:momentum_2}
    \begin{aligned}
      & \int_\Om \left(\Re{z}\Im{ \left( \nabla u \Delta \Bar{u} \right)}
      - \Re{z} \Im{ \left( \nabla u |u|^{\frac{4}{ N -2}}\Bar{u} \right)} +\Re{z}\Im{ \left(\Delta  \nabla u \Bar{u} \right)} \right) \cdot  \nabla a \,\mathrm{d}x
     -  \int_\Om  \Re{z}\Im{ \left(\Bar{u}  \nabla  \left(|u|^{\frac{4}{ N -2}}u \right) \right)} \cdot  \nabla a \,\mathrm{d}x \\
      & = 2\int_\Om \left(\Re{z}\Im{ \left(Du \Delta \Bar{u} \right)}
      - \Re{z}\Im{ \left( \nabla u |u|^{\frac{4}{ N -2}}\Bar{u} \right)} \right) \cdot  \nabla a \,\mathrm{d}x
     -\int_\Om \left(\Re{z}\Im{ \left(\Delta u \Bar{u} \right)}
      - \Re{z}\Im{ \left(\Bar{u} |u|^{\frac{4}{ N -2}}u \right)} \right)\cdot \Delta a \,\mathrm{d}x .
    \end{aligned}
  \end{align}
To estimate \eqref{eq:momentum_2}, we choose $\varepsilon>0$ small (to be determined later), and
define $u_h = P^\Omega_{\ge \varepsilon} u$ as well as $u_{l} = P^\Omega_{< \varepsilon} u$. Then the Sobolev's embedding theory gives
     \begin{align}\label{est:fre_1}
       \|u_{l} \|_{L^\infty_t \dot{H}^1_x(I \times \Om)}
       + \|u_{l} \|_{L^\infty_t L^{2 N /( N -2)}_x(I \times \Om)} \lesssim \varepsilon_1.
     \end{align}
Combining \eqref{est:fre_1} with Lemma \ref{lea.6v65}, we infer
     \begin{align}\label{est:fre_3}
     \begin{aligned}
       \|u_{h} \|_{L^\infty_t L^2_x(I \times \Om)}
       & \lesssim \left\|P^{\Omega}_{\varepsilon\le \cdot < \varepsilon^{-1}} u \right\|_{L^\infty_t L^2_x(I \times \Om)} + \left\|P^{\Omega}_{>\varepsilon^{-1}}u \right\|_{L^\infty_t L^2_x(I \times \Om)} \\
       &\lesssim \frac1{\varepsilon} \left\|P^{\Omega}_{\varepsilon_1\le \cdot < \varepsilon^{-1}} u \right\|_{L^\infty_t \dot{H}^1_x(I \times \Om)} + \varepsilon \left\|P^{\Omega}_{>\varepsilon_1^{-1}}u \right\|_{L^\infty_t \dot{H}^1_x(I \times \Om)}
       \lesssim \varepsilon + \frac1{\varepsilon}.
     \end{aligned}
     \end{align}
     Furthermore, we have
     \begin{align}\label{est:fre_4}
       \begin{aligned}
         \|u_{h} \|_{L^2_t \dot{H}^1_x(I \times \Om)} & \lesssim \left\|P^{\Omega}_{\varepsilon\le \cdot < \varepsilon^{-1}} u \right\|_{L^2_t \dot{H}^1_x(I \times \Om)} + \left\|P^{\Omega}_{>\varepsilon^{-1}}u \right\|_{L^2_t \dot{H}^1_x(I \times \Om)} \\
       &\lesssim \varepsilon^{-1}\sqrt{|I|} + \varepsilon \left\|P^{\Omega}_{>\varepsilon^{-1}}u \right\|_{L^2_t \dot{H}^2_x(I \times \Om)}
       \lesssim \varepsilon + \frac{\sqrt{|I|}}{\varepsilon}.
       \end{aligned}
     \end{align}
 Now, we write \eqref{eq:momentum_2} differently as
     \begin{align*}
    \begin{aligned}
      & \int_\Om \left(\Re{z}\Im{ \left( \nabla u \Delta \Bar{u} \right)}
      - \Re{z}\Im{ \left( \nabla u |u|^{\frac{4}{ N -2}}\Bar{u} \right)} +\Re{z}\Im{ \left(\Delta  \nabla u \Bar{u} \right)} \right) \cdot  \nabla a \,\mathrm{d}x
 -  \int_\Om \Re{z}\Im{ \left(\Bar{u}  \nabla  \left(|u|^{\frac{4}{ N -2}}u \right) \right)} \cdot  \nabla a \,\mathrm{d}x \\
      & = 2\int_\Om \left(\Re{z}\Im{ \left( \left( \nabla u_l +  \nabla u_h \right) \Delta \Bar{u} \right)}
      -  \Re{z}\Im{ \left( \nabla u |u|^{\frac{4}{ N -2}} \left(\Bar{u}_l + \Bar{u}_h \right) \right)} \right) \cdot  \nabla a \,\mathrm{d}x \\
      & \quad +\int_\Om \Re{z}\Im{ \left( \nabla  u \left(\Bar{u}_l + \Bar{u}_h \right) \right)}\cdot \Delta  \nabla a \,\mathrm{d}x .
    \end{aligned}
  \end{align*}
  By the H\"older's inequality, Sobolev's embedding and \eqref{est:fre_1}, \eqref{est:fre_3}, we deduce
  \begin{align}\label{est:disspation_1}
    \begin{aligned}
      &\int_I\int_\Om \left(\Re{z}\Im{ \left( \left( \nabla u_l +  \nabla u_h \right) \Delta \Bar{u} \right)}
      - \Re{z}\Im{ \left ( \nabla u |u|^{\frac{4}{ N -2}} \left(\Bar{u}_l + \Bar{u}_h \right) \right) } \right) \cdot  \nabla a \,\mathrm{d}x \mathrm{d}t \\
      &\lesssim \int_I\left(\|u_l\|_{\dot{H}^1_x} + \|u_h\|_{\dot{H}^1_x}\right) \|\Delta u\|_{L^2_x} + \| \nabla u\|_{L^{2 N /( N -2)}_x} \left\||u|^{4/( N -2)} \right\|_{L^{ N /2}_x} \left( \|u_l \|_{L^{2 N /( N -2)}_x} + \|u_h\|_{L^{2 N /( N -2)}_x} \right) \,\mathrm{d}t \\
      &\lesssim \varepsilon (1+\sqrt{|I|}) + \frac{\sqrt{|I|}}{\varepsilon}.
    \end{aligned}
  \end{align}
  Using the fact $|x|\ge \text{diam} \left(\Om^c \right)$ and \eqref{est:fre_1}, \eqref{est:fre_3}, we know
  \begin{align}\label{est:disspation_2}
    \begin{aligned}
     &  \int_I\int_\Om \Re{z}\Im{ \left( \nabla  u  \left(\Bar{u}_l + \Bar{u}_h \right) \right)}\cdot \Delta  \nabla a  \,\mathrm{d}x \mathrm{d}t \\
      & \lesssim \int_I\int_\Om | \nabla u| (|u_h| + |u_l|)\frac1{|x|^2} \,\mathrm{d}x \mathrm{d}t \\
      & \lesssim \int_I \| \nabla u\|_{L^2_x} \|u_h\|_{L^2_x} + \| \nabla u\|_{L^{2 N /( N -2)}_x} \|u_l\|_{L^{2 N /( N -2)}_x} \,\mathrm{d}t
      \lesssim \varepsilon + \frac{\sqrt{|I|}}{\varepsilon}.
    \end{aligned}
  \end{align}
Combining the estimates \eqref{est:dispersive}, \eqref{est:disspation_1} with \eqref{est:disspation_2}, and then integrating \eqref{eq:momentum_1} over $I$, we thus conclude
\begin{align*}
  \int_I\int_\Om \frac{|u(x)|^{2 N /( N -2)}}{|x|} \,\mathrm{d}x \mathrm{d}t  \lesssim R + \frac{|I|}{R} + \frac1{\varepsilon}\sqrt{|I|} ,
\end{align*}
where we used \eqref{est:fre_4} and $\varepsilon<1$. Our desired estimate \eqref{est:morawetz} follows by choosing  $\varepsilon \in \left(\frac1{A}, 1\right) $ and $R = A\sqrt{|I|}$.
\end{proof}

\begin{proof}[Proof of Theorem \ref{thm:nonexistence}]
Suppose that Theorem \ref{thm:nonexistence} does not hold. Then by Theorem \ref{th7.4}, there exists a critical element $u_c$ to $\CGL_{\Omega}$ that is global in time, such that $ \left\{u_c( t) : t \in \mathbb{R}_+  \right\}$ is pre-compact in $\dot{H}_D^1(\Omega)$. Moreover,  for some $ R>1 $ large enough, it obeys
    \begin{equation*}
1 \lesssim \int_{\Om \cap \{|x|\le R\}} |u_c (t,x)|^{\frac{2 N }{ N -2}} \,\mathrm{d}x\quad \text{uniformly with respect to}\ t\in  \mathbb{R}_+. %[0, \infty).
    \end{equation*}
    Integrating the previous estimate over time on $I$ with $|I|\ge 1$, we obtain
    \begin{equation}\label{est:n_1}
        |I| \lesssim R\int_I\int_{\Om \cap \{|x|\le R\}} \frac{|u_c (t,x)|^{\frac{2 N }{N -2}}}{|x|} \,\mathrm{d}x\mathrm{d}t \lesssim R\int_I\int_{\Om \cap \{|x|\le R\sqrt{|I|}\}} \frac{|u_c (t,x)|^{\frac{2 N }{ N -2}}}{|x|} \,\mathrm{d}x\mathrm{d}t.
    \end{equation}
On the other hand, the Morawetz inequality in Lemma \ref{lem:morawetz} tells us that
    \begin{equation}\label{est:n_2}
        \int_I\int_{\Om \cap \{|x|\le R\sqrt{|I|}\}} \frac{|u_c (t,x)|^{\frac{2 N }{ N -2}}}{|x|} \,\mathrm{d}x\mathrm{d}t \lesssim R\sqrt{|I|},
    \end{equation}
    with the implicit constant depending only on $E(u_c (0) )$. So \eqref{est:n_1} and \eqref{est:n_2} imply
    \begin{equation*}
|I|\lesssim R^2 \sqrt{|I|},
    \end{equation*}
    which is a contradiction if $|I|$ is chosen to be sufficiently large.
\end{proof}

\begin{proof}[Proof of Theorem \ref{th1.3}]
The global well-posedness follows immediately from Theorem \ref{th7.4} and Theorem \ref{thm:nonexistence}. The proof for the decay follows from the argument in \cite[Section 6]{CGZ}.
\end{proof}

\section{Global weak solutions and weak-strong uniqueness of NLS}

Based on the global well-posedness of energy-critical defocusing CGL equation in the strictly convex exterior domain, we shall show the existence of global weak solutions to energy-critical defocusing NLS equation in the same setting. As pointed out in the introduction, the existence of strong solution to energy-critical defocusing NLS equation in the exterior domain in three dimension $N=3$ was established in \cite{KVZ1}. However, due to technical restrictions, the existence of global strong solutions in higher dimensions $N\ge 4$ remains open. In this section, we shall prove the existence of weak solutions and the weak-strong uniqueness of defocusing NLS in exterior domains for $N = 3, 4, 5$.

\subsection{Global exsitence of weak solutions}
\begin{proof}[Proof of Corollary \ref{cor:nls}, part 1)]
		By Theorem \ref{th1.3}, there are global solutions $u_n$ of
		\begin{align}\label{eq7.1v24}
			\begin{cases}
				&  \pa_tu_n - z_n \Delta u_n=-z_n|u_n|^{\f {4}{N -2}} u_n\quad \text{in } \ \mathbb{R}_+\times\Omega,\\
                & u_n(0, x) = v_0 (x)  \text{ in } \Omega,\\
                & u_n(t,x)|_{x\in \partial \Omega} = 0,
			\end{cases}
		\end{align}
		where $|z_n| = 1$ with $\Re z_n \ge 0$ and $z_n\to i$ as $n \to \infty$.
		From the energy estimates, we have $u_n \in L^\infty_t  \dot{H}^1_D $ and $\pa_tu_n\in L^\infty_t H^{-1}_x$, with the uniform bound
		\begin{align}\label{est:uniform}
			\|u_n \|_{L^\infty_t \dot{H}^1_D} + \left\|\pa_tu_n \right\|_{L^\infty_t H^{-1}_x} \le C \left(\|u_0\|_{\dot{H}^1_D}\right).
		\end{align}
		Then there exists $u \in L^\infty_t \dot{H}^1_D$ such that, up to an extraction of subsequence, we have
		\begin{align}\label{eq:weak_star}
			\begin{aligned}
				&(1)\  u_n(t) \rightharpoonup u(t)\ \text{weakly}\ \text{in}\ \dot{H}^1_D ,\ \text{ as } n  \to \infty \ \text{for every} \ t. \\
				& (2) \ \text{For every } t\in \mathbb{R}_+, \
				\  u_{n}(t,x)  \to u (t,x) \  \text{as} \ n \to \infty
				\text{ for almost every  } \ x\in \Om,		\\
				& (3) \  u_n(t, x) \to u(t, x) \ \text{ as } n  \to \infty \ \text{for almost every} \ (t, x) \in \mathbb{R}_+ \times \Om.
			\end{aligned}
		\end{align}
		Testing the equation \eqref{eq7.1v24} by $\varphi \in C_c^\infty(\Om )$ and $\psi \in  C^\infty_c( \mathbb{R}_+)$, we get
		\begin{align*}
			-\int_0^\infty \left( \int_{\Om } \bar{z}_n   u_n  \varphi  \,\mathrm{d}x \right) \psi_t\mathrm{d}t + \int_0^\infty \left(\int_{\Om } \nabla u_n  \cdot \nabla \varphi \,\mathrm{d}x\right) \psi \mathrm{d}t = \int_0^\infty \left(\int_{\Om } |u_n |^{4/(N-2) } u_n  \varphi \,\mathrm{d}x\right) \psi\mathrm{d}t .
		\end{align*}
		Sending $n\to \infty$ and using the weak convergence in \eqref{eq:weak_star}, we infer that
		\begin{align*}
			\int_0^\infty \left(\int_{\Om } \bar{z}_n   u_n  \varphi \,\mathrm{d}x\right) \psi_t \mathrm{d}t \xrightarrow{n\to \infty} \int_0^\infty \left(\int_{\Om } i u \varphi \,\mathrm{d}x \right) \psi_t\mathrm{d}t,
		\end{align*}
		and
		\begin{align*}
			\int_0^\infty \left(\int_{\Om } \nabla u_n  \cdot \nabla \varphi \,\mathrm{d}x \right) \psi\mathrm{d}t \xrightarrow{n\to \infty} \int_0^\infty \left(\int_{\Om } \nabla u \cdot \nabla \varphi \,\mathrm{d}x \right) \psi\mathrm{d}t.
		\end{align*}
		Note that $|u_n |^{4/(N-2) } u_n  \varphi\psi$ is compactly supported and thus is bounded in $L^\infty_t L^{ \frac{2N}{N+2}}_x$. Then by the almost everywhere convergence in \eqref{eq:weak_star}, we get that $|u_n |^{4/(N-2) } u_n  \varphi \psi \to |u|^{4/(N-2) } u  \varphi \psi$ for almost every $(t, x)$, which  guarantees that $|u_n |^{4/(N-2) } u_n  \varphi \psi \to |u|^{4/(N-2) } u  \varphi \psi$ in $L^1_{t,x}$ as $n\to \infty$.
		
		Consequently, by taking $n \to \infty $, we obtain
		\begin{align*}
			\int_0^\infty\int_{\Om } i u \varphi \psi_t \,\mathrm{d}x \mathrm{d}t + \int_0^\infty\int_{\Om } \nabla u \cdot \nabla \varphi \psi\,\mathrm{d}x \mathrm{d}t  = -\int_0^\infty\int_{\Om } |u|^{4/(N-2) } u \varphi \psi\,\mathrm{d}x \mathrm{d}t .
		\end{align*}
		Hence $v$ satisfies
		\begin{equation}\label{eq7.2v87}
			iu_t +  \Delta u = - |u|^{4/(N-2) } u
		\end{equation}
		in the sense of tempered distribution.
		
		By the lower semi-continuity of $L^2$-norm and the Fatou's lemma, we have
		\begin{align} \label{est:energy_nondecreas}
			\begin{aligned}
				E(u(t))
				\le \liminf_{n \to\infty} E \left(u_{n }(t) \right) \le E(u_0).
			\end{aligned}
		\end{align}
		For any fixed $T>0$, $u \in L^2_t \left([0, T] ; \dot{H}^1_D \left(\Om\right) \right)$ and $\pa_tu \in L^2_t \left([0, T] ; H^{-1}_x \left(\Om \right) \right)$. By a standard limiting argument, we have
		\begin{align}
			\|u(t)\|_{L^2}^2 = \|u(s)\|_{L^2}^2 + \int_s^t \left< \pa_\tau u (\tau), u(\tau)\right>_{\dot{H}^1_D, H^{-1}} \,\mathrm{d}\tau, \quad 0\le s \le t \le T,
		\end{align}
		which implies $u \in C \left([0, T]; L^2 \left(\Om \right) \right)$.
		
		Multiplying $iu$ on both sides of \eqref{eq7.2v87} and integrating on $\Om$, we obtain the mass conservation $M(u(t))= M(u_0)$ as desired.

		Thus, we have obtained a global weak solution of \eqref{eq:nls}
		on $[0, \infty)$. By time reflection, it is easy to extend the solution $u$ on $[0, \infty)$ to obtain a global weak solution $v$ on $\mathbb{R}$.
	\end{proof}
	
	\subsection{Weak-strong uniqueness}
	 First, we show that if \eqref{eq:nls} admits a smooth strong solution, then the weak solution is trapped around the strong solution in the following sense.
	\begin{proposition}
		\label{thm:stability}
	For $N \in \{3, 4, 5 \}$, suppose that $\tilde{u} \in  C^0_t \dot{H}^2_D \left(\mathbb{R} \times \Om \right)$ is a strong solution to defocusing NLS \eqref{eq:nls}. Let ${u}$ be the global weak solution to \eqref{eq:nls} with initial data $u_0$ satisfying the energy inequality
		\begin{align}
			E(u(t)) \le E(u_0), \quad \text{for all}\  t\in \mathbb{R}.
		\end{align}			
		Fix any $T>0$ and let $w = u-\tilde{u} $. Then there exists a constant $C= C  \left(\tilde{u}, T \right)>0$ such that
		\begin{align}
			\|w(t)\|_{\dot{H}^1_D}^2 \le C  \|w(0)\|_{\dot{H}^1_D}^2
		\end{align}
		uniformly for all $ t \in (-T, T)$.
	\end{proposition}
	\begin{proof}
		Without loss of generality, we only consider the forward time case: $t\ge 0$. Direct computation shows that $w$ satisfies the equation
		\begin{align}\label{eq:perturbation}
			\begin{cases}
				iw_t + \Delta w = f \left(\tilde{u} \right) -   f \left(w + \tilde{u} \right) , \\
				w(0) = u_0- \tilde{u}_0 ,
			\end{cases}
		\end{align}
		where $f \left(\tilde{v} \right) = -\left|\tilde{u} \right|^{4/(N-2)} \tilde{u} $. We can expand $E(v)$ as
		\begin{align}\label{eq:decompose}
			E(u) = E \left(\tilde{u} \right) + I + II,
		\end{align}
		where
		\begin{align}
			I =  \Re \int_{\Om} \left( \nabla \tilde{u} \cdot\overline{\nabla w }
			- f \left(\tilde{u} \right) \overline{w} \right) \,\mathrm{d}x ,
		\end{align}
		and
		\begin{align}
			II = \int_{\Om} \left( \f12 |\nabla w|^2 - \left( F \left(w+ \tilde{u} \right) - F \left(\tilde{u} \right) - \Re \left( f \left(\tilde{u} \right) \bar{w}
			\right) \right) \right)
			\,\mathrm{d}x ,
		\end{align}
		with $F \left(\tilde{u} \right) = \f{N-2}{2N} \left|\tilde{u} \right|^{\f{2N}{N-2}}$.
		
		Notice first that, by the Cauchy-Schwarz inequality,
		\begin{align}
			\f{d}{dt} \|w\|_{L^2}^2 =
			- 2  \Re \int_{\Om}  \left( \overline{ f \left(\tilde{u} + w \right) - f \left(\tilde{u} \right) } \right) \cdot (iw)  \,\mathrm{d}x
			\le C_1 \int_{\Om} \left(|w|^2 + F(w) \right) \,\mathrm{d}x .
		\end{align}
		By Gronwall's lemma, we obtain
		\begin{align}\label{est:w_l2norm}
			\|w (t)\|_{L^2}^2 \le e^{C_1t} \|w (0)\|_{L^2}^2 + C_2 \int_0^t e^{C_1(t-s)} \|\nabla w(s)\|_{L^2}^2\,\mathrm{d}s,
		\end{align}
		where $C_1 = C_1 \left(\tilde{v} \right)$ and $C_2$ is a constant.
		
		Then by a standard limitation argument and \eqref{est:w_l2norm}, we have
		\begin{align}\label{est:i_0-i_t}
			\begin{aligned}
				I(0) - I(t) &= - \int_0^t \f{d}{d s}  \left(\Re \int_{\Om} \left( \nabla \tilde{u} \cdot\overline{\nabla w }
				- f \left(\tilde{u} \right) \overline{w} \right)(s,x)  \,\mathrm{d}x \right)
				\,\mathrm{d}s \\
				& \le C_3 \int_0^t \int_{\Om}
				|w|^2 + F(w)  \,\mathrm{d}x \mathrm{d}s
				\le C_4 \|w (0)\|_{L^2}^2 + C_4 \int_0^t \|\nabla w(s)\|_{L^2}^2\,\mathrm{d}s
			\end{aligned}
		\end{align}
		with $C_4 = C_4 \left(T, \tilde{u}\right)$.
		
		Since $E \left(\tilde{u}(t) \right) = E \left(\tilde{u}_0 \right)$, and $E(u(t)) \le E(u_0)$, the expansion \eqref{eq:decompose} implies
		\begin{align}\label{est:v_t}
			0 \le E(u_0) - E(u(t)) = I(0) - I(t) + II(0) - II(t).
		\end{align}
		By the Cauchy-Schwarz inequality, Sobolev inequality and energy of ground state, we infer that
		\begin{align}\label{est:ii_t}
			\begin{aligned}
				II(t) &\ge \int_{\Om} \left[ \f12|\nabla w|^2 - C F(w) - C_5 |w|^2  \right]  \,\mathrm{d}x
				\ge \f12\|\nabla w(t)\|_{L^2_x}^2  - C_5 \|w(t)\|_{L^2}^2 ,
			\end{aligned}
		\end{align}
		where $C_5 = C_5 \left(\tilde{u} \right)$. At the same time, by the Cauchy-Schwarz inequality and Sobolev inequality, we have
		\begin{align}\label{est:ii_0}
			II(0) \le C_6  \left(\|\nabla w(0)\|_{L^2}^2 + \|w(0)\|_{L^2}^2 \right)
		\end{align}
		with $C_6 = C_6 \left(\tilde{u} \right)$
		
		We note that \eqref{est:ii_t} gives
		\begin{align}
			\|\nabla w(t)\|_{L^2}^2\le C_7  \left(II(t) + \|w(t)\|_{L^2}^2\right).
		\end{align}
		Then by \eqref{est:v_t}, \eqref{est:i_0-i_t} and \eqref{est:ii_0}, we finally conclude that
		\begin{align}\label{est:gronwall}
			\begin{aligned}
				\|\nabla w(t)\|_{L^2}^2  + \|w(t)\|_{L^2}^2 & \le C_9  \left(I(0) - I(t) + II(0) + \|w(t)\|_{L^2}^2  \right)\\
				& \le C \left(\|\nabla w(0)\|_{L^2}^2 + \|w(0)\|_{L^2}^2\right) + \int_0^t \|\nabla w(s)\|_{L^2}^2 + \|w(s)\|_{L^2}^2 \,\mathrm{d}s
			\end{aligned}
		\end{align}
		with a uniform constant $C= C \left(\tilde{v}, T\right)$.
	\end{proof}
	
	\begin{proof}[Proof of Corollary \ref{cor:nls}, part 2)]
This follows directly from Proposition \ref{thm:stability}.
   \end{proof}

\bigskip
%\noindent \textbf{Acknowledgments.} We are grateful to Rowan Killip, Monica Visan, Renjin Jiang and Lifeng Zhao for helpful discussion.

\end{document}